\documentclass[11pt]{elsarticle}

\usepackage{lineno,hyperref}
\modulolinenumbers[5]

\usepackage{calrsfs,amsmath,amssymb,graphicx,array,accents,figlatex,epsfig}
\usepackage{setspace,wrapfig,colortbl,marvosym, bbm, stmaryrd,centernot, color,upgreek,diagbox}
\usepackage{amsfonts}
\usepackage{epstopdf}
\usepackage{booktabs}
\usepackage{algorithm}
\usepackage{algpseudocode}
\usepackage{algorithm,algcompatible}
\usepackage{caption}

\usepackage[T1]{fontenc}
\usepackage{ae,aecompl}

\usepackage{fancyhdr}

\usepackage{etoolbox} 
\usepackage{natbib} 
\usepackage{blindtext}  
\usepackage{enumitem}

\ifpdf
  \DeclareGraphicsExtensions{.eps,.pdf,.png,.jpg}
\else
  \DeclareGraphicsExtensions{.eps}
\fi

\usepackage{enumitem}
\setlist[enumerate]{leftmargin=.5in}
\setlist[itemize]{leftmargin=.5in}

\newlength{\kaka}

\newcommand{\ahref}[2]{}

\newcommand{\obs}{^{\text{obs}}}

\newcommand{\beq}{\begin{equation}}
\newcommand{\eeq}{\end{equation}}
\newcommand{\lb}{\label}

\newcommand{\bea}{\begin{eqnarray}}
\newcommand{\eea}{\end{eqnarray}}
\newcommand{\bxr}{\begin{array}}
\newcommand{\exr}{\end{array}}

\newcommand\exs{\hspace*{0.4mm}}
\newcommand\xxs{\hspace*{0.2mm}}

\newcommand\nxs{\hspace*{-0.2mm}}

\newcommand{\norms}[1]{\parallel\! #1 \!\parallel}

\newcommand{\bC} {\boldsymbol{C}}

\newcommand{\bV} {\boldsymbol{V}}

\newcommand{\bn} {\boldsymbol{n}}

\newcommand{\bb} {\boldsymbol{b}}
\newcommand{\beps} {\boldsymbol{\epsilon}}

\newcommand{\bF} {\boldsymbol{F}}

\newcommand{\bx} {\boldsymbol{x}}
\newcommand{\by} {\boldsymbol{y}}
\newcommand{\be} {\boldsymbol{e}}

\newcommand{\bg} {{\boldsymbol{g}}}
\newcommand{\bq} {{\boldsymbol{q}}}

\newcommand{\sip} {\!\cdot\!}

\newcommand{\bzero}{\boldsymbol{0}}

\newcommand{\bu} {\boldsymbol{u}}

\newcommand{\bv} {\boldsymbol{v}}

\newcommand{\bxi} {\boldsymbol{\xi}}

\begin{document}

\begin{frontmatter}

\title{Deep regularization networks for inverse problems with noisy operators}

\author{Fatemeh Pourahmadian$^{1,2}$\corref{cor1}} 
\author{Yang Xu$^1$}

\address{\vspace*{1.25mm}$^1$ Department of Civil, Environmental \& Architectural Engineering, University of Colorado Boulder, USA} 
\address{\vspace*{-1mm}$^2$ Department of Applied Mathematics, University of Colorado Boulder, USA}
\cortext[cor1]{Corresponding author: tel. 303-492-2027, email {\tt fatemeh.pourahmadian@colorado.edu}}

\begin{abstract}
A supervised learning approach is proposed for regularization of large inverse problems where the main operator is built from noisy data. This is germane to superresolution imaging via the sampling indicators of the inverse scattering theory. We aim to accelerate the spatiotemporal regularization process for this class of inverse problems to enable real-time imaging. In this approach, a neural operator maps each pattern on the right-hand side of the scattering equation to its affiliated regularization parameter. The network is trained in two steps which entails:~(1)~training on low-resolution regularization maps furnished by the Morozov discrepancy principle with nonoptimal thresholds, and~(2)~optimizing network predictions through minimization of the Tikhonov loss function (i.e., the imaging objective) regulated by the validation loss. Step 2 allows for tailoring of the approximate maps of Step 1 toward construction of higher quality images. This approach enables direct learning from test data and dispenses with the need for a-priori knowledge of the optimal regularization maps (or labeled datasets) for training. The network, trained on low-resolution data, is then used to quickly generate dense regularization maps for high-resolution images. We highlight the importance of the training loss function on the network's generalizability. In particular, we demonstrate that networks informed by the logic of discrepancy principle lead to images of higher contrast. In this case, the training process involves multitask optimization with many loss components. We explore the state of the art in loss balancing and propose a new method to adaptively select the appropriate loss weights during training without requiring an additional optimization process. The proposed approach is synthetically examined for imaging damage evolution in an elastic plate. The results indicate that the \emph{discrepancy-informed regularization networks} not only accelerate the imaging process, but also remarkably enhance the image quality in complex environments.   

\end{abstract}


\end{frontmatter}

\section{Introduction} \label{sec1}

Recent studies in laser ultrasonic (LU) imaging~\cite{naru2023,pour2021} reveal that sampling methods of the inverse scattering theory~\cite{cakoni2022,monk2023,recoq2025}, thanks to their rigorous mathematical foundation and use of full-waveform data, enable  reconstruction of the subsurface with a quality that exceeds state-of-the-art LU imaging solutions. The latter~\cite{zare2024,budyn2021,cui2018} mostly relies on partial data inversion which could undermine the image quality, especially in complex environments, but enables fast reconstructions which is a much desired advantage in real-time imaging. In contrast, the high-quality reconstructions of the sampling methods, e.g., the linear sampling method (LSM), come at a cost of solving a typically large inverse problem featuring a noisy operator that requires spatiotemporal (or spatiospectral) regularization. The regularization process (in its current form) creates a bottleneck for real-time imaging. This is particularly the case when dense sampling is necessary for high-resolution reconstructions, or when specimens of complex geometry are excited by broadband inputs. In this case, optimal regularization maps are crucial for high-fidelity reconstructions, but their timely calculation remains a challenge for existing data inversion tools. This work aims to fundamentally address this issue through establishing an ML-based regularization process that reduces the computational cost of LSM reconstructions. This, in turn, enables real-time superresolution imaging of the subsurface which is a critical need, for instance, in online monitoring of additive manufacturing processes for quick detection and closed-loop mitigation of manufacturing defects in high-precision or safety-sensitive components~\cite{tahe2018}. 

In what follows, we outline the available approaches to learning regularization functions and argue that the existing logics are inadequate for addressing inverse problems with noisy operators, emphasizing the necessity of developing a new paradigm for this class of problems.   

Existing approaches to learning regularization maps~\cite{haber2003,habr2024,arri2019,afkh2021,li2020,haup2024,sant2024} are mostly predicated on the following form for the inverse problem
\beq\lb{can_form}
A \exs\bx_{\text{true}} \exs+\, \beps ~=~ \bb^\epsilon, 
\eeq
where the forward operator $A$ is \emph{exact} i.e., noiseless and the observation data $\bb^\epsilon$ on the right-hand side is perturbed by an unknown additive noise $\beps$. Given $A$, $\bb^\epsilon$, and the ill-posed nature of~\eqref{can_form}, one may construct an approximate solution $\tilde{\bx}$ to $\bx_{\text{true}}$ through variational regularization by minimizing for instance the Tikhonov loss function~\cite{tikh1963,Kress1999},   
\beq\lb{app_sol}
\tilde{\bx}(\bb^\epsilon,\lambda) ~=~ \text{arg}\!\min_{\bx}  \norms{\!\exs A \exs\bx - \bb^\epsilon \exs\!}_{L^2}^2 \exs+\,\,  \mathcal{R}(\bx,\lambda), \qquad  \mathcal{R}(\bx,\lambda) ~=~ \lambda\!\norms{\nxs\bx\nxs}_{L^2}^2,  
\eeq
where the regularization parameter $\lambda > 0$ must be specified prior to the calculation of $\tilde{\bx}$. Common approaches to estimating $\lambda$~\cite{morozov1993,hans2010,hans1993,xu2023} requires solving~\eqref{app_sol} (or a closely related form) multiple times in order to satisfy a certain criteria. This may be computationally expensive for large-scale or nonlinear problems~\cite{kalt2008, benn2018,gala1992,vogel1996,mead2008}. 

In their seminal work~\cite{haber2003}, Haber and Tenorio proposed a learning approach for optimal regularization functional $\mathcal{R}(\bx,\lambda)$ which has, over the past two decades, led to a suit of new techniques for supervised learning of regularization parameters and functions~\cite{habr2024,habe2008,chun2024,gouj2023,kofl2023,elia2023,li2020,antil2020,lunz2018,lunz2018,cala2017,chun2017,chun2011}. This is accomplished through solving a bilevel optimization problem which, within the context of~\eqref{app_sol}, takes the following form
\beq \lb{opt_lambda}
\!\left\{\begin{array}{l}
\!\! \lambda_{\text{opt}} ~=~ \text{arg}\!\min\limits_{\lambda} \dfrac{1}{N} \sum\limits_{i = 1}^{N} \norms{\tilde{\bx}^i - \bx_{\text{true}}^i}_{L^2}^2, \\*[5mm]
\!\! \tilde{\bx}^i(\bb^{i},\lambda) ~=~ \text{arg}\!\min\limits_{\bx}  \norms{\!\exs A \exs\bx - \bb^{i} \exs\!}_{L^2}^2 \exs+\,\, \lambda\!\norms{\nxs\bx^i\nxs}_{L^2}^2.
\end{array}\right. 
\eeq  

Here, the optimal regularization parameter $\lambda_{\text{opt}}$ is learned from a set of training pairs $(\bb^{i},\bx_{\text{true}}^i)$, $i = 1,\ldots, N$. This logic and most of its existing variants in the literature leverage the exact nature of operator $A$ to gain access to $\bx_{\text{true}}^i$, $i = 1,\ldots, N$, in order to construct the training dataset. In this setting, the training may be conducted \emph{offline} using synthetic data. The trained operator is then used \emph{online} to generate enhanced regularization maps for unseen data and/or to expedite the data inversion. 

In this study, we aim to develop a method for supervised learning of regularization maps for another class of inverse problems of general form 
\beq\lb{can_form2}
F^\delta \bg ~=~\bu_{\text{\tiny L}}, \qquad \norms{\nxs F^\delta - \exs F \nxs}_{L^2}\,\, \leqslant \, \delta.
\eeq

Here, $F^\delta$ is an integral operator whose kernel is constructed from experimental data and thus is contaminated by noise. $\delta\!>\!0$ is a measure of noise in data typically defined by the $L^2$ distance between the noiseless and noisy operators. However, the noiseless operator $F$ is fundamentally unknown, and $\delta$ is typically uncertain in the experiments. The right-hand side $\bu_{\text{\tiny L}}$ is assumed to be \emph{exact} as it is constructed by a computer model or a known closed-form solution. \eqref{can_form2} is ill-posed, and given $F^\delta$ and $\bu_{\text{\tiny L}}$, the objective is to build an approximate solution $\tilde{\bg}_{\text{\tiny L}}$ using variational regularization. In this work, we focus on Tikhonov-type loss functions, similar to~\eqref{app_sol}, which is germane to the classical formulation of the linear sampling method~\cite{cakoni2011linear,colton2003linear}. There are more rigorous forms for regularization of~\eqref{can_form2} for instance via the generalized linear sampling method~\cite{Audibert2014,Fatemeh2017,pour2020} which we reserve for future investigations.    

In this setting, the reconstruction involves solving~\eqref{can_form2}, by minimizing the Tikhonov loss, for a dictionary of right-hand sides $\bu^{\text{n}}_{\text{\tiny L}}$, ${\text{n}} = 1, \ldots, N_{\text{\tiny RHS}}$. The image is then constructed based on $1/\!\norms{\nxs \tilde{\bg}^{\text{n}}_{\text{\tiny L}} \nxs}$. Here, the regularization parameter is in fact a vector of length $N_{\text{\tiny RHS}}$ as it needs to be independently specified for every $\bu^{\text{n}}_{\text{\tiny L}}$. The issue is that there is no direct way to obtain the optimal distribution for the regularization parameter. As such, the simplest approach to solving~\eqref{can_form2}, by way of discrepancy principle~\cite{moro2012,kirsch2011,engl1996}, requires two numerical root finding operations per right-hand side which turns into a major obstacle in real-time imaging mainly due to: (1) length of the regularization vector, (2) complexity of the Picard plots, and (3) delicate process of hyperparameter tuning. Point 1 indicates that the regularization vector is typically large even in two-dimensional (2D) imaging. For example, the length of regularization vector in single-frequency 2D reconstructions of this study is of $O(10^6)$, while the corresponding vector for time-domain reconstructions in the same configuration may be of $O(10^9)$. These scales may sharply grow in higher-frequency or larger reconstructions and in 3D tomography. Point 2 refers to the largely overlapping distributions of (a) eigenvalues of $F^\delta$ and (b) the projected right-hand-side patterns in the Picard plots which is especially the case in complex environments with many scatterers of various dimensions and geometries. This makes the selection of spectral filters not only more difficult but also critical for high-fidelity reconstructions. In this context, Point 3 highlights the challenge of careful (yet efficient) tuning of the (Tikhonov) filter parameters since the problem is high-dimensional and the optimal values for the hyperparameters, such as the Morozov threshold, may reside in a tight range. Given the above, developing a suitable approach to learning the regularization process may be the key to real-time imaging by way of the sampling indicators.
 
The fundamental challenge that impedes direct application of the current methodologies, mentioned above, for learning regularization maps to~\eqref{can_form2} is the presence of noise in $F^\delta$. In this case, the exact solution $\bg_{\text{true}}$ to~\eqref{can_form2} which is one of the main elements of training in the existing logics remains unknown. In addition, since $F^\delta$ is constructed online from test data, calculation of optimal regularization maps for training may be too expensive in high-dimensional problems.   

To help bridge the gap, we propose a learning approach for Tikhonov regularization of~\eqref{can_form2} where the training data consists of (a) experimental data used to build $F^\delta$ and (b) a subset of right-hand-side patterns $\bu^{\text{t}}_{\text{\tiny L}}$, ${\text{t}} = 1, \ldots, N_{\text{\tiny trn}}$ where $N_{\text{\tiny trn}} <\!\!\nxs< N_{\text{\tiny RHS}}$. Using eigenvalues and eigenvectors of $F^\delta$ and the reduced set $\lbrace \bu^{\text{t}}_{\text{\tiny L}} \rbrace$, a low-resolution and approximate regularization map is quickly constructed based on the Morozov discrepancy principle. Here, we use a rough estimate for the Morozov threshold and do not require an a-priori knowledge of its optimal value (or its manual tuning). The idea is to map the regularization process by a neural network that takes each right-hand side pattern $\bu^{\text{t}}_{\text{\tiny L}}$, projected onto the eigenspace of $F^\delta$, to its affiliated regularization parameter. The latter is then used to compute the minimizer $\tilde{\bg}^{\text{t}}_{\text{\tiny L}}$ of the Tikhonov loss function. The training involves two steps: (1) learning the discrepancy principle in order to quickly generate dense regularization maps germane to all $N_{\text{\tiny RHS}}$ patterns on the right-hand side of~\eqref{can_form2}, and (2) optimizing the network predictions, beyond the approximate maps of step 1, to enhance image quality through controlled minimization of the Tikhonov loss (or the imaging objective). Step 1 is accomplished by minimizing a measure of distance between network predictions and the discrepancy-based regularization maps on low-resolution data. Special attention is paid to design of appropriate loss functions for training. It is shown that the network's generalizability improves in Step 1 when the logic of discrepancy principle in included in the loss function. In this case, we propose an efficient algorithm for loss balancing which builds on the logic of GradNorm~\cite{chen2018} and Dynamic Scaling~\cite{xu2024} but does not require an additional optimization process to compute the loss weights during training. Step 2 makes use of the model trained in step 1 as an initial state. Minimizing the Tikhonov loss in this step does not involve labeled datasets and thus may lead to overfitting. To regulate the training process, we introduce a criteria to stop training in Step 2 based on the relative trajectories of training and validation loss functions and show that this criteria is critical for enhancing the network output while preserving its stability. The performance of regularization networks is synthetically examined for imaging a dynamic damage zone in an elastic plate where a single microcrack evolves into a cloud of randomly distributed fractures over a sequence of time steps. This allows us to investigate the consistency of our findings across various configurations.   

A unique advantage of the proposed approach is that Step 2 enables direct learning from test data with a key objective of enhancing the image contrast. This is analyzed by introducing a contrast metric and studying its evolution by the end of each training step across configurations. The results demonstrate that the discrepancy-informed regularization networks can significantly improve the image quality and contrast when the noise in data is significant. Remarkable gain in compute time and efficiency may arise from the fact that: (a) the proposed approach dispenses with the need for manual tuning of hyperparameters such as the Morozov threshold, and (b) the regularization network is trained on limited data (i.e.,~downsampled right-hand-side patterns) and then used for quick evaluations on a much larger sampling grid. Network generalizability is paramount for the proper execution of the latter.  

This paper is organized as follows. Section~\ref{PS} provides a more comprehensive statement of the inverse problem. Section~\ref{Method} introduces the deep regularization networks, modes and steps of training, learning objectives as well as loss balancing and training regulation strategies. Section~\ref{IR} is dedicated to implementation of the proposed approach and discussion of the results.

\section{Problem statement}\label{PS}

Imaging by way of the sampling indicators~\cite{cakoni2022,Audibert2014,audi2017,Fatemeh2017} involves minimizing loss functions of the form
\beq\lb{GCf}
\mathfrak{J}(\bg; \bu_{\text{\tiny L}}) ~:=~ \norms{\nxs F^\delta \bg \exs-\exs\bu_{\text{\tiny L}} \nxs}_{L^2}^2 +~ \alpha\, \big(\bg, B \bg \big)_{L^2}, \qquad  \alpha>0, 
\eeq
with respect to $\bg \in L^2$. Here, the operator $B$ can be expressed as the following  
\beq\lb{B}
B \in \big{\lbrace} I, \, F^\delta_{\sharp}, \, F^\delta_{\sharp} +  \alpha^{-\chi} \delta I \big{\rbrace}, \quad F^\delta_{\sharp}\,\colon \!\!\!=\, \frac{1}{2} |F^\delta+F^{\delta^*}\nxs\nxs | \:+\: \frac{1}{2 \textrm{i}} (F^\delta\nxs-F^{\delta^*}\nxs), \quad \chi \in \, ]0, 1[,
\eeq
wherein $I$ is the identity operator and $()^*$ indicates the adjoint operator.

\paragraph*{Anatomy of $F^\delta$}
In experiments, the specimen is illuminated by a set of incident waves over the excitation surface $S^{\textrm{inc}}$, while the resulting motion is captured on the observation surface $S\obs$. Let $\bu^{\text{obs}}(\bxi, \by; \omega)$ denote the spectrum of displacement measured at the detector location $\bxi \in S\obs$ and frequency $\omega$ due to excitation at $\by \in S^{\textrm{inc}}$. In parallel, let $\bu^{\text{f}}(\bxi, \by; \omega)$ be the simulated response of the intact specimen (i.e., the background) in the same sensing configuration. In this setting, the scattering operator takes the form
\beq\lb{So} 
F^\delta(\bg)(\bxi;\omega) ~=\,  \int_{S^{\text{inc}\!}} \bV^\delta(\bxi,\by;\omega) \sip \bg(\by;\omega) \,\, \text{d}S_{\by}, \qquad \bg \in L^2, \quad\!\! \bxi \in S\obs,
\eeq 
where $V^\delta_{ij}(\bxi,\by;\omega)$, $i,j\!=\!1,2,3$, in~\eqref{So} is the $i^{\textrm{th}}$ component of scattered field $[u_i^{\text{obs}} \nxs- u_i^{\text{f}}\exs ](\bxi,\by;\omega)$ at $\bxi \in S\obs$ due to excitation at $\by \in S^{\textrm{inc}}$ in the $j^{\textrm{th}}$ direction. Here, we assume that $F^\delta$ is compact, injective, and has a dense range. Within the context of ultrasonic imaging, these properties are rigorously established in~\cite{Fatemeh2017} and~\cite{pour2020} for imaging from far- and near- field data respectively. See~\cite{cakoni2022} for the corresponding analysis in electromagnetic inverse scattering.  

\paragraph*{Imaging process}

In a computer model of the background, the image support is specified and sampled at $N_{\text{p}}$ points. Next, a set of trial scatterers $\lbrace \text{L}_{\text{s}} \rbrace$, $\text{s} = 1,\ldots, N_{\text{s}}$ is defined; each of which is then placed at every sampling point, one at a time, and the affiliated scattering signature is computed over the observation surface $S\obs$. The trial scatterers may include cracks and pores that induce dipole and monopole footprints on $S\obs$.  As such, one may create a library of simulated patterns ${\bu}_{\text{\tiny L}}^{\text{n}}$, $\text{n} = 1,2,\ldots, N_{\text{p}} N_{\text{s}}$ that include {(i)} monopole signatures created by planting a set of point sources with various polarizations at every sampling point, and {(ii)} dipole patterns constructed by nucleating infinitesimal fractures of varying unit normal vectors on the same sampling grid. Upon discretizing $F^\delta$, the reconstruction involves solving the highly {\emph{ill-posed}} scattering equation $\boldsymbol{F}^{\delta} {\bg}^{\text{n}} = {\bu}_{\text{\tiny L}}^{\text{n}}$ for all $\text{n} = 1,2,\ldots, N_{\text{p}} N_{\text{s}}$ through a regularization process. For this purpose, in this study, we consider the case of $B = I$ in~\eqref{GCf}, corresponding to the classical linear sampling method (LSM), and aim to minimize the resulting Tikhonov loss function
\vspace*{-0.5mm}
 \beq\label{dllsm}
\mathfrak{J}_{\text{LSM}}^{\text{n}}(\bg^{\text{n}}) \,\,\colon \!\!\!= \,\, \norms{\nxs \boldsymbol{F}^\delta \bg^{\text{n}} \exs-\exs\bu^{\text{n}}_{\text{\tiny L}} \nxs}_{L^2}^2 +\,\,  \alpha^{\text{n}} \!\nxs\norms{{\bg}^{\text{n}}}_{L^2}^2, \qquad \text{n} = 1,2,\ldots, N_{\text{p}} N_{\text{s}}, 
\vspace*{-0.5mm}
\eeq 
where $\alpha^{\text{n}}$ is the regularization parameter that needs to be computed separately for every right-hand side ${\bu}_{\text{\tiny L}}^{\text{n}}$. The fastest (and perhaps the simplest) approach to estimating $\alpha^{\text{n}}$ is the Morozov discrepancy principle~\cite{Kress1999} that minimizes $\mathfrak{J}_{\text{LSM}}^{\text{n}}$ for a fixed residual threshold $\norms{\nxs \boldsymbol{F}^\delta \bg^{\text{n}} \exs-\exs\bu^{\text{n}}_{\text{\tiny L}} \nxs} \,=\nxs \eta \!\nxs\norms{\!{\bg}^{\text{n}}\!\nxs}$. The challenge is that the optimal value for $\eta\!>\!0$ is unknown and this parameter is typically manually tuned. For a given $\eta$, finding $\alpha^{\text{n}}$ requires one numerical root-finding operation for every ${\text{n}}$. This could make manual tuning of $\eta$ quite computationally expensive for large sampling grids, or more precisely, when $N_{\text{p}} N_{\text{s}}$ is large.  

Once the minimizers $\bg^{\text{n}}$, $\text{n} = 1,2,\ldots, N_{\text{p}} N_{\text{s}}$, are computed for every trial scatterer $\text{L}_{\text{s}}$ and sampling point $\bx^{\text{p}}$ i.e., $\bg^{\text{n}} = \bg^{\text{n}}(\bx^{\text{p}}, \text{L}_{\text{s}})$, $\text{p} = 1,\ldots, N_{\text{p}}$ and ${\text{s}} = 1,\ldots, N_{\text{s}}$. The LSM imaging indicator is formed as the following 
\beq\lb{LSM}
\mathcal{L}(\bx^{\text{p}}) \,\, := \,\, \frac{1}{\norms{\bg^{\text{p}}}^2}, \qquad
\textcolor{black}{
\bg^{\exs\text{p}}(\bx^{\text{p}}) \,\,\colon \!\!\!= \,\,\, \text{arg\hspace{-11.5mm}}\min_{\lbrace \bg^{\text{n}}(\bx^{\text{p}}, \text{L}_{\text{s}}) \rbrace_{s = 1,\ldots, N_{\text{s}}}}\!\! \norms{\bg^{\text{n}}(\bx^{\text{p}}, \text{L}_{\text{s}})}^2_{L^2}, ~\text{p} = 1,\ldots, N_{\text{p}}.}
\eeq

In this work, we aim to accelerate calculation of the imaging indicator $\mathcal{L}(\bx^{\text{p}})$ by automating and optimizing the regularization process.

\section{Deep regularization networks for inversion of noisy operators}\label{Method}

\begin{figure}[!tp]
\vspace*{-0mm}
\center\includegraphics[width=0.97\linewidth]{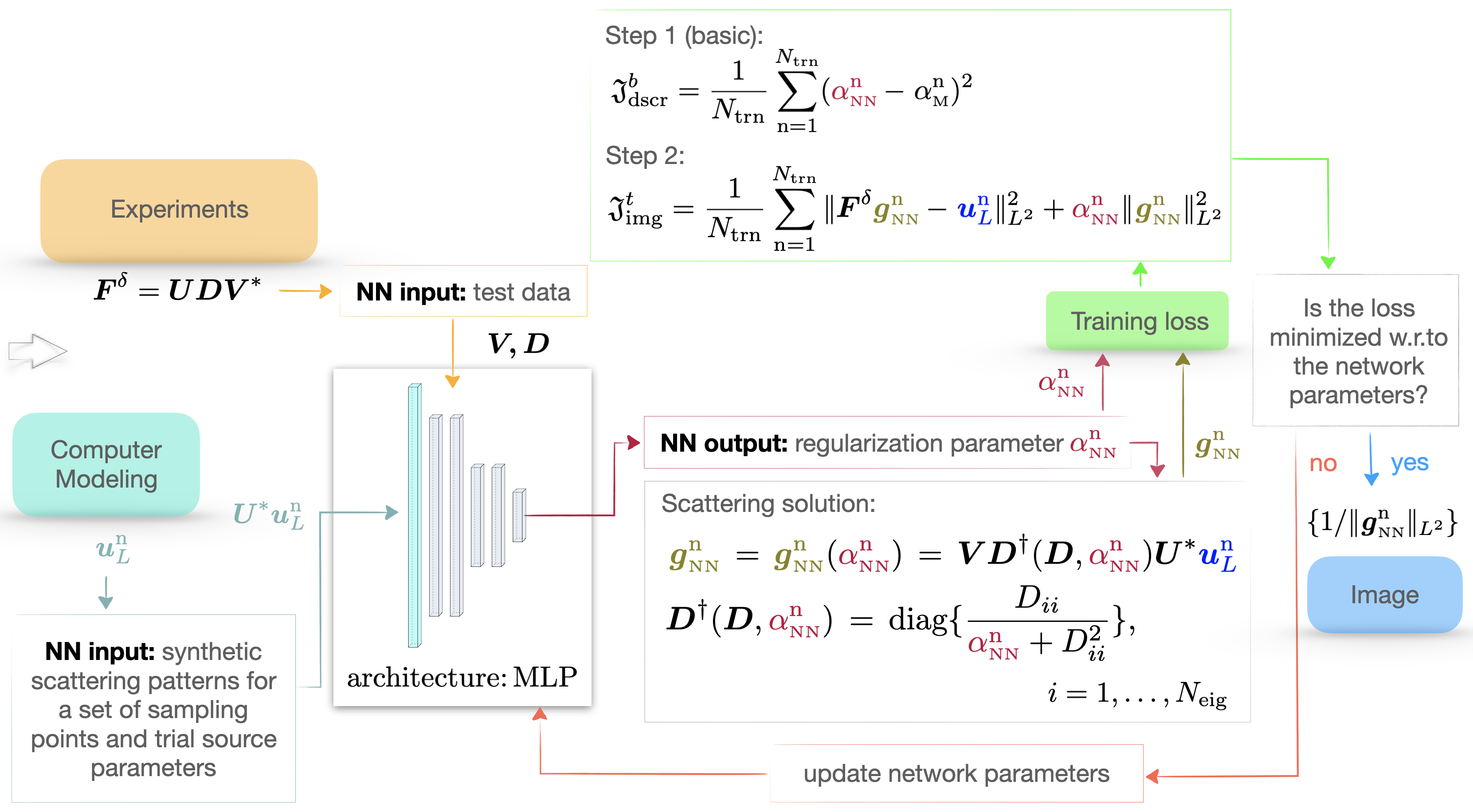} \vspace*{-2mm}
\caption{\small{The proposed learning logic for the Regularization Network. The network {(i)} takes the SVD of scattering operator $\boldsymbol{F}^\delta = \boldsymbol{U} \boldsymbol{D} \boldsymbol{V}^*$, {(ii)} projects the right-hand side $\bu_{\text{\tiny L}}^\text{n}$ onto the eigenspace of measurements $\boldsymbol{U}$, and {(iii)} outputs the affiliated regularization parameters $\alpha_{\text{\tiny NN}}^\text{n}$. The output is then used to compute the spectral filter factors ${D}^\dagger_{ii}$, $i = 1,2,\ldots,N_{\text{eig}}$, germane to Tikhonov regularization and the scattering solution $\boldsymbol{g}_{\text{\tiny NN}}^\text{n}$. The learning objective is two-fold: (a) to learn the Morozov discrepancy principle in Step 1 to accelerate the regularization process, and (b) to optimize the network predictions in Step 2 to further enhance the quality of reconstructions.}}
\label{DLt2}
\vspace*{-5.5mm}
\end{figure}

The data used for training and evaluation of the regularization network (R-Net) includes {(i)} the scattering operator $\boldsymbol{F}^\delta$, and {(ii)} the library of patterns $\lbrace \bu^{\text{n}}_{\text{\tiny L}} \rbrace$, $\text{n} = 1,2,\ldots, N_{\text{p}} N_{\text{s}}$, on the right-hand side of~\eqref{dllsm} which henceforth is referred to as the RHS patterns. In this study, R-Net takes the form of a multilayer perceptron (MLP) owing to its simplicity and universal approximation property~\cite{horn1989}. The network makes use of the singular value decomposition of $\boldsymbol{F}^\delta = \boldsymbol{U} \boldsymbol{D} \boldsymbol{V}^*$ to create the appropriate learning space for Tikhonov regularization as shown in Fig.~\ref{DLt2}. More specifically, $\boldsymbol{U}$ is used to project the RHS pattern $\bu^{\text{n}}_{\text{\tiny L}}$ onto the eigenspace of measurements. This modal representation of input $\boldsymbol{U}^*\bu^{\text{n}}_{\text{\tiny L}}$ along with the eigenvalues in $\boldsymbol{D}$ furnish the (Picard) design space for the spectral filter factors~\cite{hans2010,chun2017} in $\boldsymbol{D}^\dagger = \boldsymbol{D}^\dagger(\boldsymbol{D},\alpha_{\text{\tiny{NN}}}^{\text{n}})$ that will be used along with $\boldsymbol{V}$ to form the regularized solution $\bg_{\text{\tiny{NN}}}^{\text{n}}$. In this setting, R-Net is a map from $\boldsymbol{U}^*\bu^{\text{n}}_{\text{\tiny L}}$ to $\alpha_{\text{\tiny{NN}}}^{\text{n}}$ which is trained in two steps. The objective of Step 1 is to learn the Morozov discrepancy principle in order to accelerate the regularization process, while Step 2 seeks to further optimize the network predictions (beyond the discrepancy principle) in search for a min-norm solution which, within the context of sampling indicators, may be interpreted as an effort to enhance the image contrast. 

It is worth mentioning that the proposed regularization process in Fig.~\ref{DLt2} is highly compatible with the logic of kernel methods in operator learning in that the scattering solution $\boldsymbol{g}_{\text{\tiny NN}}^\text{n}$ may be represented by an integral operation whose kernel contains the unknown filter factors ${D}^\dagger_{ii}$, $i = 1,2,\ldots,N_{\text{eig}}$. In particular, the recently developed implicit Fourier neural operators (IFNO)~\cite{li2020(2),kova2021,you2022} may be advantageous over the classical MLP due to their resolution independence and superior extrapolation capability which could be a subject of future investigation.

\subsection{R-Net training:~Step 1}\label{TS1}

This step makes use of low-resolution regularization maps based on the discrepancy principle for training. For this purpose, the image support is first downsampled by a factor of $m$ in every direction so that the number of RHS patterns is reduced to $N_{\text{r}} = N_{\text{p}} N_{\text{s}}/m^d$ with $d = 2,3$ indicating the image dimensions. Next, on adopting an estimate for the Morozov threshold $\eta_\circ\!>\!0$, the regularization map $\lbrace \alpha_{\text{\tiny{M}}}^{\text{r}} \rbrace$ and associated scattering solutions $\lbrace \boldsymbol{g}_{\text{\tiny{M}}}^{\text{r}} \rbrace$ in~\eqref{dllsm} are computed such that $\norms{\nxs F^\delta \bg_{\text{\tiny{M}}}^{\text{r}} \exs-\exs \bu^{\text{r}}_{\text{\tiny L}} \nxs} \,=\nxs \eta_\circ \!\nxs\norms{\nxs{\bg}_{\text{\tiny{M}}}^{\text{r}}\!}$ on the reduced grid ${\text{r}} = 1,2,\ldots, N_{\text{r}}$. Regarding the choice of $\eta_\circ$, we learned from the synthetic experiments of Section~\ref{IR} that an overestimation is generally safe and results in an approximate regularization map that will be optimized in Step 2, whereas a significantly underestimated $\eta_\circ$ may lead to unstable reconstructions in both training steps. In light of this and in absence of an optimal estimate, one may set $\eta_\circ \in [0.3 \,\,\, 0.5)$ to generate the data for training in Step 1. At this point, one may further reduce the size of training dataset by selecting one trial scatterer $\text{L}_\text{s}$, $\text{s} = 1,2,\ldots, N_{\text{s}}$, per sampling point as the following   
\beq\lb{Lp}
\text{L}_{\exs\text{p}}(\bx^{\text{p}}) \,\,\colon \!\!\!= \,\, \text{arg\hspace{-6mm}}\min_{\lbrace \text{L}_{\text{s}} \rbrace_{s = 1,\ldots, N_{\text{s}}}}\hspace{-6mm} \norms{\boldsymbol{g}_{\text{\tiny{M}}}^{\text{r}}(\bx^{\text{p}}, \text{L}_{\text{s}})}^2_{L^2},  \quad {\text{r}} = {\text{r}}({\text{p}},{\text{s}}), \quad \text{p} = 1,\ldots, N_{\text{p}}/m^d.
\eeq
In this setting, a small subset of RHS patterns $\lbrace \bu^{\text{t}}_{\text{\tiny L}} \rbrace$, $\text{t} = 1,2,\ldots, N_{\text{trn}} \leqslant N_{\text{p}}/m^d$, participate in training. The trained R-Net will then be used to generate the regularization map for the original library of patterns $\lbrace \bu^{\text{n}}_{\text{\tiny L}} \rbrace$, $\text{n} = 1,2,\ldots, N_{\text{p}} N_{\text{s}}$, by a single forward pass through the network. 
Given $\boldsymbol{F}^\delta = \boldsymbol{U} \boldsymbol{D} \boldsymbol{V}^*$, $\lbrace \bu^{\text{t}}_{\text{\tiny L}} \rbrace$ and $\lbrace \alpha_{\text{\tiny{M}}}^{\text{t}} \rbrace$ for $\text{t} = 1,2,\ldots, N_{\text{trn}}$, training R-Net in Step 1 may be conducted in two modes, namely: basic and informed. 

\subsubsection{Basic mode of training}\label{BMT}

In the basic mode, network parameters are optimized such that the misfit between $\lbrace \alpha_{\text{\tiny{M}}}^{\text{t}} \rbrace$ and the network-predicted regularization parameters $\lbrace \alpha_{\text{\tiny{NN}}}^{\text{t}} \rbrace$ is minimized for $\text{t} = 1,2,\ldots, N_{\text{trn}}$. As such, the training loss function takes the form   
\beq\lb{Jbds}  
\mathfrak{J}^b_{\text{dscr}} \,=\, \frac{1}{{N_{\text{trn}}}}\sum_{\text{t} = 1}^{N_{\text{trn}}} (\textcolor{black}{\alpha_{\text{\tiny{NN}}}^{\text{t}}} \hspace{-0.2mm} - \hspace{0.2mm} \textcolor{black}{\alpha_{\text{\tiny{M}}}^{\text{t}}})^2. 
\eeq

In this case, on setting $\lbrace {\alpha_{\text{\tiny{M}}}^{\text{v}}} \rbrace_{\text{v} = 1,2,\ldots, N_{\text{vld}}} \subset  \lbrace {\alpha_{\text{\tiny{M}}}^{\text{n}}} \rbrace_{\text{n} = 1,2,\ldots, N_{\text{p}}N_{\text{s}}} \!\setminus \lbrace {\alpha_{\text{\tiny{M}}}^{\text{t}}} \rbrace_{\text{t} = 1,2,\ldots, N_{\text{trn}}}$, the validation loss may be computed as follows,  
\beq\lb{Vbds}  
\mathfrak{V}^b_{\text{dscr}} \,=\, \frac{1}{{N_{\text{vld}}}}\sum_{\text{n} = 1}^{N_{\text{vld}}} (\textcolor{black}{\alpha_{\text{\tiny{NN}}}^{\text{v}}} \hspace{-0.2mm} - \hspace{0.2mm} \textcolor{black}{\alpha_{\text{\tiny{M}}}^{\text{v}}})^2, 
\eeq
wherein $\lbrace{\alpha_{\text{\tiny{NN}}}^{\text{v}}\rbrace}_{\text{v} = 1,2,\ldots, N_{\text{vld}}}$ is the network output when the input $$\lbrace {\boldsymbol{U}^*\bu_{\text{\tiny{L}}}^{\text{v}}} \rbrace_{\text{v} = 1,2,\ldots, N_{\text{vld}}} \subset \, \lbrace {\boldsymbol{U}^*\bu_{\text{\tiny{L}}}^{\text{n}}} \rbrace_{\text{n} = 1,2,\ldots, N_{\text{p}}N_{\text{s}}} \setminus \lbrace {\boldsymbol{U}^*\bu_{\text{\tiny{L}}}^{\text{t}}} \rbrace_{\text{t} = 1,2,\ldots, N_{\text{trn}}}$$ is affiliated with the validation dataset. The numerical experiments in Section~\ref{IR} indicate that while the basic mode of training is remarkably fast and leads to reconstructions comparable to that obtained via the discrepancy principle, its generalization capability is limited, and thus, image enhancement by basic R-Nets may be difficult. 

\subsubsection{Discrepancy-informed regularization networks}\label{iMT}

To help improve R-Net's generalizability, we propose the informed mode of training. Here, the idea is to include the logic of Morozov discrepancy principle~\cite{moro2012,kirsch2011,engl1996} in the learning process. For this purpose, the training loss function is recast as 
\beq\lb{Jids}  
\mathfrak{J}^i_{\text{dscr}} = \frac{1}{{N_{\text{trn}}}}\sum_{\text{t} = 1}^{N_{\text{trn}}} w_1^{\text{t}} J_1^{\xxs\text{t}} + w_2^{\text{t}} J_2^{\xxs\text{t}}
\eeq
where $J_1^{\xxs\text{t}}$, $\text{t} = 1,2,\ldots, N_{\text{trn}}$, is the square of normal misfit between network-predicted and discrepancy-based regularization parameters, i.e.,
\beq\lb{J1n}  
J_1^{\xxs\text{t}} =  \Big| \dfrac{\alpha_{\text{\tiny{NN}}}^{\text{t}}-\alpha_{\text{\tiny{M}}}^{\text{t}}}{\alpha_{\text{\tiny{N}}}^{\text{t}}} \Big|^2, \quad   \alpha_{\text{\tiny{N}}}^{\text{t}} \,=\, \alpha_{\text{\tiny{M}}}^{\text{t}} \,\,\, \text{or}\,\, \max(\alpha_{\text{\tiny{M}}}^{\text{t}}), \quad \forall {\text{t}} | \,\, \alpha_{\text{\tiny{N}}}^{\text{t}} > 0, 
\eeq
while $J_2^{\xxs\text{t}}$, $\text{t} = 1,2,\ldots, N_{\text{trn}}$, is the square of discrepancy functional with a-priori estimate $\eta_\circ$ 
\beq\lb{J2n}  
 J_2^{\xxs\text{t}}  = \big( \sum_{j = 1}^{N_{\text{eig}}} f_{j}^{\xxs\text{t}}(\alpha_{\text{\tiny{NN}}}^{\text{t}}, \eta_\circ) \big)^2, \quad f_{j}^{\xxs\text{t}}(\alpha_{\text{\tiny{NN}}}^{\text{t}}, \eta_\circ) = \frac{\alpha_{\text{\tiny{NN}}}^{\text{t}} - \eta_\circ^2 D_{jj}^2}{(\alpha_{\text{\tiny{NN}}}^{\text{t}} + D_{jj}^2)^2} \big{|} (\boldsymbol{u}^*_j, \textcolor{black}{\boldsymbol{u}^{\text{t}}_{\text{\tiny L}}}) \big{|}^2,
\eeq
wherein $\boldsymbol{u}_j$ is the $j^{\text{th}}$ column of $\boldsymbol{U}$ i.e., the $j^{\text{th}}$ left eigenvector of $\boldsymbol{F}^\delta$. Keep in mind that $J_2^{\xxs\text{t}} = 0$ provides a map from $\eta_\circ$ to $\alpha_{\text{\tiny{NN}}}^{\text{t}}$ for every RHS pattern $\boldsymbol{u}^{\text{t}}_{\text{\tiny L}}$. It is critical to properly balance the loss components during training by adaptively adjusting $w_1^{\text{t}}$ and $w_2^{\text{t}}$ for every $\text{t} = 1,2,\ldots, N_{\text{trn}}$. In this vein, one may show that $\forall\text{t}$, $  J_1^{\xxs\text{t}}$ and $J_2^{\xxs\text{t}}$ are of $O(1)$. Nonetheless, setting $w_1^{\text{t}} = w_2^{\text{t}} = 1$ did not result in balanced trainings in our numerical experiments. In most cases, the optimizer would prioritize one of the loss components and lose sensitivity to the other during training. To address this issue, we considered the state-of-the-art solutions for adaptive loss balancing~\cite{heyd2019,Kend2017,chen2018,wang2021,bisc2021,xu2024}, namely:~softmax adaptive weights (SoftAdapt)~\cite{heyd2019,Kend2017},~gradient normalization (GradNorm)~\cite{chen2018},~learning rate annealing (LRA)~\cite{wang2021}, relative loss balancing with random lookback (ReLoBRaLo)~\cite{bisc2021}, and dynamic scaling (DynScl)~\cite{xu2024}. In what follows, a brief discerption of each method is  included along with our findings on their performance when applied to~\eqref{Jids}. 

SoftAdapt weights each loss component based on a statistical measure of its rate of change relative to other objectives such that the convergence of all loss components in the parameter space is approximately isotropic~\cite{heyd2019}. On denoting the iteration step by $t$, let us define the rate of change of loss components by 
\begin{equation}
\label{sadapt}
\Delta J_k^{\xxs\text{t}}~=~J_k^{\xxs\text{t}}(t) \,-\, J_k^{\xxs\text{t}}(t-1), \quad k ~=~ 1,2, \quad \text{t} ~=~ 1,2,\ldots, N_{\text{trn}}, \quad t ~=~ 1,\ldots, N_{\text{{epoch}}},
\end{equation}
where $N_{\text{{epoch}}}$ signifies the number of epochs. In this setting, the loss weights are specified by
\begin{equation}
\label{eq:softadapt}
w_k^{\text{t}}~=~\displaystyle \frac{\text{e}^{\beta \exs \left(\Delta J_k^{\xxs\text{t}}-\max\left(\{\Delta J_k^{\xxs\text{t}}\}_{k = 1,2}\right)\right)}}{\sum_{k=1}^2 \text{e}^{\beta  \left(\Delta J_k^{\xxs\text{t}}-\max\left(\{\Delta J_k^{\xxs\text{t}}\}_{k = 1,2}\right)\right)}}, \quad \text{t} ~=~ 1,2,\ldots, N_{\text{trn}},
\end{equation}
where $\beta$ is a tunable hyperparameter whose default value is $\beta = 0.1$. Softmax adaptive weights are advantageous as they are computationally efficient and flexible, and thus, may be applied across a broad range of network architectures and loss functions, particularly when all components are on the same scale in a normalized space. However, there are three caveats when using this approach for balancing complex multi-objective systems: (i) the heuristic nature of the Softmax function, which may not effectively balance loss functions with numerous components due to its limited range, (ii) the highly oscillatory behavior of SoftAdapt weights especially in earlier epochs due to errors and irregular variations of loss components; which may interfere with the optimization process and lead to instability and lack of convergence, and (iii) tendency of the Softmax function to assign higher weights to components that exhibit slower rates of change. Keep in mind that the decay rates are equally influenced by the scale of each loss component at every epoch so that a fast-converging task could have a much smaller decay rate compared to a slower-converging objective of greater magnitude. This may mislead the optimizer into prioritizing objectives that are already converging, thereby reducing its sensitivity to other loss components that still require attention.

Our implementations of SoftAdapt for balancing~\eqref{Jids} were consistently unsuccessful mainly due to the first and second caveats mentioned above. In particular, the large number of objectives $2N_{\text{trn}}$, which is of $O(10^3)$ in the numerical experiments of Section~\ref{IR}, and the unstable behavior of SoftAdapt weights repeatedly led to the failure of optimization process.

LRA balances the objectives through modifying the learning rates on the basis of the gradients of various loss terms~\cite{wang2021}. LRA does not require a separate optimization for loss balancing which is a desired attribute for R-Nets. However, LRA weights are not bounded and may straddle across several scales during training which may result in repeated overshooting of various objectives or the optimizer might prioritize easier-to-achieve objectives leading to training instability. This issue is a critical risk in training the regularization networks since, here, we are concerned with the regularization of noisy operators, and LRA relies on the estimated derivatives by automatic differentiation to determine the loss weights at every epoch. Given the above, LRA may not be an appropriate choice for balancing~\eqref{Jids}.    

GradNorm balances the gradient of weighted objectives with respect to network parameters to ensure that all loss components train at similar rates~\cite{chen2018}. This method tackles the common problem of gradient imbalances by applying penalties to tasks with overly large or small gradients. Let $\text{\bf w}$ denote (a subset of) R-Net's weights. GradNorm finds the optimal weights $w_k^{\text{t}}(t)$, at every epoch $t$, by minimizing the $L_1$ distance between the actual and average norms of loss gradients with respect to $\text{\bf w}$. More specifically,
\begin{equation}
\label{eq:gradnorm}
\begin{aligned}
&\displaystyle L_{\nabla}~=\,\,\displaystyle \sum_{{\text{t}}=1}^{N_{\text{trn}}}\sum_{k=1}^{2}\left|\left\|\nabla_{\text{\bf w}} w_k^{\text{t}}(t) J_k^{\xxs\text{t}}(t) \right\Vert_2-\,\frac{1}{{2N_{\text{trn}}}}\sum_{{\text{t}}=1}^{N_{\text{trn}}}\sum_{k=1}^{2} \left(r_k^{\text{t}}(t)\right)^{\tilde{\upeta}} \left\|\nabla_{\text{\bf w}} w_k^{\text{t}}(t) J_k^{\xxs\text{t}}(t) \right\|_2 \right|_1, \\*[0.5mm]
&r_k^{\text{t}}(t)~=~\frac{J_k^{\xxs\text{t}}(t)/J_k^{\xxs\text{t}}(t_\circ)}{ \sum_{{\text{t}}=1}^{N_{\text{trn}}}\sum_{k=1}^{2}{J_k^{\xxs\text{t}}(t)/J_k^{\xxs\text{t}}(t_\circ)}},  \quad k ~=~ 1,2, \quad \text{t} ~=~ 1,2,\ldots, N_{\text{trn}},
\end{aligned}
\end{equation}
where $\tilde{\upeta}$ is a hyperparameter gauging the intensity of imposing a common convergence rate on all objectives. When loss components are of similar scale, a smaller $\tilde{\upeta}$ would be appropriate. 

GradNorm is logically sound but entails (1) automatic differentiation in order to evaluate $L_{\nabla}$ at every epoch, and (2) a separate minimization of $L_{\nabla}$ at every epoch while optimizing $\mathfrak{J}^i_{\text{dscr}}$. This when combined with noisy data and the many-objective loss in~\eqref{Jids} leads to an unstable and tardy training process which is inconsistent with the philosophy of R-Nets. Our implementations of GradNorm confirmed that the large vector of unknown weights at every epoch creates a quite complex inner optimization problem that ultimately disrupts minimization of the main objective $\mathfrak{J}^i_{\text{dscr}}$.     

ReLoBRaLo integrates the strengths of SoftAdapt and GradNorm to weight the loss components according to their past behavior and contributions~\cite{bisc2021}. Similar to SoftAdapt, LRA, and GradNorm, ReLoBRaLo involves hyperparameter tuning that in addition to the above-mentioned impediments further complicates the training process.

DynScl dynamically balances the loss by adjusting the scale of each objective during training. This is achieved by incorporating the specific logic (or physics) of each loss component, rather than relying on generic functions such as Softmax which may exhibit inconsistent behavior~\cite{xu2024}. In DynScl, the proper weight for each objective is identified relatively quickly and remains stable throughout the optimization process. In addition, In contrast to LRA and ReLoBRaLo, DynScl does not require a predefined lookback window, separate learning rate adjustment, or extensive hyperparameter tuning.

Inspired by the logic of GradNorm and advantages of DynScl, we explicitly analyzed the derivatives of $J_1^{\xxs\text{t}}$ and $J_2^{\xxs\text{t}}$ with respect to R-Net's parameters (denoted by ${\text{\bf w}}$ in~\eqref{eq:gradnorm}) and realized that by setting 
\beq\lb{w12n} 
\begin{aligned}
&w_1^{\text{t}} = \alpha_{\text{\tiny{N}}}^{\text{t}}, \quad \alpha_{\text{\tiny{N}}}^{\text{t}} \,=\, \alpha_{\text{\tiny{M}}}^{\text{t}} \,\,\, \text{or}\,\, \max(\alpha_{\text{\tiny{M}}}^{\text{t}}), \quad  \text{t} ~=~ 1,2,\ldots, N_{\text{trn}},   \\*[0.5mm] 
&w_2^{\text{t}} = \Big({\sum_{j = 1}^{N_{\text{eig}}} \frac{\partial f_{j}^{\xxs\text{t}}}{\partial \alpha_{\text{\tiny{NN}}}^{\text{t}}}(\alpha_{\text{\tiny{NN}}}^{\text{t}}, \eta_\circ) +\epsilon}\Big)^{-1}, \,\, \dfrac{\partial f_{j}^{\xxs\text{t}}}{\partial \alpha_{\text{\tiny{NN}}}^{\text{t}}}(\alpha_{\text{\tiny{NN}}}^{\text{t}}, \eta_\circ) \,=\, \frac{2D_{jj}^2(\alpha_{\text{\tiny{NN}}}^{\text{t}} + \eta_\circ^2)}{(\alpha_{\text{\tiny{NN}}}^{\text{t}} + D_{jj}^2)^3} \big{|} (\boldsymbol{u}^*_j, \textcolor{black}{\boldsymbol{u}^{\text{t}}_{\text{\tiny L}}}) \big{|}^2,
\end{aligned}
\eeq 
one may show that
\beq\lb{DJ12n} 
\begin{aligned}
& \nabla_{\!\text{w}} [w_1^{\text{t}} J_1^{\xxs\text{t}}](t) \,=\, \frac{2w_1^{\text{t}}}{\alpha_{\text{\tiny{N}}}^{\text{t}}} \, \Big( \dfrac{\alpha_{\text{\tiny{NN}}}^{\text{t}}-\alpha_{\text{\tiny{M}}}^{\text{t}}}{\alpha_{\text{\tiny{N}}}^{\text{t}}} \Big) \, \frac{\partial \alpha_{\text{\tiny{NN}}}^{\text{t}}}{\partial \text{w}} \,=\, O\big( 2 \frac{\partial \alpha_{\text{\tiny{NN}}}^{\text{t}}}{\partial \text{w}}\big),  \\*[0.5mm] 
&  \nabla_{\!\text{w}} [w_2^{\text{t}} J_2^{\xxs\text{t}}](t)  \,=\,  {2w_2^{\text{t}}} \, \sum_{j = 1}^{N_{\text{eig}}} f_{j}^{\xxs\text{t}}(\alpha_{\text{\tiny{NN}}}^{\text{t}}, \eta_\circ)  \, \Big({\sum_{j = 1}^{N_{\text{eig}}} \frac{\partial f_{j}^{\xxs\text{t}}}{\partial \alpha_{\text{\tiny{NN}}}^{\text{t}}}(\alpha_{\text{\tiny{NN}}}^{\text{t}}, \eta_\circ) }\Big) \, \frac{\partial \alpha_{\text{\tiny{NN}}}^{\text{t}}}{\partial \text{w}} \,=\, O\big( 2 \frac{\partial \alpha_{\text{\tiny{NN}}}^{\text{t}}}{\partial \text{w}}\big).
\end{aligned}
\eeq 

The discrepancy function $J_2^{\xxs\text{t}}=0$ is monotonically increasing with respect to $\alpha_{\text{\tiny{NN}}}^{\text{t}}$. Nonetheless, an infinitesimal threshold $\epsilon$ is used to avoid dividing by zero in the second of~\eqref{w12n} due to potential numerical errors. This approach automatically achieves the GradNorm's objective of normalizing the derivatives without requiring a separate minimization procedure per epoch or use of automatic differentiation. This is accomplished by explicitly approximating the scale of relevant derivatives instead of their numerical calculation which may involve significant computational cost and complications during the optimization process. In the synthetic experiments of Section~\ref{IR}, the proposed weights in~\eqref{w12n} seem to effectively balance the loss components across all cases considered. To monitor the training performance and gauge the number of training epochs, one may form the validation dataset $\lbrace {{\boldsymbol{U}^*\bu_{\text{\tiny{L}}}^{\text{v}}}, \alpha_{\text{\tiny{M}}}^{\text{v}}}, \alpha_{\text{\tiny{NN}}}^{\text{v}} \rbrace_{\text{v} = 1,2,\ldots, N_{\text{vld}}}$, similar to data used in~\eqref{Vbds}, to compute the validation loss as follows
\beq\lb{Vids}  
\begin{aligned}
&\mathfrak{V}^i_{\text{dscr}} = \frac{1}{{N_{\text{vld}}}}\sum_{\text{v} = 1}^{N_{\text{vld}}} w_1^{\text{v}} V_1^{\xxs\text{v}} + w_2^{\text{v}} V_2^{\xxs\text{v}}, \\*[0.5mm]
&V_1^{\xxs\text{v}} =  \Big| \dfrac{\alpha_{\text{\tiny{NN}}}^{\text{v}}-\alpha_{\text{\tiny{M}}}^{\text{v}}}{\alpha_{\text{\tiny{N}}}^{\text{v}}} \Big|^2, \quad   \alpha_{\text{\tiny{N}}}^{\text{v}} \,=\, \alpha_{\text{\tiny{M}}}^{\text{v}} \,\,\, \text{or}\,\, \max(\alpha_{\text{\tiny{M}}}^{\text{v}}), \quad \forall {\text{v}} | \,\, \alpha_{\text{\tiny{N}}}^{\text{v}} > 0, \\*[0.5mm]
&  V_2^{\xxs\text{v}}  = \big( \sum_{j = 1}^{N_{\text{eig}}} f_{j}^{\xxs\text{v}}(\alpha_{\text{\tiny{NN}}}^{\text{v}}, \eta_\circ) \big)^2, \quad f_{j}^{\xxs\text{v}}(\alpha_{\text{\tiny{NN}}}^{\text{v}}, \eta_\circ) = \frac{\alpha_{\text{\tiny{NN}}}^{\text{v}} - \eta_\circ^2 D_{jj}^2}{(\alpha_{\text{\tiny{NN}}}^{\text{v}} + D_{jj}^2)^2} \big{|} (\boldsymbol{u}^*_j, \textcolor{black}{\boldsymbol{u}^{\text{v}}_{\text{\tiny L}}}) \big{|}^2, \\*[0.5mm] 
&w_1^{\text{v}} = \alpha_{\text{\tiny{N}}}^{\text{v}}, \quad \alpha_{\text{\tiny{N}}}^{\text{v}} \,=\, \alpha_{\text{\tiny{M}}}^{\text{v}} \,\,\, \text{or}\,\, \max(\alpha_{\text{\tiny{M}}}^{\text{v}}), \quad  \text{v} ~=~ 1,2,\ldots, N_{\text{vld}},   \\*[0.5mm] 
&w_2^{\text{v}} = \Big({\sum_{j = 1}^{N_{\text{eig}}} \frac{2D_{jj}^2(\alpha_{\text{\tiny{NN}}}^{\text{v}} + \eta_\circ^2)}{(\alpha_{\text{\tiny{NN}}}^{\text{v}} + D_{jj}^2)^3} \big{|} (\boldsymbol{u}^*_j, \textcolor{black}{\boldsymbol{u}^{\text{v}}_{\text{\tiny L}}}) \big{|}^2 +\epsilon}\Big)^{-1}, 
\end{aligned}
\eeq

When $\eta_\circ$ is sufficiently close to the optimal threshold or in case the imaged area is relatively simple (e.g., when early reconstructions show a few well-separated scatterers), the training of R-Net may stop here (at the end of Step 1). Step 2 aims to further optimize the R-Net's output using network parameters obtained in Step 1 as an initial state.

\subsection{R-Net training:~Step 2}\label{TS2}

Given $\boldsymbol{F}^\delta$, $\lbrace{\boldsymbol{U}^*\boldsymbol{u}^{\text{t}}_{\text{\tiny L}}} \rbrace_{{\text{t}} = 1,2,\ldots,N_{\text{trn}}\!}$, and model parameters from Step 1, the learning objective in this step is to furnish regularization maps $\lbrace\alpha^{\text{t}}_{\text{\tiny NN}} \rbrace_{{\text{t}} = 1,2,\ldots,N_{\text{trn}}\!}$ that optimize the imaging (i.e.,~Tikhonov) loss function within the Bayes risk minimization framework,
\beq\lb{Jimg}
\mathfrak{J}_{\text{img}} = \frac{1}{{N_{\text{trn}}}}\sum_{\text{t} = 1}^{N_{\text{trn}}} |\hspace{-0.4mm}|\boldsymbol{F}^\delta \textcolor{black}{\boldsymbol{g}_{\text{\tiny{NN}}}^{\text{t}}} \hspace{-0.2mm} - \hspace{0.2mm} \textcolor{black}{\boldsymbol{u}^{\text{t}}_{\text{\tiny L}}}|\hspace{-0.4mm}|^2_{L^2} + \textcolor{black}{\alpha_{\text{\tiny{NN}}}^{\text{t}}} |\hspace{-0.4mm}|\hspace{0.2mm}  \textcolor{black}{\boldsymbol{g}_{\text{\tiny{NN}}}^{\text{t}}} |\hspace{-0.4mm}|^2_{L^2},
\eeq
wherein $\forall {\text{t}} \in 1,2,\ldots, N_{\text{trn}}$,
\beq\lb{gtNN}
\textcolor{black}{\boldsymbol{g}_{\text{\tiny{NN}}}^{\text{t}}}(\textcolor{black}{\alpha_{\text{\tiny{NN}}}^{\text{t}}}) \,=\, \boldsymbol{V} \textcolor{black}{\boldsymbol{D}^\dagger} \boldsymbol{U}^* \textcolor{black}{\boldsymbol{u}^{\text{t}}_{\text{\tiny L}}}, \quad \textcolor{black}{\boldsymbol{D}^\dagger}(\boldsymbol{D},\textcolor{black}{\alpha_{\text{\tiny{NN}}}^{\text{t}}}) \,=\, \text{diag}\lbrace \frac{D_{jj}}{\textcolor{black}{\alpha_{\text{\tiny{NN}}}^{\text{t}}}+D_{jj}^2} \rbrace, \quad j = 1,\ldots, {N_{\text{eig}}}.
\eeq

Here, $\boldsymbol{g}_{\text{\tiny{NN}}}^{\text{t}}$ is the well-known minimizer of~\eqref{Jimg} given the regularization parameter $\alpha_{\text{\tiny{NN}}}^{\text{t}}$. Observe that by optimizing~\eqref{Jimg}, R-Net aims to minimize the residual of scattering equation, while maximizing the LSM indicator functional. In this setting, Step 2 may enhance the image contrast through its search for a better approximation of min-norm solution at every sampling point. Unlike Step 1, Step 2 does not involve labeled datasets for training. As such, a careful regulation of the training process is required to prevent overfitting. For this purpose, we propose a criteria to stop training as outlined in Algorithm~\ref{AL2} which is based on the relative behavior of training and validation loss functions. In this step, the validation dataset may be defined by $\boldsymbol{F}^\delta$ and $\lbrace{\boldsymbol{U}^*\boldsymbol{u}^{\upnu}_{\text{\tiny L}}} \rbrace_{\upnu = 1,2,\ldots,N_{\upnu\!}} = \lbrace{\boldsymbol{U}^*\boldsymbol{u}^{\text{n}}_{\text{\tiny L}}} \rbrace_{{\text{n}} = 1,2,\ldots,N_{\text{p}}N_{\text{s}}\!} \setminus\lbrace{\boldsymbol{U}^*\boldsymbol{u}^{\text{t}}_{\text{\tiny L}}} \rbrace_{{\text{t}} = 1,2,\ldots,N_{\text{trn}}\!}$, whereby the validation loss is described by   
\beq\lb{Vimg}
\mathfrak{V}_{\text{img}} = \frac{1}{{N_{\upnu}}}\sum_{\upnu = 1}^{N_{\text{trn}}} |\hspace{-0.4mm}|\boldsymbol{F}^\delta \textcolor{black}{\boldsymbol{g}_{\text{\tiny{NN}}}^{\upnu}} \hspace{-0.2mm} - \hspace{0.2mm} \textcolor{black}{\boldsymbol{u}^{\upnu}_{\text{\tiny L}}}|\hspace{-0.4mm}|^2_{L^2} + \textcolor{black}{\alpha_{\text{\tiny{NN}}}^{\upnu}} |\hspace{-0.4mm}|\hspace{0.2mm}  \textcolor{black}{\boldsymbol{g}_{\text{\tiny{NN}}}^{\upnu}} |\hspace{-0.4mm}|^2_{L^2},
\eeq
where $\lbrace\alpha^{\upnu}_{\text{\tiny NN}} \rbrace_{{\upnu} = 1,2,\ldots,N_{\upnu}\!}$ is the R-Net output at every epoch when the input is $\lbrace{\boldsymbol{U}^*\boldsymbol{u}^{\upnu}_{\text{\tiny L}}} \rbrace_{\upnu = 1,2,\ldots,N_{\upnu\!}}$. It should also be mentioned that $\boldsymbol{g}_{\text{\tiny{NN}}}^{\upnu} = \boldsymbol{g}_{\text{\tiny{NN}}}^{\upnu}(\textcolor{black}{\alpha_{\text{\tiny{NN}}}^{\upnu}})$ is calculated similar to~\eqref{gtNN}. To compare the training and validation loss trajectories, i.e.,~$\mathfrak{J}_{\text{img}}(t)$ and $\mathfrak{V}_{\text{img}}(t)$, both are normalized at every epoch by their respective values at the first epoch of Step 2, denoted by epoch1 in Algorithm~\ref{AL2},
\beq\lb{hJV}
\hat{\mathfrak{J}}_{\text{img}}(t) = \mathfrak{J}_{\text{img}}(t)/\mathfrak{J}_{\text{img}}(\text{epoch1}), \quad \hat{\mathfrak{V}}_{\text{img}}(t) = \mathfrak{V}_{\text{img}}(t)/\mathfrak{V}_{\text{img}}(\text{epoch1}).
\eeq

Next, normal loss variations per epoch are computed as the following
\beq\lb{DhJV}
\Delta\hat{\mathfrak{J}}_{\text{img}}(t) = \hat{\mathfrak{J}}_{\text{img}}(t) - \hat{\mathfrak{J}}_{\text{img}}(t-1), \quad \Delta\hat{\mathfrak{V}}_{\text{img}}(t) = \hat{\mathfrak{V}}_{\text{img}}(t) - \hat{\mathfrak{V}}_{\text{img}}(t-1).
\eeq

Now, let $\langle\Delta\hat{\mathfrak{J}}_{\text{img}}\rangle(t)$ and $\langle\Delta\hat{\mathfrak{V}}_{\text{img}}\rangle(t)$ denote the root mean square (rms) of normal loss variations within $[t_\circ \,\,\, t)$ wherein $t_\circ = \max(t-N_{\text{rms}}, \text{epoch1})$. This provides a more stable metric of the overall trend of training and validation losses around $t$. On denoting the absolute and relative thresholds by $\sigma_{\text{a}}$ and $\sigma_{\text{r}}$, respectively, the training in Step 2 stops if the rms of (either training or validation) loss variations falls below $\sigma_{\text{a}}$, i.e.,
\beq\lb{SC1}
\langle\Delta\hat{\mathfrak{V}}_{\text{img}}\rangle(t) < \sigma_{\text{a}} \,\,\, \lor \,\,\, \langle\Delta\hat{\mathfrak{J}}_{\text{img}}\rangle(t) < \sigma_{\text{a}}.
\eeq

This criteria is concerned with the overall trajectory of each loss separately. The training also stops if the relative variation $$\varrho(t) = \langle\Delta\hat{\mathfrak{V}}_{\text{img}}\rangle(t) / \langle\Delta\hat{\mathfrak{J}}_{\text{img}}\rangle(t)$$
falls under $\sigma_{\text{r}}$ or grows beyond $1/\sigma_{\text{r}}$, i.e., 
\beq\lb{SC2}
 \varrho(t) < \sigma_{\text{r}} \,\,\, \lor \,\,\, \varrho(t) > 1/\sigma_{\text{r}},
 \eeq
which respectively indicate overfitting and underfitting. This criteria is critical for efficient image enhancement in Step 2 as training over a large number of epochs in this step does not necessary lead to an improved image and may in fact deteriorate the original reconstruction largely due to overfitting. The pseudocode for training R-Nets, including both steps, is provided in Algorithm~\ref{AL1}.   

\begin{algorithm}[h!]
\caption{Deep regularization networks
 for inverse problems with noisy operators}
\label{AL1}
  \begin{algorithmic}
    \REQUIRE (1)~discretized operator $\bF^\delta$ and its SVD decomposition $\bF^\delta=\boldsymbol{U} \boldsymbol{D} \boldsymbol{V}^*$ \\\hspace{1.35cm}~(2)~right-hand side dictionary $\lbrace\bu^{\text{n}}_{\text{\tiny L}} \rbrace_{{\text{n}} = 1,2,\ldots,N_{\text{p}}N_{\text{s}}}$ \\\hspace{1.35cm}~(3)~mode of training (basic or informed) in Step 1  
    \vspace*{0.75 mm}
    \STATE Downsample the right-hand side to $\lbrace\bu^{\text{t}}_{\text{\tiny L}} \rbrace_{{\text{t}} = 1,2,\ldots,N_{\text{trn}}\!}$ for training 
    \vspace*{0.75 mm}
    \STATE Compute the low-resolution and approximate regularization maps $\lbrace\alpha^{\text{t}}_{\text{\tiny M}} \rbrace_{{\text{t}} = 1,2,\ldots,N_{\text{trn}}\!}$ using the discrepancy principle with a rough estimate $\eta_\circ$ for the residual misfit 
    \vspace*{0.75 mm}
    \STATE {\bf{Step 1:}}
    \FOR{$t = 0$ \textbf{to} epoch1}
    \vspace*{0.75 mm}
     \IF{training mode is basic}
      \vspace*{0.75 mm}
      \STATE Based on the network output $\lbrace\alpha^{\text{t}}_{\text{\tiny NN}} \rbrace_{{\text{t}} = 1,2,\ldots,N_{\text{trn}}\!}$, update \vspace*{-2 mm} $$ \mathfrak{J}^b_{\text{dscr}} = \frac{1}{{N_{\text{trn}}}}\sum_{\text{t} = 1}^{N_{\text{trn}}} (\textcolor{black}{\alpha_{\text{\tiny{NN}}}^{\text{t}}} \hspace{-0.2mm} - \hspace{0.2mm} \textcolor{black}{\alpha_{\text{\tiny{M}}}^{\text{t}}})^2 \vspace*{-7 mm}  $$ 
      \vspace*{0 mm}  
      \STATE Optimize $\mathfrak{J}^b_{\text{dscr}}$
      \vspace*{0.75 mm}
    \ELSE
    \vspace*{0.75 mm}
    \STATE Based on the network output $\lbrace\alpha^{\text{t}}_{\text{\tiny NN}} \rbrace_{{\text{t}} = 1,2,\ldots,N_{\text{trn}}\!}$, compute $\forall {\text{t}}$,
       \vspace*{0.75 mm}
    \vspace*{0.75 mm}
        \STATE  $J_1^{\xxs\text{t}} =  | \frac{\alpha_{\text{\tiny{NN}}}^{\text{t}}-\alpha_{\text{\tiny{M}}}^{\text{t}}}{\alpha_{\text{\tiny{N}}}^{\text{t}}} |^2, \quad   \alpha_{\text{\tiny{N}}}^{\text{t}} \,=\, \alpha_{\text{\tiny{M}}}^{\text{t}} \,\,\, \text{or}\,\, \max(\alpha_{\text{\tiny{M}}}^{\text{t}}), \quad \forall {\text{t}} | \,\, \alpha_{\text{\tiny{N}}}^{\text{t}} > 0 $ 
      \vspace*{0.75 mm}
      \STATE $ J_2^{\xxs\text{t}}  = \big( \sum_{j = 1}^{N_{\text{eig}}} f_{j}^{\xxs\text{t}}(\alpha_{\text{\tiny{NN}}}^{\text{t}}, \eta_\circ) \big)^2, \quad f_{j}^{\xxs\text{t}}(\alpha_{\text{\tiny{NN}}}^{\text{t}}, \eta_\circ) = \frac{\alpha_{\text{\tiny{NN}}}^{\text{t}} - \eta_\circ^2 D_{jj}^2}{(\alpha_{\text{\tiny{NN}}}^{\text{t}} + D_{jj}^2)^2} \big{|} (\boldsymbol{u}^*_j, \textcolor{black}{\boldsymbol{u}^{\text{t}}_{\tiny \text{L}}}) \big{|}^2$
      \vspace*{0.75 mm}
      \STATE $w_1^{\text{t}} = \alpha_{\text{\tiny{N}}}^{\text{t}}, \quad w_2^{\text{t}} = \big({\sum_{j = 1}^{N_{\text{eig}}} \frac{\partial f_{j}^{\xxs\text{t}}}{\partial \alpha_{\text{\tiny{NN}}}^{\text{t}}}(\alpha_{\text{\tiny{NN}}}^{\text{t}}, \eta_\circ) +\epsilon}\big)^{-1}$
      \vspace*{0.75 mm}
       \vspace*{0.75 mm}
      \STATE Update \vspace*{-4 mm} $$ \mathfrak{J}^i_{\text{dscr}} = \frac{1}{{N_{\text{trn}}}}\sum_{\text{n} = 1}^{N_{\text{trn}}} w_1^{\text{t}} J_1^{\xxs\text{t}} + w_2^{\text{t}} J_2^{\xxs\text{t}} \vspace*{-7 mm}  $$ 
      	\vspace*{0.75 mm}  
      \STATE Optimize $\mathfrak{J}^i_{\text{dscr}}$
   	\ENDIF 
	\ENDFOR 
    \vspace*{0.5 mm}
     \STATE {\bf{Step 2:}}
    \FOR{$t =$ epoch1$+1$ \textbf{to} epoch2}
   \STATE Based on the network output $\lbrace\alpha^{\text{t}}_{\text{\tiny NN}} \rbrace_{{\text{t}} = 1,2,\ldots,N_{\text{trn}}\!}$, compute $\forall {\text{t}}$, \vspace*{-2 mm} $$\textcolor{black}{\boldsymbol{g}_{\text{\tiny{NN}}}^{\text{t}}}(\textcolor{black}{\alpha_{\text{\tiny{NN}}}^{\text{t}}}) \,=\, \boldsymbol{V} \textcolor{black}{\boldsymbol{D}^\dagger} \boldsymbol{U}^* \textcolor{black}{\boldsymbol{u}^{\text{t}}_L}, \quad \textcolor{black}{\boldsymbol{D}^\dagger}(\boldsymbol{D},\textcolor{black}{\alpha_{\text{\tiny{NN}}}^{\text{t}}}) \,=\, \text{diag}\lbrace \frac{D_{jj}}{\textcolor{black}{\alpha_{\text{\tiny{NN}}}^{\text{t}}}+D_{jj}^2} \rbrace, \quad j = 1,\ldots, {N_{\text{eig}}} \vspace*{-5 mm} $$
   \STATE Update \vspace*{-2 mm} $$ \mathfrak{J}_{\text{img}} = \frac{1}{{N_{\text{trn}}}}\sum_{\text{t} = 1}^{N_{\text{trn}}} |\hspace{-0.4mm}|\boldsymbol{F}^\delta \textcolor{black}{\boldsymbol{g}_{\text{\tiny{NN}}}^{\text{t}}} \hspace{-0.2mm} - \hspace{0.2mm} \textcolor{black}{\boldsymbol{u}^{\text{t}}_L}|\hspace{-0.4mm}|^2_{L^2} + \textcolor{black}{\alpha_{\text{\tiny{NN}}}^{\text{t}}} |\hspace{-0.4mm}|\hspace{0.2mm}  \textcolor{black}{\boldsymbol{g}_{\text{\tiny{NN}}}^{\text{t}}} |\hspace{-0.4mm}|^2_{L^2} \vspace*{-5 mm} $$
    \vspace*{0.5 mm}
   \STATE Optimize $\mathfrak{J}_{\text{img}}$
    \vspace*{0.5 mm}
  \STATE Compute the validation loss $\mathfrak{V}_{\text{img}}$ according to~\eqref{Vimg}
   \vspace*{0.5 mm}
   \IF{stop training flag $s_t = 1$ per Algorithm~\ref{AL2}}
   \STATE STOP
   \ENDIF  
    \ENDFOR
  \end{algorithmic} 
\end{algorithm}

\begin{algorithm}[h!]
\caption{Stop training criteria at epoch $t$}
\label{AL2}
  \begin{algorithmic}
   \REQUIRE (1)~discretized operator $\bF^\delta$ and its SVD decomposition $\bF^\delta=\boldsymbol{U} \boldsymbol{D} \boldsymbol{V}^*$\\\hspace{1.35cm}~(2)~right-hand side dictionary $\lbrace\bu^{\text{n}}_{\text{\tiny L}} \rbrace_{{\text{n}} = 1,2,\ldots,N_{\text{p}}N_{\text{s}}}$\\\hspace{1.35cm}~(3)~training and validation loss values $\mathfrak{J}_{\text{img}}, \mathfrak{V}_{\text{img}}$ at epochs $t-N_{\text{rms}}$ to $t$\\\hspace{1.35cm}~(4)~relative and absolute thresholds $\sigma_{\text{r}}$ and $\sigma_{\text{a}}$  
  \STATE Compute normal loss values at epoch $t$, \vspace*{-1 mm} $$\hat{\mathfrak{J}}_{\text{img}}(t) = \mathfrak{J}_{\text{img}}(t)/\mathfrak{J}_{\text{img}}(\text{epoch1}), \,\, \hat{\mathfrak{V}}_{\text{img}}(t) = \mathfrak{V}_{\text{img}}(t)/\mathfrak{V}_{\text{img}}(\text{epoch1})  \vspace*{-1 mm} $$
  \STATE Compute normal loss variations, \vspace*{-1 mm} $$\Delta\hat{\mathfrak{J}}_{\text{img}}(t) = \hat{\mathfrak{J}}_{\text{img}}(t) - \hat{\mathfrak{J}}_{\text{img}}(t-1), \,\, \Delta\hat{\mathfrak{V}}_{\text{img}}(t) = \hat{\mathfrak{V}}_{\text{img}}(t) - \hat{\mathfrak{V}}_{\text{img}}(t-1)  \vspace*{-1 mm}  $$
  \STATE Compute root mean square of variations given $t_\circ = \max(t-N_{\text{rms}},\text{epoch1})$, \vspace*{-1 mm} $$\langle\Delta\hat{\mathfrak{J}}_{\text{img}}\rangle(t) = \sqrt{\frac{1}{t-t_\circ}\sum_{\tau =\xxs t_\circ}^{t-1} \big(\Delta\hat{\mathfrak{J}}_{\text{img}}(\tau)\big)^2}, \,\,  \langle\Delta\hat{\mathfrak{V}}_{\text{img}}\rangle(t) = \sqrt{\frac{1}{t-t_\circ}\sum_{\tau = \xxs t_\circ}^{t-1} \big(\Delta\hat{\mathfrak{V}}_{\text{img}}(\tau)\big)^2} \vspace*{-1 mm} $$
  \STATE Compute relative variation $\varrho(t) = \langle\Delta\hat{\mathfrak{V}}_{\text{img}}\rangle(t) / \langle\Delta\hat{\mathfrak{J}}_{\text{img}}\rangle(t)$  \vspace*{1 mm} 
  \IF{$\varrho(t) < \sigma_{\text{r}}$ \,or\, $\varrho(t) > 1/\sigma_{\text{r}}$ \,or\, $\langle\Delta\hat{\mathfrak{V}}_{\text{img}}\rangle(t) < \sigma_{\text{a}}$ \,or\, $\langle\Delta\hat{\mathfrak{J}}_{\text{img}}\rangle(t) < \sigma_{\text{a}}$}
   \STATE stop training flag $s_t = 1$
   \ELSE
   \STATE stop training flag $s_t = 0$
  \ENDIF
  \end{algorithmic} 
\end{algorithm}
    
\section{Implementation and results}\lb{IR}   

This section examines the performance of R-Nets according to Algorithm~\ref{AL1} and~\ref{AL2} through a set of numerical experiments. Here, images constructed by R-Net are compared to those obtained by the Morozov discrepancy principle. In what follows the synthetic sensory data, namely the scattered fields and RHS patterns, are simulated by a computational platform based on the elastodynamic boundary integral equations, see~\cite{pour2020,Bon1999} for details of the computational method and simulations.

\subsection{Testing configuration}

As illustrated in Fig.~\ref{SetNum1}, an elastic plate of dimensions $3$ $\!\times\!$ $3$ $\!\times\!$ $0.02$ featuring a randomly cracked damage zone is modeled. The shear modulus, mass density, and Poisson's ratio of the plate are $\mu = 1$, $\rho = 1$ and $\nu = 0.25$, whereby the shear and compressional wave speeds are calculated as $c_s = 1$ and $c_p = 1.73$. Evolving in five time steps $t_1-t_5$, the damage zone is comprised of randomly distributed cracks $\Gamma_{1}-\Gamma_{15}$ that are hidden within the thickness of the plate. As such,~3D simulations are required to model wave motion in the specimen. A detailed description of scatterers including the center $(x_c, y_c)$, length $\ell$, and orientation $\phi$ (with respect to $x$ axis) of each crack $\Gamma_{\kappa}$, $\kappa = \lbrace 1, 2, ..., 15 \rbrace$ is provided in Table~\ref{Num1}. All fractures in this configuration are traction-free.

\begin{figure}[!tp]
\center\includegraphics[width=0.75\linewidth]{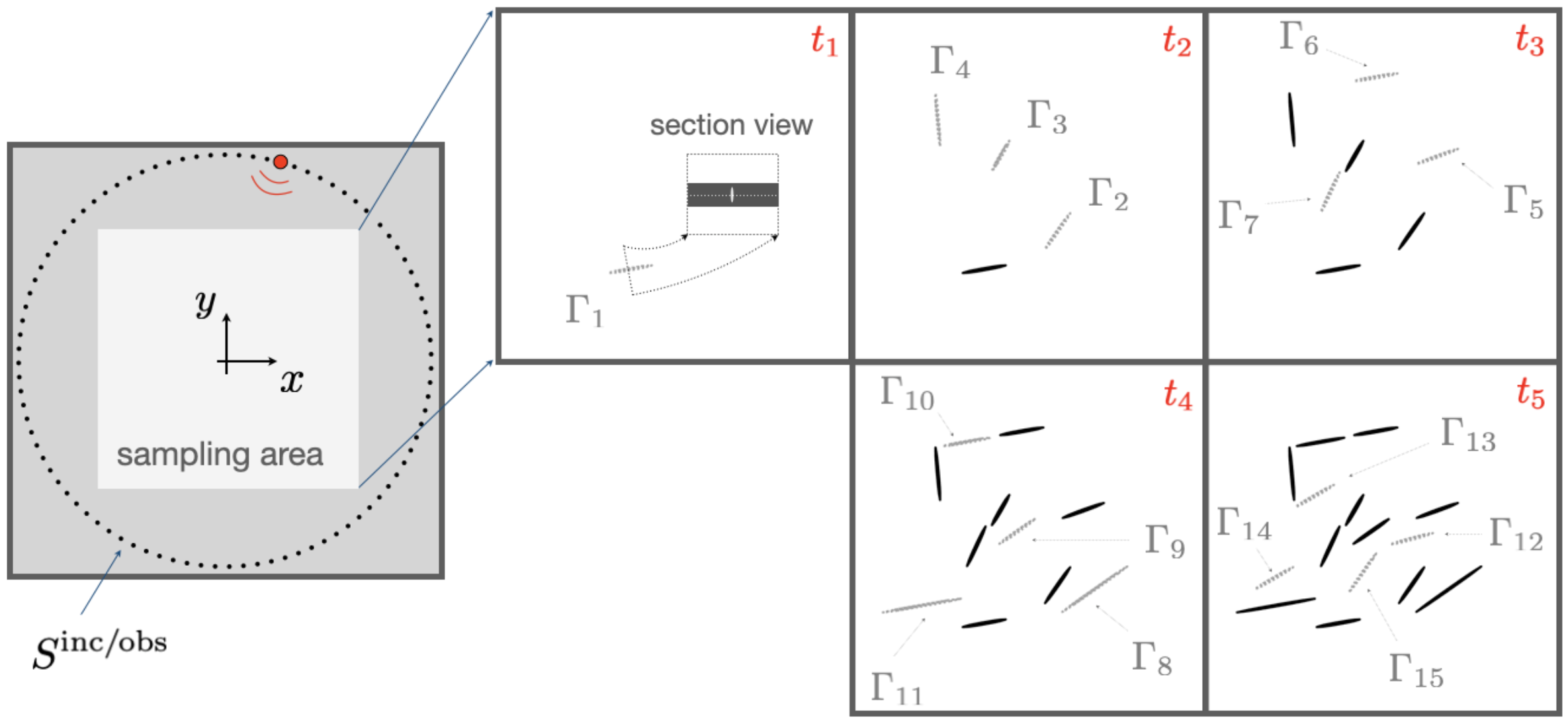} \vspace*{-3mm} 
\caption{Sensing configuration of synthetic experiments on an elastic plate (left) featuring a damage zone comprised of randomly distributed cracks $\Gamma_{1}-\Gamma_{15}$ evolving in five time steps ($t_1-t_5$) within the thickness of the specimen according to the sectional view shown at $t_1$. } \lb{SetNum1}
\vspace*{-1.5mm}
\end{figure} 

\begin{table}[!h]
\begin{center}
\caption{\small Damage zone configuration:~center $(x_c, y_c)$, length $\ell$, and orientation $\phi$ (with respect to $x$ axis) of cracks $\Gamma_{\kappa}$, $\kappa = \lbrace 1, 2, ..., 15 \rbrace$, shown in Fig.~\ref{SetNum1}.} \vspace*{-2mm}
\label{Num1}
\resizebox{0.7\textwidth}{!}{%
\begin{tabular}{|c|c|c|c|c|c|c|c|c|} \hline
\!\!$\kappa$\!\! & 1 & 2 & 3 & 4&5 & 6 & 7 & 8  \\ \hline\hline 
\!\!$x_{\text{c}}(\Gamma_\kappa)$\!\!    & \!\!$-0.33$\!\!  & \!\!$0.21$\!\!  & \!\!$-0.21$\!\! & \!\!$-0.68$\!\! & \!\!$0.4$\!\! & \!\!$-0.05$\!\! & \!\!$-0.39$\!\!  & \!\!$0.49$\!\!  \\ 
 \hline
 \!\!$y_{\text{c}}(\Gamma_\kappa)$\!\!    & \!\!$-0.62$\!\!  & \!\!$-0.34$\!\!  & \!\!$0.22$\!\! & \!\!$0.49$\!\! & \!\!$0.21$\!\! & \!\!$0.8$\!\! & \!\!$-0.05$\!\!  & \!\!$-0.37$\!\!  \\ 
 \hline
\!\!$\ell \exs (\Gamma_\kappa)$\!\!   & \!\!$1/3$\!\!  & \!\!$1/3$\!\! & \!\!$1/4$\!\! & \!\!$2/5$\!\! & \!\!$1/3$\!\! & \!\!$1/3$\!\! & \!\!$1/3$\!\!  &  \!\!$3/5$\!\!  \\ 
 \hline
 \!\!$\phi \exs (\Gamma_\kappa)$\!\!   & \!\!$\pi/18$\!\!  & \!\!$11\pi/36$\!\! & \!\!$\pi/3$\!\! & \!\!$19\pi/36$\!\! & \!\!$\pi/9$\!\!  & \!\!$\pi/18$\!\!  & \!\!$13\pi/36$\!\!  &  \!\!$7\pi/36$\!\! \\ 
 \hline\hline
\!\!$\kappa$\!\! & 9 & 10 & 11 & 12 & 13 & 14 & 15 &  \\ \hline\hline 
\!\!$x_{\text{c}}(\Gamma_\kappa)$\!\!  & \!\!$-0.09$\!\! & \!\!$-0.46$\!\! & \!\!$-0.8$\!\! & \!\!$0.21$\!\!  & \!\!$-0.5$\!\!  & \!\!$-0.8$\!\!  & \!\!$-0.15$\!\! &  \\ 
 \hline
 \!\!$y_{\text{c}}(\Gamma_\kappa)$\!\!  & \!\!$0.06$\!\! & \!\!$0.72$\!\! & \!\!$-0.5$\!\! & \!\!$0$\!\!  & \!\!$0.32$\!\!  & \!\!$-0.29$\!\!  & \!\!$-0.25$\!\! &  \\ 
 \hline
\!\!$\ell \exs (\Gamma_\kappa)$\!\! & \!\!$1/3$\!\! & \!\!$2/5$\!\! & \!\!$3/5$\!\! & \!\!$1/3$\!\!  & \!\!$1/3$\!\!  & \!\!$1/3$\!\! & \!\!$7/20$\!\! & \\ 
 \hline
 \!\!$\phi \exs (\Gamma_\kappa)$\!\!  & \!\!$7\pi/36$\!\! & \!\!$\pi/18$\!\! & \!\!$\pi/18$\!\! & \!\!$\pi/12$\!\!  & \!\!$\pi/6$\!\!  & \!\!$\pi/6$\!\! & \!\!$11\pi/36$\!\! & \\ 
 \hline
\end{tabular}}
\end{center}
\vspace*{-1mm}
\end{table}

\subsection{Forward simulations}
Numerical experiments are conducted in five steps at $t = \lbrace t_1, t_2, ..., t_5 \rbrace$ when the specimen assumes the associated configurations shown in Fig.~\ref{SetNum1}~($t_1-t_5$). Every sensing step involves in-plane harmonic excitations on a set of source points over the incident grid $S^{\text{inc}}$. The excitation frequency \mbox{$\omega = 72$ rad/s} induces shear wavelength $\lambda_s = 0.08$ in the specimen. The incident wave interacts with fractures at each $t_k$, $k = 1,2,\ldots, 5$, giving rise to a scattered field governed by the Navier equations -- whose footprint ${\bv_k^\text{obs}} = [\bu{\obs} - \exs\bu^{\text{f}}]_k$ over the observation grid $S^{\text{obs}}$ is computed. As mentioned earlier, due to the specimen's asymmetry through the thickness, our simulations are performed in three dimensions via an elastodynamic code rooted in the boundary element method~\cite{Bon1999, pour2020}. However, only the in-plane components of the computed scattered fields, in \mbox{the $x-y$} plane, are used for the reconstructions. The incident/observation grid is a circle of radius $1.45$ as shown in Fig.~\ref{SetNum1}.


\subsection{Data Inversion}

The discretized scattering equation is formed in two steps, by:~(1) constructing the discrete scattering operators ${\text{\bf F}}_k$ and ${\text{\bf F}}_k^\delta$, $k = 1,2,\ldots, 5$, and (2) computing the RHS patterns $\bu^{\text{n}}_\text{\tiny L}$, $\text{n} = 1,2,\ldots, N_{\text{p}}N_{\text{s}}$. This is followed by the calculation of LSM imaging indicator through minimization of the Tikhonov imaging loss as elucidated in the sequel. 

\subsubsection*{Construction of the discretized scattering operator}
For both illumination and sensing purposes, $S^{\text{inc/\nxs obs}}$ is sampled by a uniform grid of $N = 1000$ excitation/observation points. On denoting the polarization amplitude of excitation by $\bq$ and the orthonormal bases in \mbox{the $x-y$} plane by $(\be_1,\be_2)$, the discretized scattering operator ${\text{\bf F}}_k$, $k = 1,2,\ldots, 5$, takes the form of a $2N\!\times 2N$ matrix with components    
\beq\lb{mat2}
\textrm{\bf{F}}_{k}(2j\nxs+\nxs1\!:\!2j\nxs+\nxs2, \,2i\nxs+\nxs1\!:\!2i\nxs+\nxs2) ~=\, 
\left[\begin{array}{cccc}
\!V_{k}^{11} \!&\! V_{k}^{12}\!\!  \\*[1mm]
\!V_{k}^{21} \!&\! V_{k}^{22} \!\! \\
\end{array}\right] (\bxi_j,\bx_i),  \qquad i,j = 0,\ldots N-1,
\eeq
where~$V_{k}^{\iota\upsilon}(\bxi_j,\bx_i)$ $(\iota,\upsilon\!=\!1,2)$ is the $\iota^{\textrm{th}}$ component of the scattered displacement field measured at $\bxi_j \in S^{\text{obs}}$ due to a unit harmonic excitation applied at $\bx_i \in S^{\text{inc}}$ along the coordinate direction $\upsilon$ such that 
\beq\label{mat1}\nonumber
\bv_{k}\obs(\bxi_j) \,=\,
\textrm{\bf{V}}_{k} (\bxi_j,\bx_i)\exs \bq(\bx_i)
\eeq

To account for the presence of noise in measurements, we consider the perturbed operators
\begin{equation}\label{DFN}
\textrm{\bf{F}}_{k}^\delta \,\, \colon \!\!\!= \, (\boldsymbol{I} + \boldsymbol{N}^{\delta} ) \exs \textrm{\bf{F}}_{k}, \quad k = 1,2,\ldots, 5,
\end{equation}
where $\boldsymbol{I}$ is the $2N \times 2N$ identity matrix, and $\boldsymbol{N}^{\delta}$ is the noise matrix of commensurate dimensions whose components are uniformly-distributed (complex) random variables in $[-\delta, \, \delta]^2$. In this study, $\delta \in \lbrace 0.1, 0.25 \rbrace$. 

\subsubsection*{A physics-based library of patterns}

Shown in Fig.~\ref{SetNum1}, the sampling region is a square $[-0.8,0.8]^2$, concentric with the specimen $\Uppi$, which is probed by a uniform $N_{\text{p}} = 100 \!\times\! 100$ grid of sampling points~$\bx_{\small \circ}\!$ where the imaging indicator functional is evaluated. At every sampling point, a set of infinitesimal trial dislocations is nucleated characterized by their unit normal direction $\textrm{\bf{n}}_{\small \circ}$. The latter samples the unit semicircle of possible orientations at $N_{\text{s}} = 72$ directions. Accordingly, the right-hand side library involves $N_{\text{p}}N_{\text{s}} = 10000 \!\times\! 72$ patterns affiliated with the trial pairs $L^{\text{n}} = L^{\text{n}}(\bx_{\small \circ}^{\text{n}},\textrm{\bf{n}}_{\small \circ}^{\text{n}})$, ${\text{n}} = 1,2,\ldots, N_{\text{p}}N_{\text{s}}$. Each trial signature $\bu_{\text{\tiny L}}^{\text{n}}(\bxi_j)$ over the observation grid $\bxi_j \in S^{\text{obs}}$, $j = 1,2,\ldots, N$, is computed separately by solving
\beq\lb{PhiL2}
\begin{aligned}
&\nabla \sip (\bC \colon \! \nabla \bu_{\text{\tiny L}}^{\text{n}}) \,+\, \rho \exs \omega^2\bu_{\text{\tiny L}}^{\text{n}} ~=~ \bzero \,\, &\text{in}& \,\, {\Uppi}\,\backslash L^{\text{n}}(\bx_{\small \circ}^{\text{n}},\textrm{\bf{n}}_{\small \circ}^{\text{n}}), \\*[0.25mm]
&\bn \cdot \bC \exs \colon \!  \nabla  \bu_{\text{\tiny L}}^{\text{n}} ~=~ \delta (\bxi-\bx_{\small \circ}^{\text{n}}) \exs \textrm{\bf{n}}_{\small \circ}^{\text{n}}  \quad &\,\text{on}& \,\, L^{\text{n}}(\bx_{\small \circ}^{\text{n}},\textrm{\bf{n}}_{\small \circ}^{\text{n}}), \\*[0.25mm]
&\bn \cdot \bC \exs \colon \!  \nabla  \bu_{\text{\tiny L}}^{\text{n}} ~=~ \bzero  \quad &\,\text{on}& \,\, \partial \Uppi.
\end{aligned}     
\eeq
wherein $\bC = \bC(\mu,\nu)$ is the specimen's fourth-order elasticity tensor.

Once the RHS library $\lbrace \bu_{\text{\tiny L}}^{\text{n}} \rbrace_{{\text{n}} = 1,2\ldots, N_{\text{p}}N_{\text{s}}}$ is computed, the reconstruction follows by solving the scattering equation
${\boldsymbol{F}}_{\!\! k}^{\delta} {\bg}^{\text{n}} = {\bu}_{\text{\tiny L}}^{\text{n}}$ for all $\text{n} = 1,2,\ldots, N_{\text{p}} N_{\text{s}}$
through minimizing~\eqref{dllsm}. 

\subsection{Results \& discussion}

 \begin{figure}[!tp]
\center\includegraphics[width=0.92\linewidth]{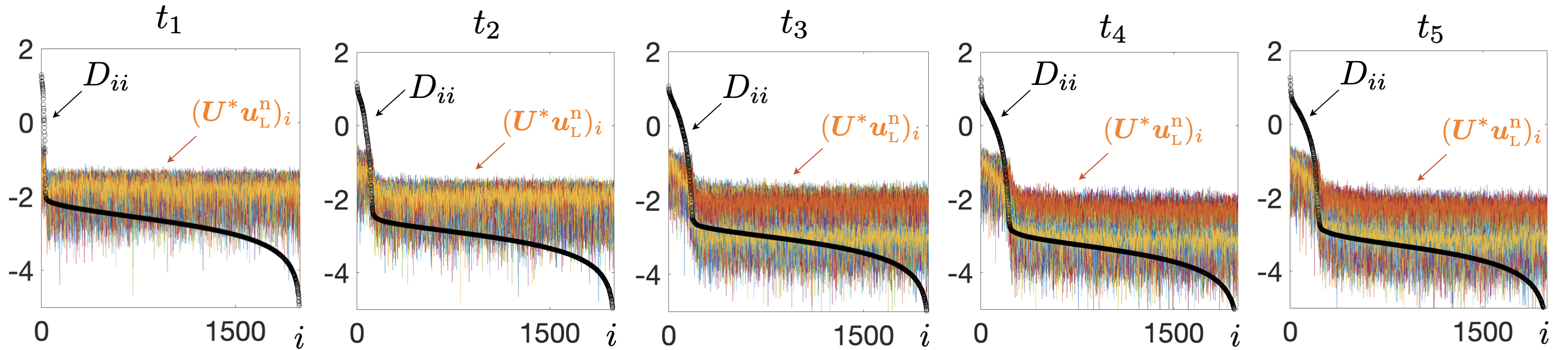} \vspace*{-3mm} 
\caption{Picard plots at sensing steps $t_k$, $k = 1,2,\ldots, 5$, wherein~$D_{ii}$, $i = 1,2,\ldots, 2000$, represents the $i^{\text{th}}$ eigenvalue of scattering operator $\boldsymbol{F}_k$, constructed from noiseless data at every $t_k$, while $(\boldsymbol{U}^*\bu^\text{n}_{\text{\tiny L}})_i = (\bu^*_i,\bu^\text{n}_{\text{\tiny L}})$ with $\text{n} = 1,2,\ldots, N_\text{trn}$ is the projected RHS patterns onto the $i^{\text{th}}$ left eigenvector of $\boldsymbol{F}_k$.} \lb{picard}
\vspace*{-1mm}
\end{figure} 

This section provides a comparative analysis of LSM reconstructions through minimizing~\eqref{dllsm} by way of three regularization approaches, namely: (i) Morozov discrepancy principle, (ii) basic R-Nets, and (iii) discrepancy-informed R-Nets. For future reference, all R-Nets in this section are architecturally similar MLPs, with input/output dimension $2N = 2000$, comprised of two fully-connected, ReLU-activated layers with $1000$ and $500$ neurons, respectively. In addition, Fig.~\ref{picard} shows the Picard plots affiliated with the scattering equation  
$$ \boldsymbol{F}_{k} \exs {\bg}^{\text{n}} = {\bu}_{\text{\tiny L}}^{\text{n}}, \quad \text{n} = 1,2,\ldots, N_{\text{trn}}, $$
for all $k = 1,2,\ldots,5$, depicting the distribution of eigenvalues of $\boldsymbol{F}_{k}$ against the projected RHS patterns and their evolution over time. Fig.~\ref{picard} highlights the critical role and complexity of selecting spectral filters for robust data inversion especially at later sensing steps when the domain becomes more complex and the two distributions in the Picard plots become less discernible.    

\subsubsection*{Reconstructions via the Morozov discrepancy principle}
 
The Morozov regularization maps $\alpha_{\text{\tiny M}} = \alpha_{\text{\tiny M}}(\bx_{\small \circ},\textrm{\bf{n}}_{\small \circ})$ are obtained for a given threshold $\eta$ by numerically solving 
 \beq\lb{MM}  
 \sum_{j = 1}^{N_{\text{eig}}} \frac{\alpha_{\text{\tiny{M}}}^{\text{n}} - \eta^2 D_{jj}^2}{(\alpha_{\text{\tiny{M}}}^{\text{n}} + D_{jj}^2)^2} \big{|} (\boldsymbol{u}^*_j, \textcolor{black}{\boldsymbol{u}^{\text{n}}_{\text{\tiny L}}}) \big{|}^2 ~=~0, \quad \boldsymbol{u}^{\text{n}}_{\text{\tiny L}} = \boldsymbol{u}^{\text{n}}_{\text{\tiny L}}(\bx_{\small \circ}^{\text{n}},\textrm{\bf{n}}_{\small \circ}^{\text{n}}), \quad N_{\text{eig}} = 2000,
 \eeq
for all $\text{n} = 1,2,\ldots, N_{\text{p}} N_{\text{s}} = 7.2 \!\times\! 10^5$. Here, $\eta$ is manually gauged for every reconstruction by computing the Morozov maps $\alpha_{\text{\tiny M}}$ for a set of thresholds $\eta = \lbrace 0.01, 0.02, 0.03, \ldots, 0.4 \rbrace$ and their associated LSM reconstructions $\mathcal{L}_{\text{\tiny M}}$. A proper estimate for $\eta$ is then selected through a qualitative (visual) comparison of the reconstructions. Table~\ref{MT} provides the manually tuned values of $\eta_k^{\text{\tiny np}}$ for every $t_k$, $k = 1,2,\ldots, 5$, at three noise levels $\text{np} = \lbrace 0, 10, 25 \rbrace \%$. Figs.~\ref{M0},~\ref{M10pn}, and~\ref{M25pn} illustrate the corresponding Morozov maps $\alpha_{\text{\tiny M}}$ and LSM reconstructions $\mathcal{L}_{\text{\tiny M}}$.

\renewcommand{\arraystretch}{1.2}
 \begin{table}[!h]
\begin{center}
\caption{\small{Manually gauged Morozov thresholds $\eta_k^{\text{\tiny np}}$ at sensing steps $t_k$, $k = 1,2,\ldots, 5$, for three noise levels at $\text{np} = \lbrace 0, 10, 25 \rbrace \%$.}} \vspace*{-2mm}
\label{MT}
\begin{tabular}{c||c|c|c|c|c} 
\diagbox{{\small{threshold}}}{{\!\small{time}}} & $t_1$\! & $t_2$\! & $t_3$\! & $t_4$\! & $t_5$\! \\ \hline\hline  
$\eta_k^{\text{\tiny 0}}$    & $0.3$  & $0.2$  & $0.08$ & $0.08$ & $0.09$ \\  
 \hline
$\eta_k^{\text{\tiny 10}}$    & \!\!$0.3$\!\!  & \!\!$0.23$\!\!  & \!\!$0.14$\!\! & \!\!$0.14$\!\! & \!\!$0.15$\!\! \\ 
 \hline
$\eta_k^{\text{\tiny 25}}$   & \!\!$0.3$\!\!  & \!\!$0.25$\!\! & \!\!$0.27$\!\! & \!\!$0.26$\!\! & \!\!$0.26$\!\! \\
\end{tabular}
\end{center}
\vspace*{-2mm}
\end{table}

 \begin{figure}[!tp]
\center\includegraphics[width=0.95\linewidth]{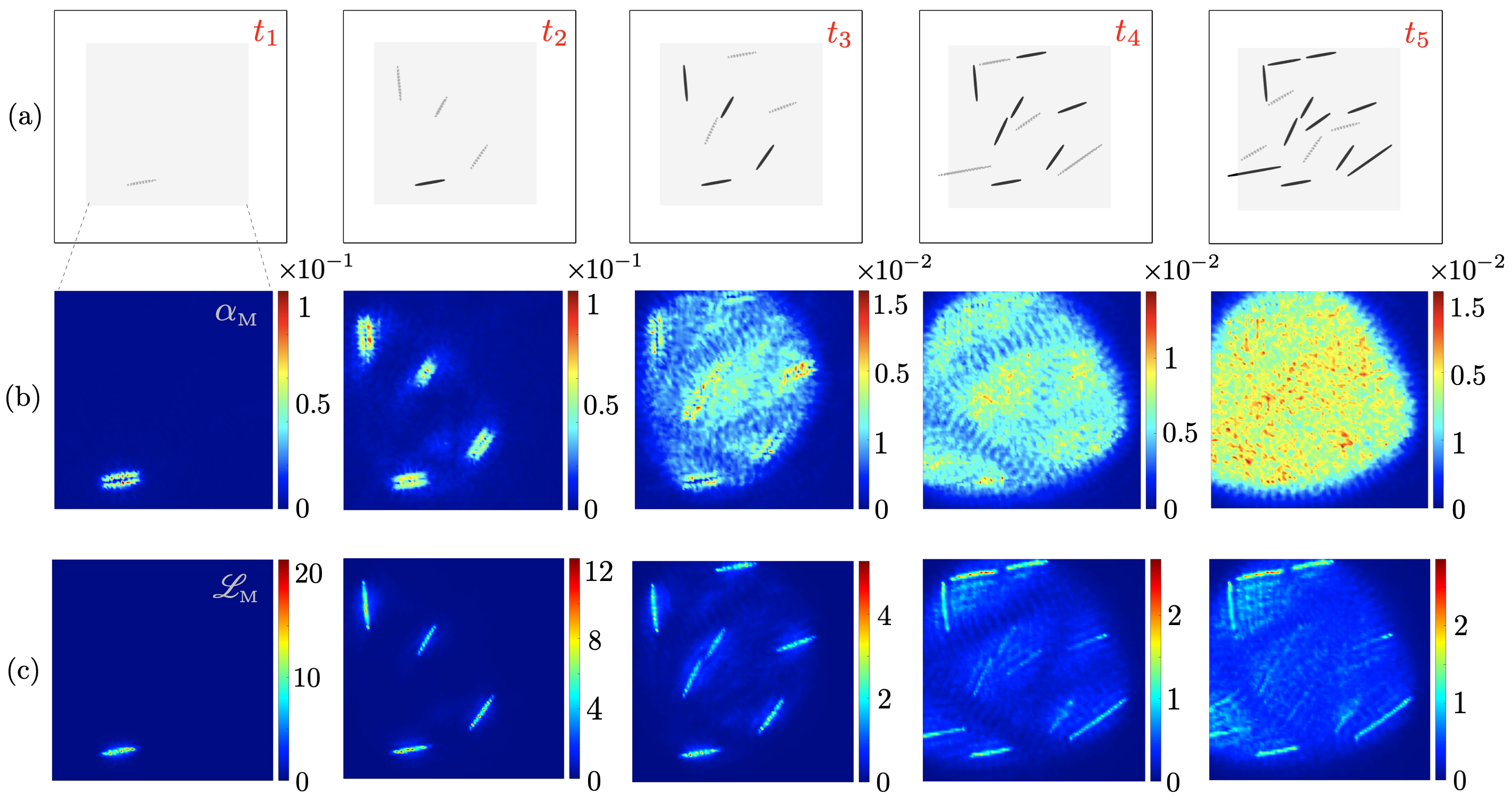} \vspace*{-3mm} 
\caption{LSM reconstructions from noiseless data by way of the Morozov discrepancy principle:~(a)~ground-truth configurations at sensing steps $t_k$, $k = 1,2,\ldots, 5$, (b) manually optimized Morozov regularization maps $\alpha_{\text{M}}$ over the dense search grid (of $100 \!\times\! 100$ sampling points) at $t_k$, and (c) the corresponding LSM indicator maps $\mathcal{L}_{\text{M}}$. } \lb{M0}
\vspace*{-4mm}
\end{figure} 

Observe that as the noise level increases, the regularization maps evolve and the colormaps are rescaled to reflect higher maximum values.

\subsubsection*{Reconstructions via basic R-Nets}

With reference to Algorithm~\ref{AL1}, basic R-Nets are trained by minimizing~\eqref{Jbds} and~\eqref{Jimg} in Steps 1 and 2, respectively. To generate the training dataset, assuming that the manually optimized residual misfits of Table~\ref{MT} are not available, the \emph{approximate} Morozov maps $\lbrace {\alpha_{\text{\tiny{M}}}^{\text{n}}} \rbrace$ are computed for every sensing step on the original grid of trial pairs $(\bx_{\small \circ}^\text{n},\textrm{\bf{n}}_{\small \circ}^\text{n})$, $\text{n} = 1,2,\ldots, 7.2 \!\times\! 10^5$ according to~\eqref{MM} by setting $\eta = \eta_\circ = 0.3$. The grid is then reduced by (a) selecting one trial scatterer per sampling point following the logic of~\eqref{Lp}, and (b) uniformly downsampling the search region by a factor of 2 in both $x$ and $y$ directions so that the training grid includes $50 \!\times\! 50$ sampling points. As such, the low-resolution Morozov maps for training are generated via $\lbrace {\alpha_{\text{\tiny{M}}}^{\text{t}}} \rbrace_{\text{t} = 1,2,\ldots, N_{\text{trn}}}$, $N_{\text{trn}} = 2500$. In this case, only $0.35\%$ of ($7.2 \!\times\! 10^5$) RHS patterns participate in network training. Keeping the selected trial scatterer for every sampling point, one may then select the complement of the training grid on the original ($100 \!\times\! 100$) image support to form the validation dataset $\lbrace {\alpha_{\text{\tiny{M}}}^{\text{v}}} \rbrace_{\text{v} = 1,2,\ldots, N_{\text{vld}}}$, $N_{\text{vld}} = 7500$. 
\begin{figure}[!tp]
\center\includegraphics[width=0.95\linewidth]{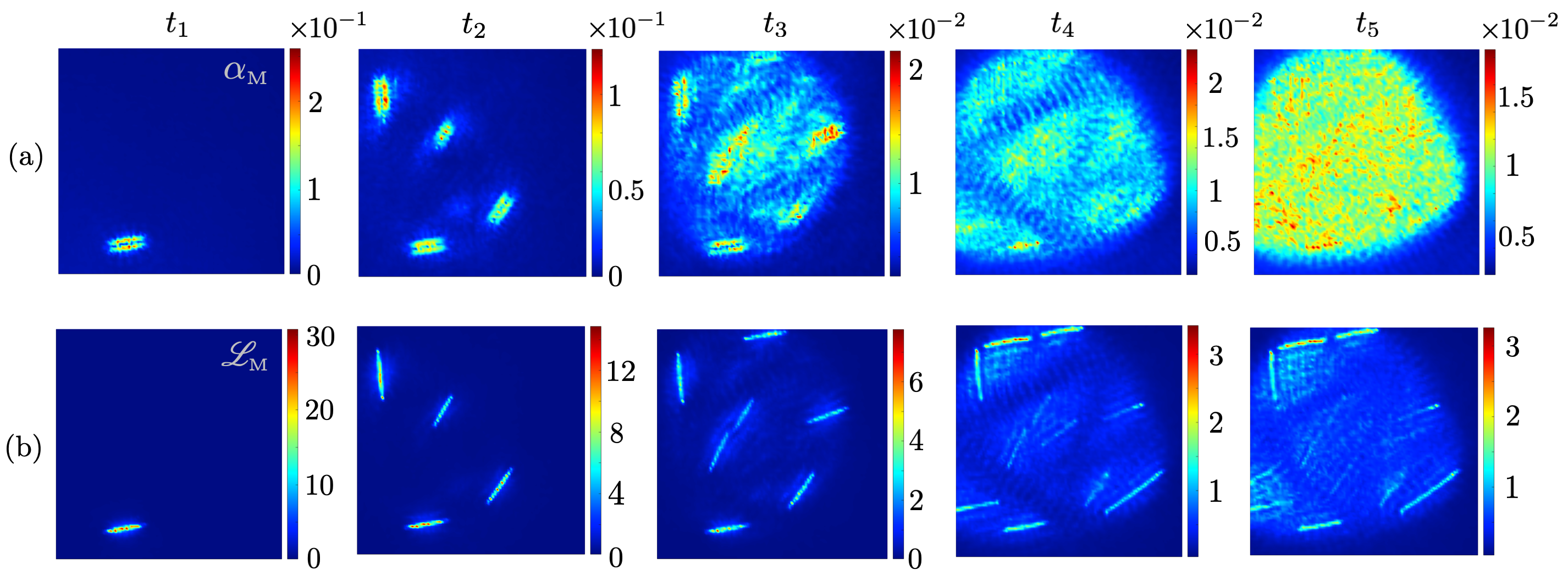} \vspace*{-3mm} 
\caption{LSM reconstructions from $10\%$ noisy data using the Morozov discrepancy principle:~(a) manually optimized Morozov regularization maps $\alpha_{\text{M}}$ over the dense search grid (of $100 \!\times\! 100$ sampling points) at $t_k$, $k = 1,2,\ldots, 5$, and (b) the corresponding LSM indicator maps $\mathcal{L}_{\text{M}}$.  } \lb{M10pn}
\vspace*{1mm}
\center\includegraphics[width=0.95\linewidth]{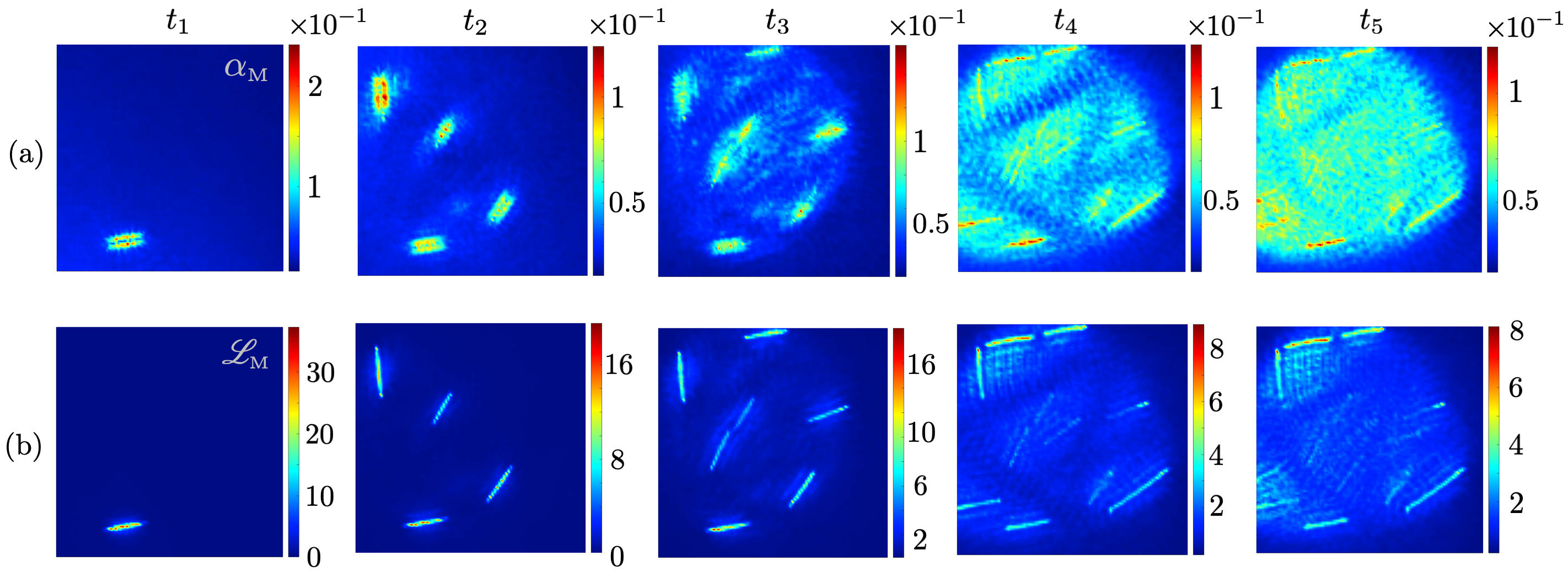} \vspace*{-3mm} 
\caption{LSM reconstructions from $25\%$ noisy data using the Morozov discrepancy principle:~(a) manually optimized Morozov regularization maps $\alpha_{\text{M}}$ over the dense search grid (of $100 \!\times\! 100$ sampling points) at $t_k$, $k = 1,2,\ldots, 5$, and (b) the corresponding LSM indicator maps $\mathcal{L}_{\text{M}}$. } \lb{M25pn}
\vspace*{-1mm}
\end{figure} 
Thereby, Step 1 of training is conducted over the interval of epochs $\text{t} \in [1 \,\,\, 2000]$ with the learning rate of $10^{-5}$ at every $t_k$. Fig.~\ref{lvb} reports the loss trajectories in Step~1 of training R-Nets in the basic mode using noiseless data. Fig.~\ref{Da} presents the network-predicted regularization maps $\alpha_{\text{\tiny \emph{NN}}}^{b_1}$ at $\text{t} = 2000$ on the low-resolution grid at all sensing steps $t_k$. Fig.~\ref{Da} also provides the misfit $\Xi^t\!$ between the approximate Morozov maps $\alpha_{\text{\tiny{M}}}^t$ used for training and the R-Net predictions via 
\beq\lb{Xi}
\Xi^t ~=~ | \alpha_{\text{\tiny \emph{NN}}}^{b_1} - \alpha_{\text{\tiny{M}}}^t |.
\eeq    
        
This plot confirms the training success and that the selected architecture for R-Nets provides sufficient complexity to capture the regularization maps at all time steps. Moreover, the validation loss trajectory in Fig.~\ref{lvb} indicates that $\text{t} = 1000$ may be a good point to stop training in Step 1, although overfitting is not observed over longer training periods. Step 2 of training is then followed over the interval of epochs $\text{t} \in [2001 \,\,\, N_{s_t}]$ with the learning rate of $5\!\times\!10^{-8}$. Here, $N_{s_t}$ signifies the epoch where the stop training flag $s_t$ turns one according to Algorithm~\ref{AL2}.
\begin{figure}[!tp]
\center\includegraphics[width=0.999\linewidth]{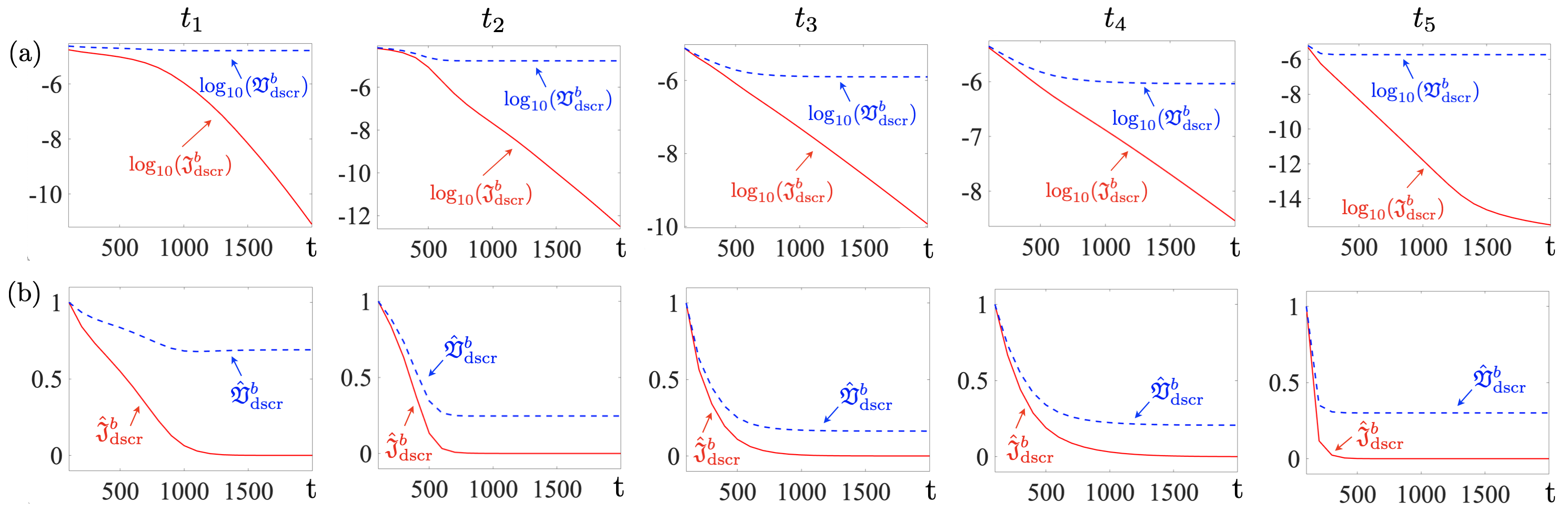} \vspace*{-7mm} 
\caption{Loss trajectories in Step~1 of training R-Nets in the basic mode using noiseless data:~(a) training loss $\log_{10}(\mathfrak{J}_{\text{dscr}}^b)$ (solid red line) and validation loss $\log_{10}(\mathfrak{V}_{\text{dscr}}^b)$ (dashed blue line) against the number of epochs $\text{\emph{t}}$ at every sensing step $t_k$, $k = 1,2,\ldots, 5$, and (b) normal training loss $\hat{\mathfrak{J}}_{\text{dscr}\!}^b$ (solid red line) and normal validation loss $\hat{\mathfrak{V}}_{\text{dscr}\!}^b$ (dashed blue line) versus $\text{\emph{t}}$.} \lb{lvb}
\vspace*{1mm}
%
\center\includegraphics[width=0.999\linewidth]{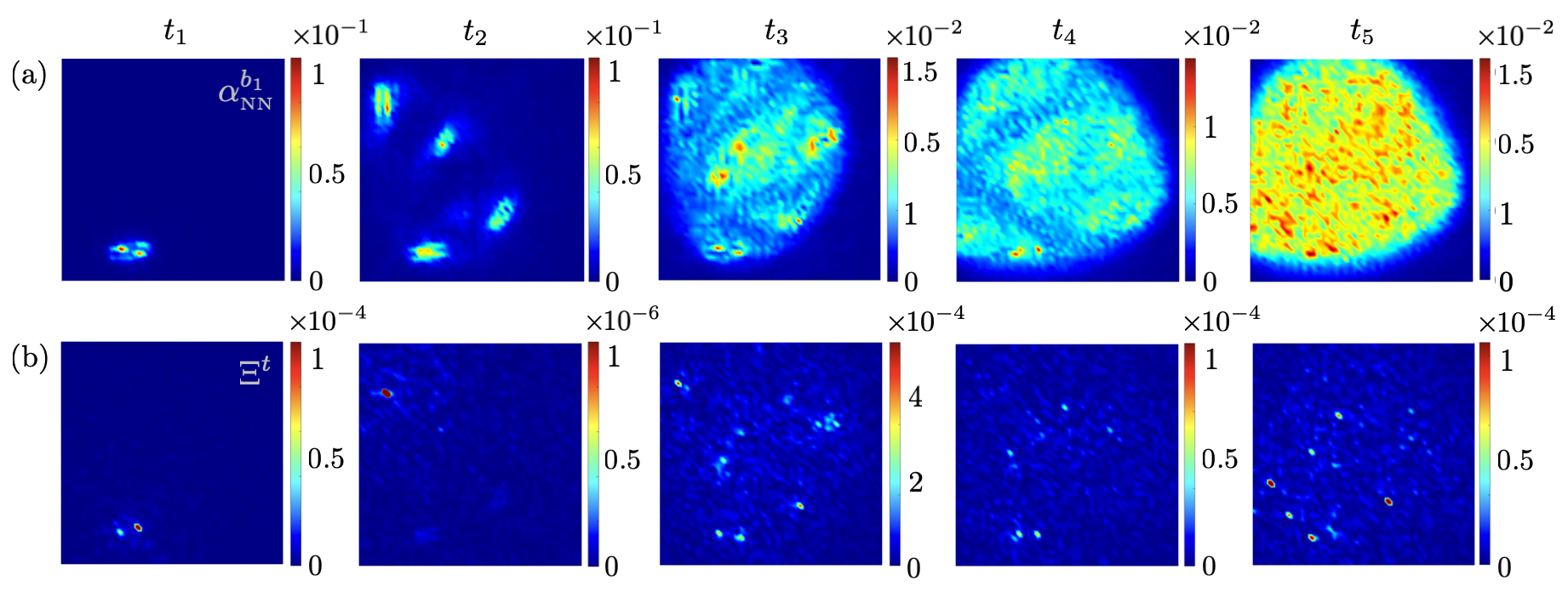} \vspace*{-7mm} 
\caption{R-Net's performance in Step 1 of training in the basic mode using noiseless data:~(a) network-predicted regularization maps $\alpha_{\text{\tiny \emph{NN}}}^{b_1}$ on the training grid (of $50 \!\times\! 50$ sampling points) at every sensing step $t_k$, $k = 1,2,\ldots, 5$, and (b) the corresponding misfits~$\Xi^t\!$, according to~\eqref{Xi}, between the R-Net's outputs and labeled Morozov maps used for training. } \lb{Da}
\vspace*{-2mm}
\end{figure} 
Fig.~\ref{bil} reports the (training and validation) loss trajectories in Step 2 of training basic R-Nets along with the trajectories of all relevant (absolute and relative) measures to the stop training criteria. Fig.~\ref{bil} highlights that, in Step 2, as $\text{t}$ increases, the validation loss shows diminishing returns which opposes the training loss behavior at sensing steps $t_3$, $t_4$, and $t_5$ which points to overfitting. To prevent the latter, the proposed criteria in Algorithm~\ref{AL2} stops saving the model as soon as 

\beq\lb{SCN}
\begin{aligned}
&\langle\Delta\hat{\mathfrak{V}}_{\text{img}}\rangle(t) < \sigma_{\text{a}} = 10^{-4} \,\,\, \lor \,\,\, \langle\Delta\hat{\mathfrak{J}}_{\text{img}}\rangle(t) < \sigma_{\text{a}} = 10^{-4}, \\*[0.25mm]
& \varrho(t) < \sigma_{\text{r}} = 5 \,\,\, \lor \,\,\, \varrho(t) > 1/\sigma_{\text{r}} = 1/5.
\end{aligned} 
 \eeq 
   
Here, $\langle\Delta\hat{\mathfrak{J}}_{\text{img}}\rangle(t)$ and $\langle\Delta\hat{\mathfrak{V}}_{\text{img}}\rangle(t)$ signify the rms of normal loss variations within $[t_\circ \,\,\, t)$ wherein $t_\circ = \max(t-N_{\text{rms}}, \text{epoch1})$ with $N_{\text{rms}} = 10000$ and $\text{epoch1} = 2000$. Note that when the domain is simple (e.g., at $t_1$) and the approximate Morozov threshold $\eta_\circ$ is close to its optimal value (see Table~\ref{MT}), Step 2 of training stops at the first epoch as the validation loss does not reflect any improvements even though the training loss decreases by $20\%$. 
\begin{figure}[!tp]
\center\includegraphics[width=0.999\linewidth]{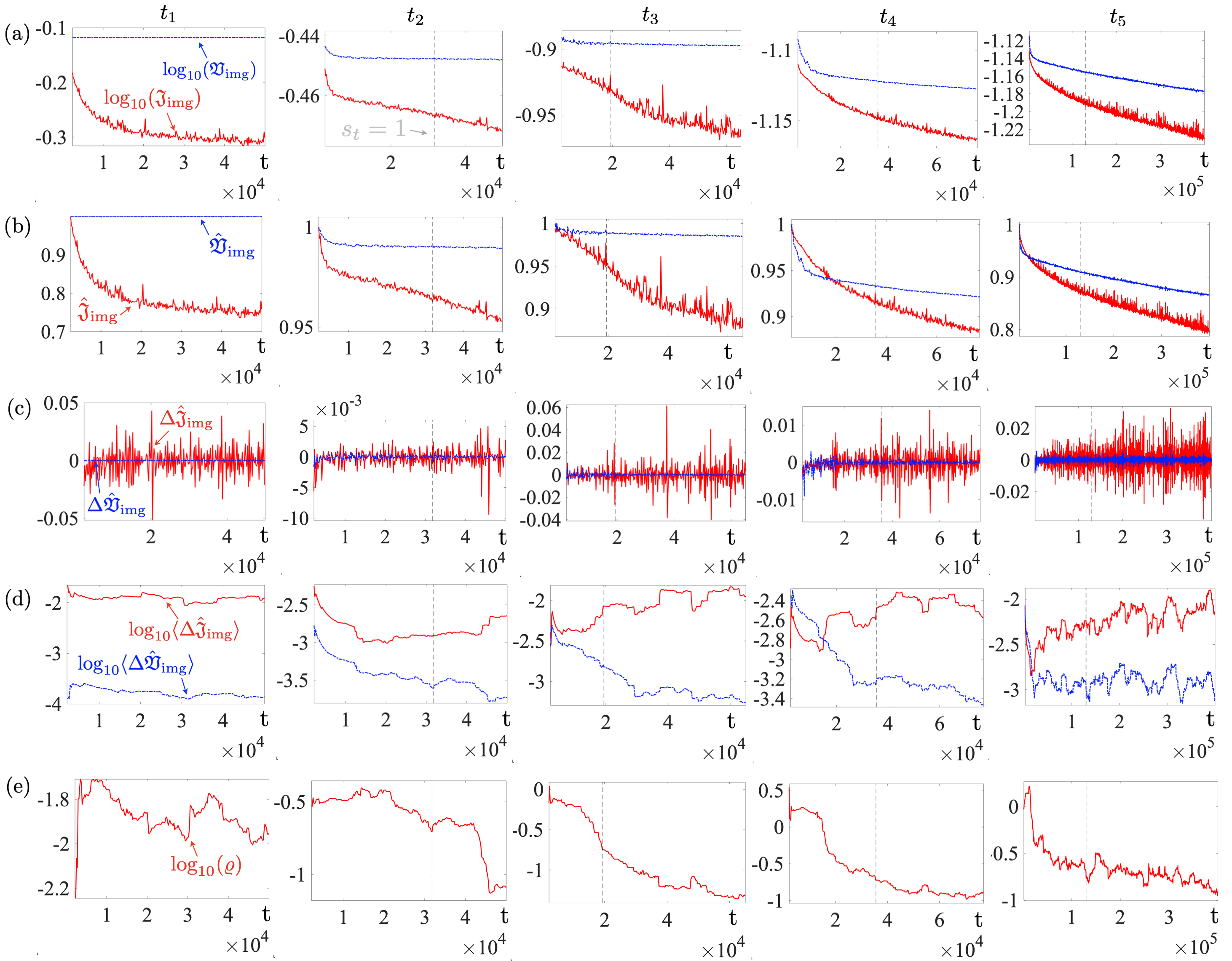} \vspace*{-7mm} 
\caption{ Convergence plots in Step~2 of training the \emph{basic R-Nets} on noiseless data:~(a) training loss $\log_{10}(\mathfrak{J}_{\text{img}})$ (solid red line) and validation loss $\log_{10}(\mathfrak{V}_{\text{img}})$ (dash-dotted blue line) against the number of epochs $\text{\emph{t}}$ at every sensing step $t_k$, $k = 1,2,\ldots, 5$, (b) normal training loss $\hat{\mathfrak{J}}_{\text{img}\!}$ (solid red line) and normal validation loss $\hat{\mathfrak{V}}_{\text{img}\!}$ (dash-dotted blue line) versus $\text{\emph{t}}$, (c) variation of normal training loss $\Delta\hat{\mathfrak{J}}_{\text{img}\!}$ (solid red line) and variation of normal validation loss $\Delta\hat{\mathfrak{V}}_{\text{img}\!}$ (dash-dotted blue line) against $\text{\emph{t}}$, (d) rms of variations of normal training loss $\log_{10} \langle \Delta\hat{\mathfrak{J}}_{\text{img}} \rangle$ (solid red line) and rms of variations of normal validation loss  $\log_{10} \langle\Delta\hat{\mathfrak{V}}_{\text{img}}\rangle$ (dash-dotted blue line) versus the number of epochs $\text{\emph{t}}$, and (e) relative loss trajectory $\log_{10} ( \varrho ) (\text{\emph{t}})$. In all panels, the vertical dashed line indicates where the stop training criteria per Algorithm~\ref{AL2} is satisfied i.e., where $s_t = 1$.} \lb{bil}
\vspace*{-2mm}
\end{figure} 
Fig.~\ref{abMNN} compares the optimal Morozov maps $\alpha_{\text{\tiny \emph{M}}}$ of Fig.~\ref{M0} on the dense search grid (of $100 \times 100$ sampling points) with the regularization maps $\alpha_{\text{\tiny \emph{NN}}}^{b_1}$ and $\alpha_{\text{\tiny \emph{NN}}}^{b_2}$ generated by basic R-Nets on the same grid at the end of training Steps 1 and 2, respectively. The affiliated LSM reconstructions are shown in Fig.~\ref{lbMNN}. It is evident from Fig.~\ref{abMNN} (b) that basic R-Nets have a mediocre generalizability. In other words, the network, whose prediction perfectly matches $\alpha_{\text{\tiny \emph{M}}}$ on the reduced training grid by the end of Step 1 (at $\text{t} = 2000$), produces a noisy regularization map on the dense search grid. The artifacts are suppressed in Step 2 of training to a limited extent. The basic R-Net's generalization errors are reflected in the LSM images of Fig.~\ref{lbMNN} in the form of elevated background fluctuations. To have a more quantitative measure of image quality, we introduce the contrast metrics. For this purpose, each LSM image is first dissected into the ``defect" and ``background" regions as depicted in Fig.~\ref{cl}. The ``background" refers to the area excluding the immediate vicinity of the reconstructed scatterers, while the ``defect" is its complement on the search grid i.e., image support. In this setting, the contrast measure ${\mathfrak{C}_{\text{\tiny mn}}}$ (\emph{resp}.~${\mathfrak{C}_{\text{\tiny mx}}}$) is defined by the ratio of the mean (\emph{resp}.~maximum) indicator value within the defect region (i.e., a neighborhood of the recovered scatterers) to the root-mean-square of the LSM imaging indicator over the background region.
\begin{figure}[!tp]
\center\includegraphics[width=0.8\linewidth]{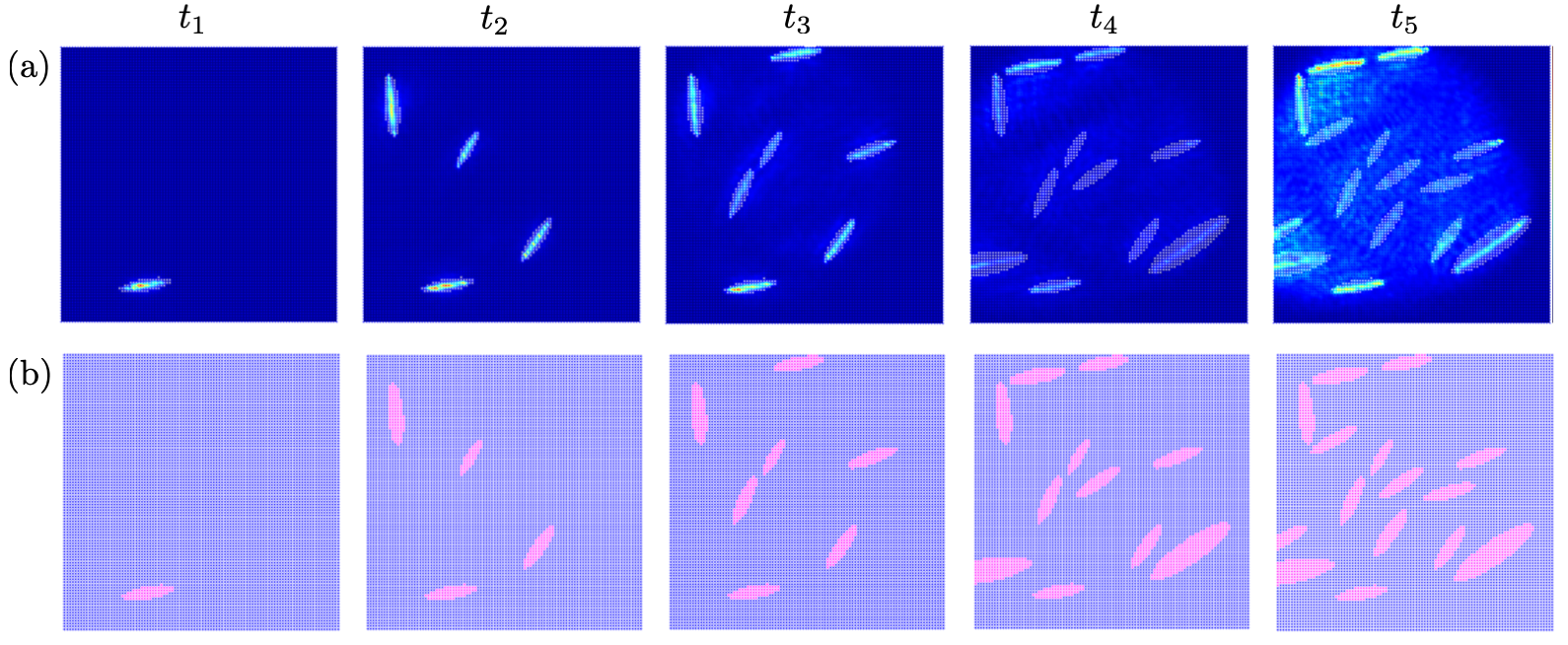} \vspace*{-2.5mm} 
\caption{Defect and background regions used to compute the contrast metrics for the reconstructions at time steps $t_k$, $k = 1,2,\ldots, 5$: (a) LSM images dissected into defect (dark blue) and background (light gray) neighborhoods, and (b) the support of defect (light purple) and background (dark blue) regions at each $t_k$.} \lb{cl}
\vspace*{-1mm}
\end{figure}
 \begin{figure}[!tp]
\center\includegraphics[width=0.999\linewidth]{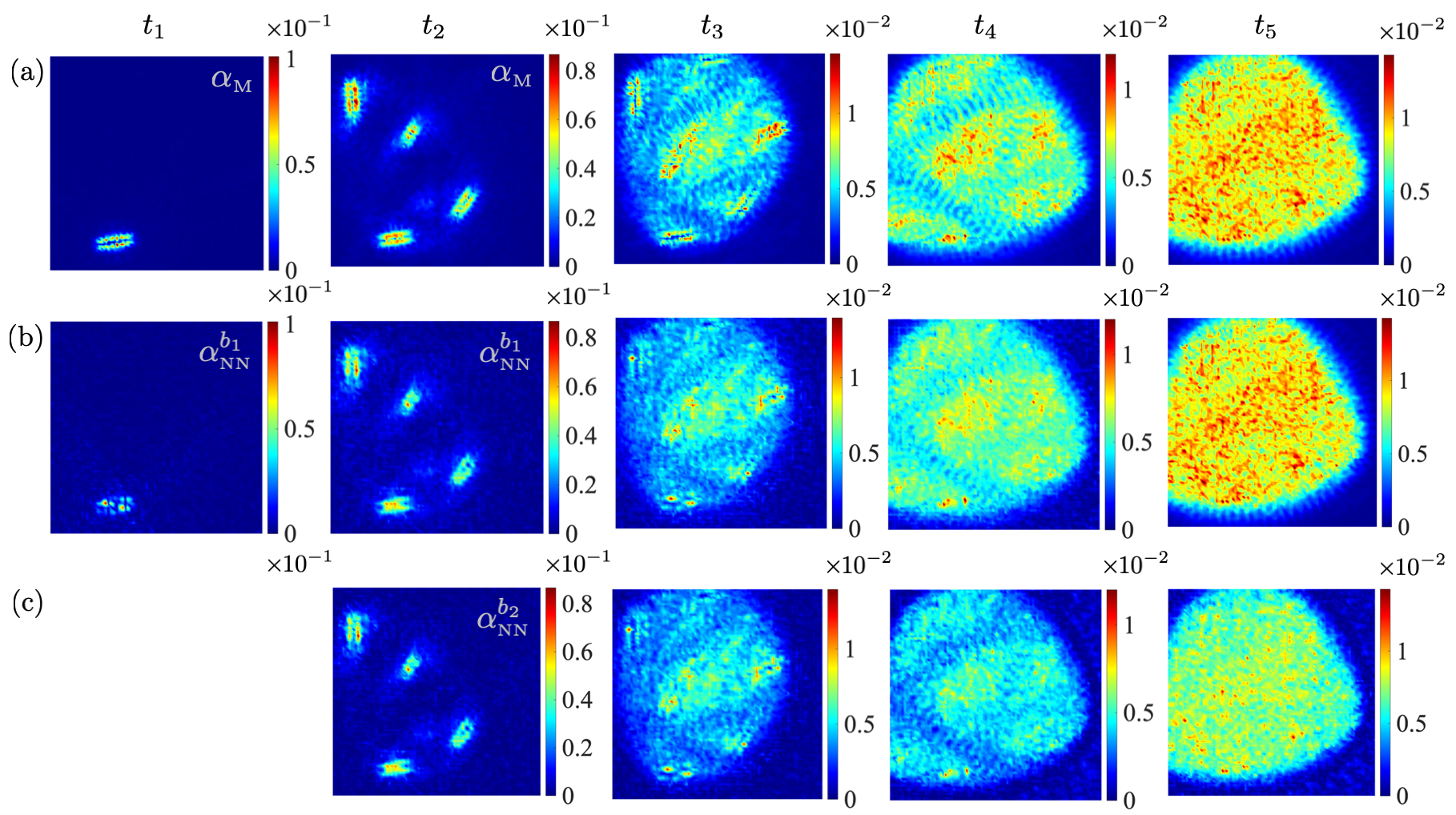} \vspace*{-7mm} 
\caption{Performance of the \emph{basic R-Nets} trained on noiseless data:~(a) manually optimized Morozov regularization maps $\alpha_{\text{\tiny \emph{M}}}$ on the dense grid (of $100 \!\times\! 100$ sampling points) at every sensing step $t_k$, $k = 1,2,\ldots, 5$, (b) network-predicted regularization maps $\alpha_{\text{\tiny \emph{NN}}}^{b_1}$ by the end of Step 1 (at $\text{\emph{t}}=2000$) on the same sampling grid, and (c) R-Net-generated regularization maps $\alpha_{\text{\tiny \emph{NN}}}^{b_2}$ by the end of Step 2 (when $s_t = 1$ as  shown in Fig.~\ref{bil}). } \lb{abMNN}
\vspace*{0mm}
\center\includegraphics[width=0.999\linewidth]{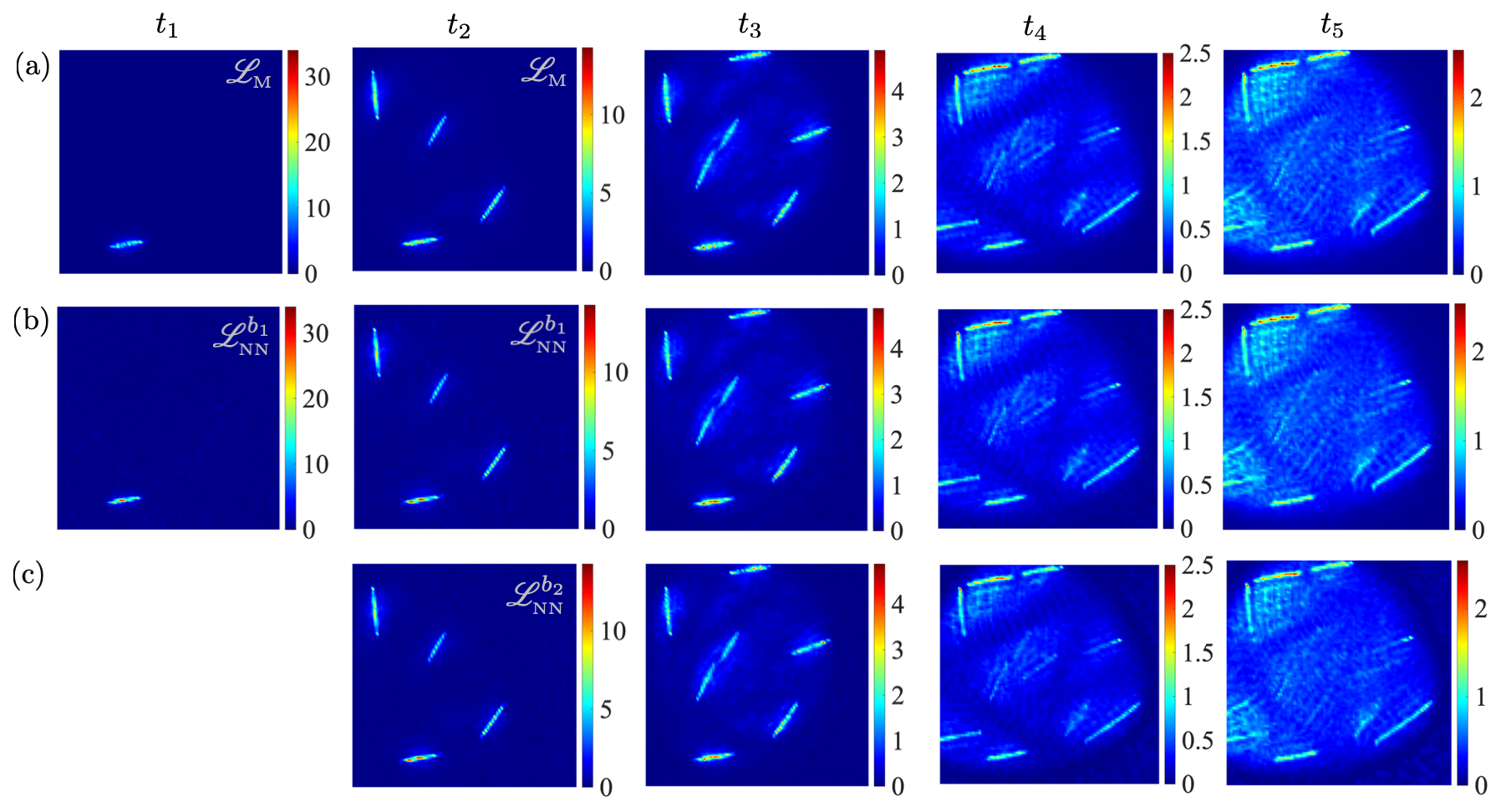} \vspace*{-7mm} 
\caption{LSM reconstructions corresponding to the regularization maps of Fig.~\ref{abMNN}:~(a) LSM images $\mathcal{L}_{\text{\tiny \emph{M}}}$ based on the manually optimized Morozov maps $\alpha_{\text{\tiny \emph{M}}}$,~(b) LSM reconstructions $\mathcal{L}_{\text{\tiny \emph{NN}}}^{b_1}$ by the network-generated regularization maps $\alpha_{\text{\tiny \emph{NN}}}^{b_1}$ at the end of Step 1 (at $\text{\emph{t}}=2000$), and~(c) LSM indicator $\mathcal{L}_{\text{\tiny \emph{NN}}}^{b_2}$ via the R-Net regularization maps $\alpha_{\text{\tiny \emph{NN}}}^{b_2}$ at the end of training (when $s_t = 1$ in Step 2).    } \lb{lbMNN}
\vspace*{-1mm}
\end{figure}
This definition may be interpreted as an approximation of the Weber's contrast~\cite{avat2023,otsu1975},
\[
W_c ~=~ \frac{I-I_b}{I_b},
\]
wherein $I$ represents the intensity of the object, while $I_b$ is the background intensity. On assuming $I_b <\!\!< I$, $W_c \sim {I}/{I_b}$. Now, if $I_b$ is specified by the rms of indicator values in the background region and $I$ represents the indicator value at every sampling point in the defect region, then the adopted contrast metrics may be described by ${\mathfrak{C}_{\text{\tiny mn}}} = \text{mean}({I}/{I_b})$ and ${\mathfrak{C}_{\text{\tiny mx}}} = \norms{\!{I}/{I_b}\!}_{\infty}$. This metric is intuitive and is inspired by the human visual system's relative sense of contrast~\cite{avat2023,otsu1975}. A higher value of the contrast metric means that the object is more distinguishable from its background. rms is a common statistical metric for quantifying signal fluctuations relative to background noise in imaging applications~\cite{zhang2018,szab2013}.

Table~\ref{T_lbMNN} provides the contrast metrics germane to Fig.~\ref{lbMNN} which points to similar conclusions mentioned above. More specifically, R-Nets trained by Step 1 generate regularization maps that lead to images with lower contrast compared to the ones based on the discrepancy principle. Step 2 recovers the image quality in terms of both (mean/max) contrast metrics. This improvement is however is neither significant nor consistent across all sensing steps.   

\subsubsection*{Reconstructions via discrepancy-informed R-Nets}

Discrepancy-informed R-Nets are trained, according to Algorithm~\ref{AL1}, by minimizing~\eqref{Jids}-\eqref{J2n} and~\eqref{w12n} in Step 1 and~\eqref{Jimg} in Step 2. 
 \begin{table}[!tp]
\begin{center}
\caption{\small{Contrast metric computed for the LSM reconstructions in Fig.~\ref{lbMNN} germane to \emph{basic R-Nets} trained on noiseless data.}} \vspace*{-2mm}
\label{T_lbMNN}
\begin{tabular}{c||c|c|c|c|c} 
\diagbox{{\small{contrast}}}{{\!\small{time}}} & $t_1$\! & $t_2$\! & $t_3$\! & $t_4$\! & $t_5$\! \\ \hline\hline  
$\mathfrak{C}_{\text{\tiny mn}}(\mathcal{L}_\text{\tiny M})$    & $99.94$  & $30.85$  & $10.41$ & $3.69$ & $2.49$ \\  
 \hline
$\mathfrak{C}_{\text{\tiny mn}}(\mathcal{L}_\text{\tiny NN}^{b_1})$    & \!\!$69.28$\!\!  & \!\!$27.33$\!\!  & \!\!$11.02$\!\! & \!\!$3.57$\!\! & \!\!$2.42$\!\! \\ 
 \hline
$\mathfrak{C}_{\text{\tiny mn}}(\mathcal{L}_\text{\tiny NN}^{b_2})$   & \!\!$69.28$\!\!  & \!\!$30.18$\!\! & \!\!$11.15$\!\! & \!\!$3.72$\!\! & \!\!$2.47$\!\! \\
 \hline 
$\mathfrak{C}_{\text{\tiny mx}}(\mathcal{L}_\text{\tiny M})$    & $166.46$  & $54.58$  & $24.64$ & $10.40$ & $7.52$ \\  
 \hline
$\mathfrak{C}_{\text{\tiny mx}}(\mathcal{L}_\text{\tiny NN}^{b_1})$    & \!\!$124.47$\!\!  & \!\!$54.58$\!\!  & \!\!$24.84$\!\! & \!\!$11.53$\!\! & \!\!$7.42$\!\! \\ 
 \hline
$\mathfrak{C}_{\text{\tiny mx}}(\mathcal{L}_\text{\tiny NN}^{b_2})$   & \!\!$124.47$\!\!  & \!\!$61.09$\!\! & \!\!$24.40$\!\! & \!\!$12.07$\!\! & \!\!$7.40$\!\! \\
\end{tabular}
\end{center}
\vspace*{-1mm}
\end{table} 
\begin{figure}[!h]
\center\includegraphics[width=0.999\linewidth]{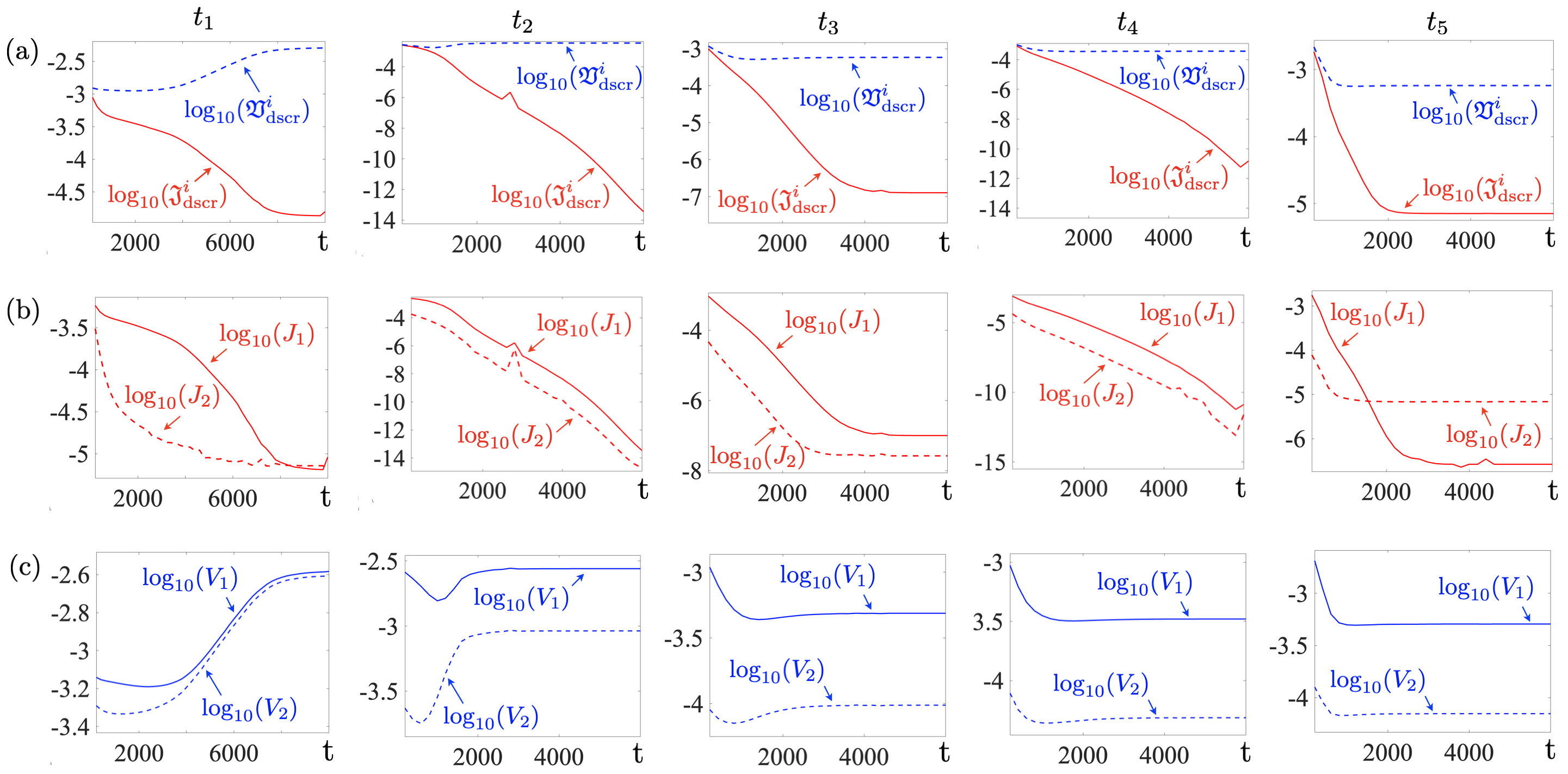} \vspace*{-7mm} 
\caption{Loss trajectories in Step~1 of training the \emph{discrepancy-informed R-Nets} on noiseless data:~(a) total training loss $\log_{10}(\mathfrak{J}_{\text{dscr}}^i)$ (solid red line) and total validation loss $\log_{10}(\mathfrak{V}_{\text{dscr}}^i)$ (dashed blue line) against the number of epochs $\text{\emph{t}}$ at every sensing step $t_k$, $k = 1,2,\ldots, 5$, (b) training loss components $\log_{10}({J}_{1})$ (solid red line) and $\log_{10}({J}_{2})$ (dashed blue line) versus $\text{\emph{t}}$ according to~\eqref{Jids2}, and (c) validation loss components $\log_{10}({V}_{1})$ (solid red line) and $\log_{10}({V}_{2})$ (dashed blue line) versus $\text{\emph{t}}$. } \lb{lvi}
\vspace*{-1mm}
\end{figure} 
\renewcommand{\arraystretch}{1.25}
The training and validation datasets are similar to that of basic R-Nets. The only difference here is the learning objective in Step 1 where the logic of discrepancy principle is included in the loss function. Given that Step 1 in this case involves multitask optimization, the training is conducted over the interval of epochs $\text{t} \in [1 \,\,\, 1000]$ with the learning rate of $5 \times 10^{-6}$ for all $t_k$, $k = 1,2,\dots,5$. Fig.~\ref{lvi} reports the convergence plots in Step~1 of training R-Nets in the informed mode using noiseless data. The validation loss trajectories in early sensing steps $t_1$ and $t_2$ indicate a slight overfitting that may be due to the fact that the R-Net model is more complex than required at these time steps due to the simplicity of the domain. However, since in general the geometry of the domain (under inspection) is unknown, and in order to remain consistent in our comparative analysis, the R-Net architecture remains unaltered at all time steps. It should be noted that the loss components $J_1$, $J_2$, $V_1$, and $V_2$ in Fig.~\ref{lvi} are defined by
\beq\lb{Jids2}  
\begin{aligned}
& J_1 = \frac{1}{{N_{\text{trn}}}}\sum_{\text{t} = 1}^{N_{\text{trn}}} w_1^{\text{t}} J_1^{\xxs\text{t}}, \quad J_2 = \frac{1}{{N_{\text{trn}}}}\sum_{\text{t} = 1}^{N_{\text{trn}}} w_2^{\text{t}} J_2^{\xxs\text{t}}, \\*[0.25mm]
& V_1 = \frac{1}{{N_{\text{vld}}}}\sum_{\text{v} = 1}^{N_{\text{vld}}} w_1^{\text{v}} J_1^{\xxs\text{v}}, \quad V_2 = \frac{1}{{N_{\text{vld}}}}\sum_{\text{v} = 1}^{N_{\text{vld}}} w_2^{\text{v}} V_2^{\xxs\text{v}}.
\end{aligned}
\eeq  

It appears from Fig.~\ref{lvi}~(b) that the proposed weights in~\eqref{w12n} effectively balance the loss components. The distribution of loss weights over the training grid are provided at $\text{t} = 10000$ in Figs.~\ref{w1t} and~\ref{w2t}. 
\begin{figure}[!tp]
\center\includegraphics[width=0.999\linewidth]{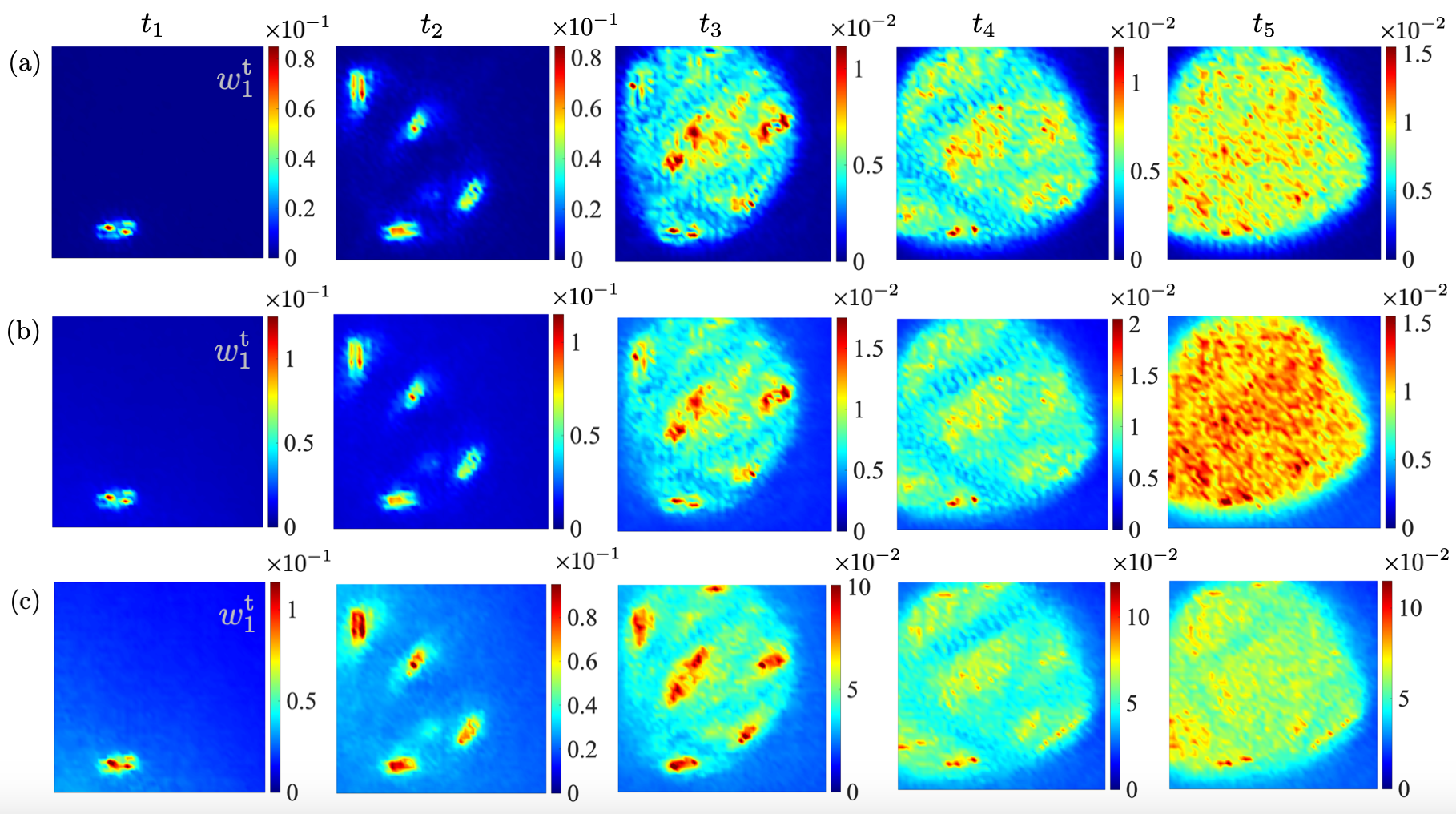} \vspace*{-7mm} 
\caption{Spatial distribution of loss weight $w_1^{\text{t}}$ in~\eqref{w12n} used for training \emph{discrepancy-informed R-Nets} on a reduced grid (of $50 \!\times\! 50$ sampling points) when the training is conducted using:~(a) noiseless data,~(b) $10\%$ noisy data, and~(c) $25\%$ noisy data. } \lb{w1t}
\vspace*{0mm}
\center\includegraphics[width=0.999\linewidth]{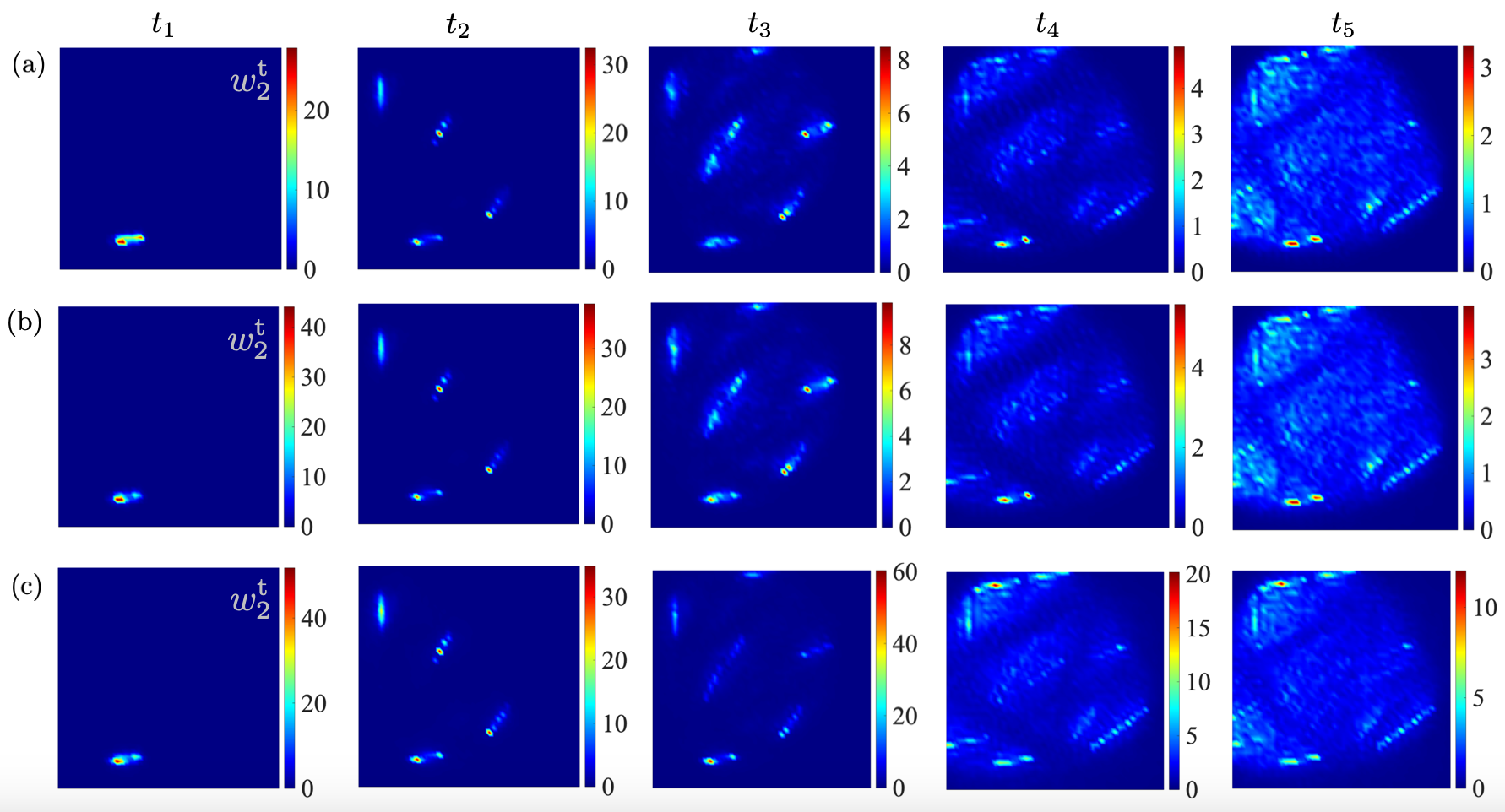} \vspace*{-7mm} 
\caption{Spatial distribution of loss weight $w_2^{\text{t}}$ in~\eqref{w12n} used for training \emph{discrepancy-informed R-Nets} on a reduced grid (of $50 \!\times\! 50$ sampling points) when the training is conducted using:~(a) noiseless data,~(b) $10\%$ noisy data, and~(c) $25\%$ noisy data. } \lb{w2t}
\vspace*{-1mm}
\end{figure} 
Step 2 of training is then followed over the interval of epochs $\text{t} \in [10001 \,\,\, N_{s_t}]$ with the learning rate of $5\!\times\!10^{-8}$. Fig.~\ref{iil} presents the convergence plots in Step 2 of training informed R-Nets on noiseless data. The stop training criteria is similar to that of basic R-Nets according to~\eqref{SCN}.
\begin{figure}[!tp]
\center\includegraphics[width=0.999\linewidth]{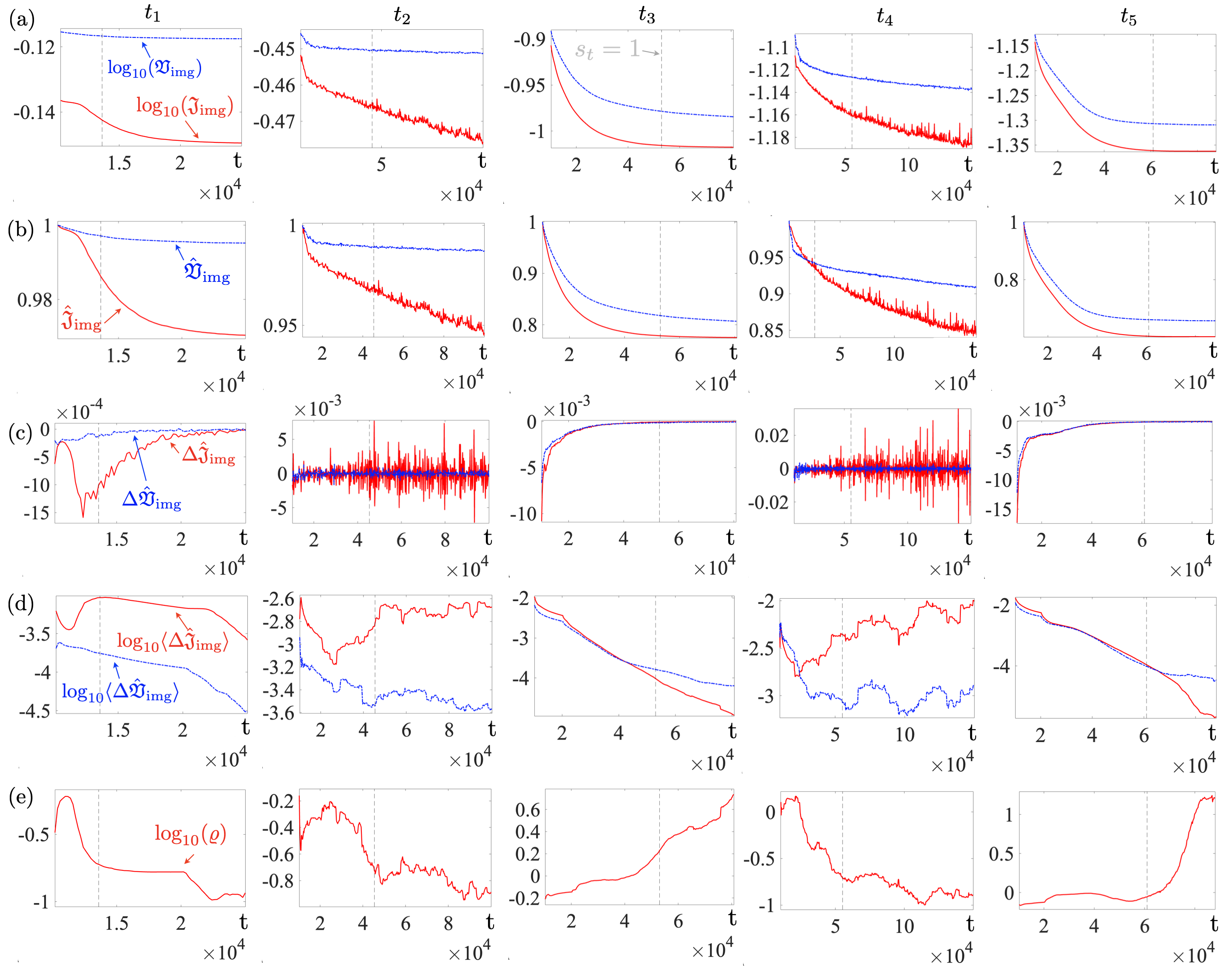} \vspace*{-7mm} 
\caption{Loss trajectories in Step~2 of training the \emph{discrepancy-informed R-Nets} on noiseless data:~(a) training loss $\log_{10}(\mathfrak{J}_{\text{img}})$ (solid red line) and validation loss $\log_{10}(\mathfrak{V}_{\text{img}})$ (dash-dotted blue line) against the number of epochs $\text{\emph{t}}$ at every sensing step $t_k$, $k = 1,2,\ldots, 5$, (b) normal training loss $\hat{\mathfrak{J}}_{\text{img}\!}$ (solid red line) and normal validation loss $\hat{\mathfrak{V}}_{\text{img}\!}$ (dash-dotted blue line) versus $\text{\emph{t}}$, (c) variation of normal training loss $\Delta\hat{\mathfrak{J}}_{\text{img}\!}$ (solid red line) and variation of normal validation loss $\Delta\hat{\mathfrak{V}}_{\text{img}\!}$ (dash-dotted blue line) against $\text{\emph{t}}$, (d) rms of variations of normal training loss $\log_{10} \langle \Delta\hat{\mathfrak{J}}_{\text{img}} \rangle$ (solid red line) and rms of variations of normal validation loss  $\log_{10} \langle\Delta\hat{\mathfrak{V}}_{\text{img}}\rangle$ (dash-dotted blue line) versus the number of epochs $\text{\emph{t}}$, and (e) relative loss trajectory $\log_{10} ( \varrho ) (\text{\emph{t}})$. In all panels, the vertical dashed line indicates where the stop training criteria per Algorithm~\ref{AL2} is satisfied i.e., where $s_t = 1$.} \lb{iil}
\vspace*{-1mm}
\end{figure} 
Fig.~\ref{aiMNN} provides the comparison among the optimal Morozov maps $\alpha_{\text{\tiny \emph{M}}}$ of Fig.~\ref{M0} and the regularization maps $\alpha_{\text{\tiny \emph{NN}}}^{i_1}$ and $\alpha_{\text{\tiny \emph{NN}}}^{i_2}$ furnished by informed R-Nets in Steps 1 and 2 of training, respectively.
\begin{figure}[!tp]
\center\includegraphics[width=0.999\linewidth]{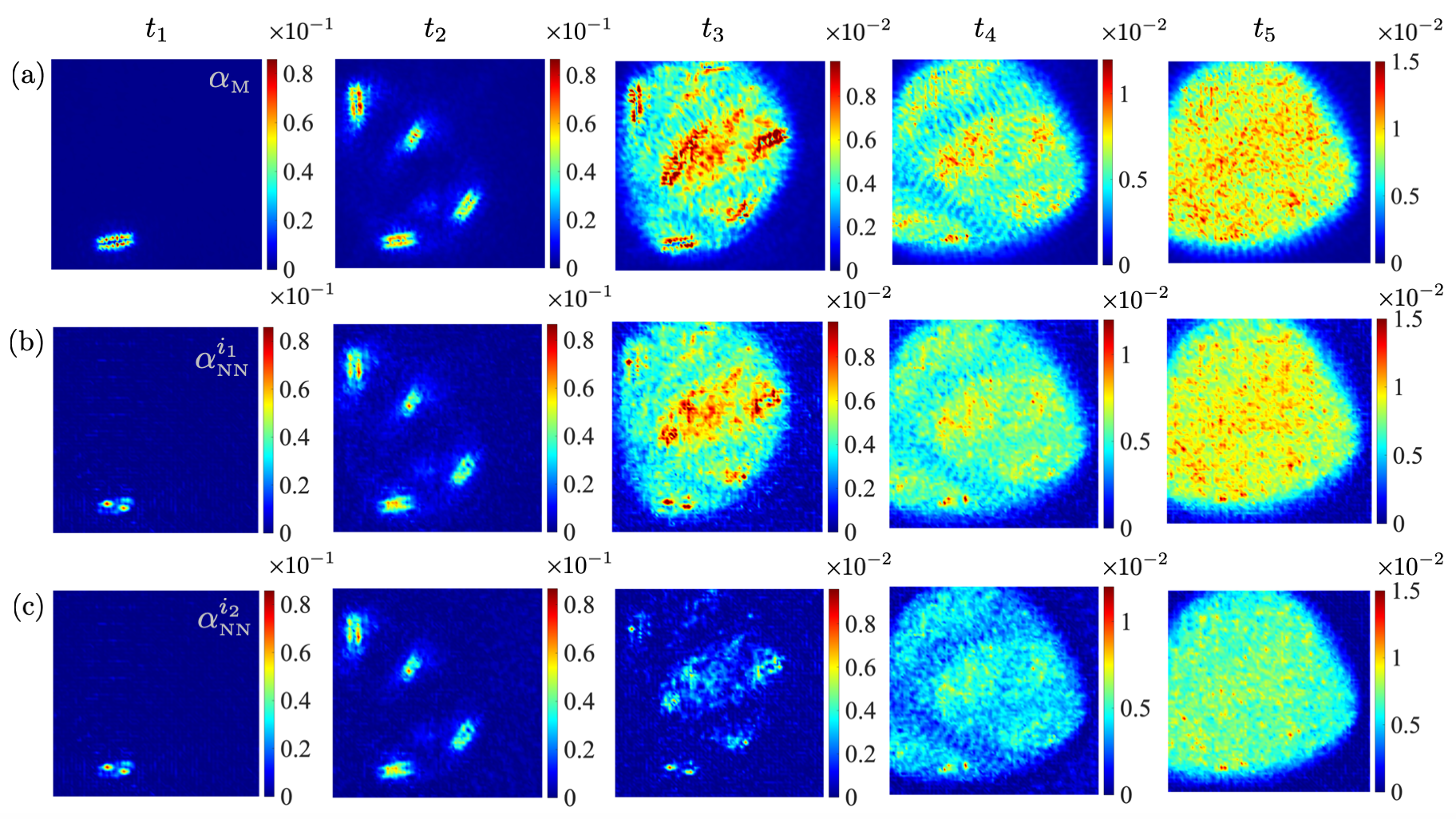} \vspace*{-7mm} 
\caption{Performance of the \emph{discrepancy-informed R-Nets} trained on noiseless data:~(a) manually optimized Morozov regularization maps $\alpha_{\text{\tiny \emph{M}}}$ on the dense grid (of $100 \!\times\! 100$ sampling points) at every sensing step $t_k$, $k = 1,2,\ldots, 5$, (b) network-predicted regularization maps $\alpha_{\text{\tiny \emph{NN}}}^{i_1}$ by the end of Step 1 (at $\text{\emph{t}}=10^4$) on the same sampling grid, and (c) R-Net-generated regularization maps $\alpha_{\text{\tiny \emph{NN}}}^{i_2}$ by the end of Step 2 (when $s_t = 1$ as shown in Fig.~\ref{iil}). Maps corresponding to each sensing step are plotted on the same color scale.} \lb{aiMNN}
\vspace*{2mm}
\center\includegraphics[width=0.999\linewidth]{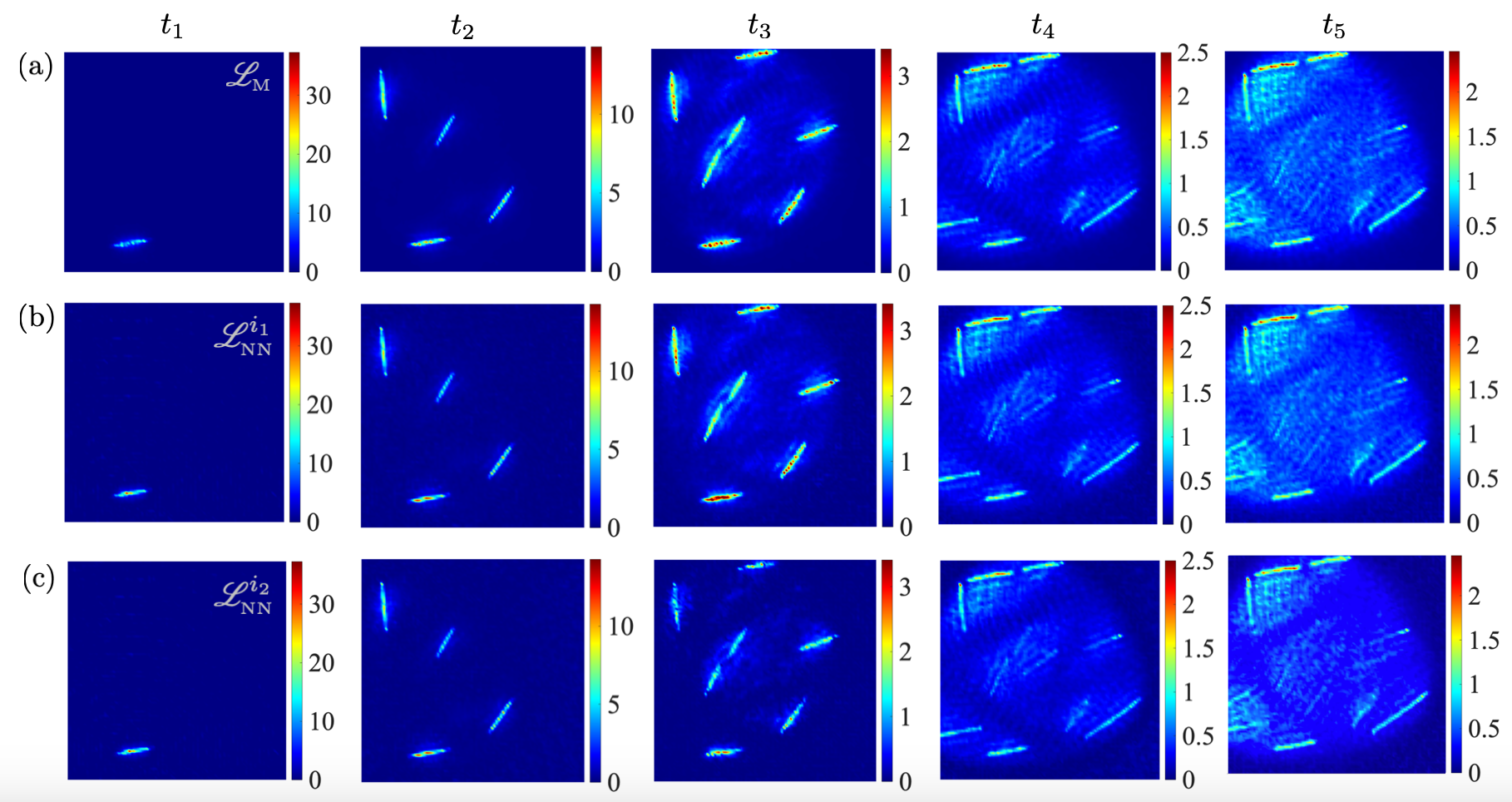} \vspace*{-7mm} 
\caption{LSM reconstructions corresponding to the regularization maps of Fig.~\ref{aiMNN}:~(a) LSM images $\mathcal{L}_{\text{\tiny \emph{M}}}$ based on the manually optimized Morozov maps $\alpha_{\text{\tiny \emph{M}}}$,~(b) LSM reconstructions $\mathcal{L}_{\text{\tiny \emph{NN}}}^{i_1}$ by the network-generated regularization maps $\alpha_{\text{\tiny \emph{NN}}}^{i_1}$ at the end of Step 1 (at $\text{\emph{t}}=10^4$), and~(c) LSM indicator $\mathcal{L}_{\text{\tiny \emph{NN}}}^{i_2}$ via the R-Net regularization maps $\alpha_{\text{\tiny \emph{NN}}}^{i_2}$ at the end of training (when $s_t = 1$ in Step 2). Maps corresponding to each sensing step are plotted on the same color scale.} \lb{liMNN}
\vspace*{-1mm}
\end{figure} 
The affiliated LSM reconstructions are shown in Fig.~\ref{liMNN}. Table~\ref{T_liMNN} reports the contrast metrics corresponding to Fig.~\ref{liMNN}. Observe from Fig.~\ref{aiMNN} and Fig.~\ref{liMNN} that (a) discrepancy-informed R-Nets exhibit superior generalization performance compared to basic R-Nets, and (b) Step 2 of training leads to a meaningful enhancement of image quality both in terms of contrast and reducing artifacts, in particular, at sensing steps $\geqslant t_3$. By comparing Tables~\ref{T_lbMNN} and~\ref{T_liMNN}, one may note that images created by informed R-Nets have consistently higher contrasts than those generated by basic R-Nets across all the reconstructions. Additionally, by the end of Step 2, informed R-Nets lead to images with marked improvement compared to the reconstructions based on the optimized Morozov maps.    
\renewcommand{\arraystretch}{1.25}
 \begin{table}[!tp]
\begin{center}
\caption{\small Contrast metric computed for the LSM reconstructions in Fig.~\ref{liMNN} corresponding to \emph{discrepancy-informed R-Nets} trained on noiseless data.} \vspace*{-1.5mm}
\label{T_liMNN}
\begin{tabular}{c||c|c|c|c|c} 
\diagbox{{\small{contrast}}}{{\!\small{time}}} & $t_1$\! & $t_2$\! & $t_3$\! & $t_4$\! & $t_5$\! \\ \hline\hline  
$\mathfrak{C}_{\text{\tiny mn}}(\mathcal{L}_\text{\tiny M})$    & $99.94$  & $30.85$  & $10.41$ & $3.69$ & $2.49$ \\  
 \hline
$\mathfrak{C}_{\text{\tiny mn}}(\mathcal{L}_\text{\tiny NN}^{i_1})$    & \!\!$95.28$\!\!  & \!\!$27.95$\!\!  & \!\!$11.22$\!\! & \!\!$3.58$\!\! & \!\!$2.43$\!\! \\ 
 \hline
$\mathfrak{C}_{\text{\tiny mn}}(\mathcal{L}_\text{\tiny NN}^{i_2})$   & \!\!$97.01$\!\!  & \!\!$31.05$\!\! & \!\!$16.12$\!\! & \!\!$3.78$\!\! & \!\!$3.21$\!\! \\
 \hline
$\mathfrak{C}_{\text{\tiny mx}}(\mathcal{L}_\text{\tiny M})$    & $166.46$  & $54.58$  & $24.64$ & $10.40$ & $7.52$ \\ 
 \hline
$\mathfrak{C}_{\text{\tiny mx}}(\mathcal{L}_\text{\tiny NN}^{i_1})$    & \!\!$189.21$\!\!  & \!\!$55.84$\!\!  & \!\!$26.04$\!\! & \!\!$11.89$\!\! & \!\!$7.16$\!\! \\ 
 \hline
$\mathfrak{C}_{\text{\tiny mx}}(\mathcal{L}_\text{\tiny NN}^{i_2})$   & \!\!$193.42$\!\!  & \!\!$63.62$\!\! & \!\!$48.26$\!\! & \!\!$12.75$\!\! & \!\!$9.22$\!\! \\
\end{tabular}
\end{center}
\vspace*{-3mm}
\end{table}
The above observations remain consistent when discrepancy-informed R-Nets are trained on noisy data. Figs.~\ref{lvi10pn} and~\ref{lvi25pn} demonstrate the loss trajectories of Step 1 of training at $10\%$ and $25\%$ noise levels. Figs.~\ref{iil10pn} and~\ref{iil25pn} illustrate the associated convergence plots related to Step 2 of training. The importance of Algorithm~\ref{AL2} to stop training in Step 2 is more evident in the case of learning from noisy data. Fig.~\ref{aiMNN10pn} compares the optimal Morozov maps $\alpha_{\text{\tiny \emph{M}}}$ of Fig.~\ref{M10pn} and the regularization maps $\alpha_{\text{\tiny \emph{NN}}}^{i_1}$ and $\alpha_{\text{\tiny \emph{NN}}}^{i_2}$ furnished by informed R-Nets trained on $10\%$ noisy data. The affiliated LSM images are provided in Fig.~\ref{liMNN10pn}. In parallel, Fig.~\ref{aiMNN25pn} compares the optimal Morozov maps $\alpha_{\text{\tiny \emph{M}}}$ of Fig.~\ref{M25pn} and the R-Net outputs $\alpha_{\text{\tiny \emph{NN}}}^{i_1}$ and $\alpha_{\text{\tiny \emph{NN}}}^{i_2}$ when trained on $25\%$ noisy data. The corresponding images are depicted in Fig.~\ref{liMNN25pn}. Tables~\ref{T_liMNN10pn} and~\ref{T_liMNN25pn} report the associated contrast metrics at $10\%$ and $25\%$ noise levels, respectively.  The main findings from these results are as follows: (1) discrepancy-informed R-Nets trained by Step 1 can effectively capture the logic of Morozov discrepancy principle and lead to reconstructions of similar or better quality, and (2) Step 2 consistently improves upon the discrepancy-based reconstructions of Step 1 which is more pronounced in the reconstructions from noisy data. To provide a quantitative comparison of the computational cost, Table~\ref{computational_cost} lists the average training time per epoch for different training steps on an Apple M2 Max processor.

\begin{figure}[!tp]
\center\includegraphics[width=0.999\linewidth]{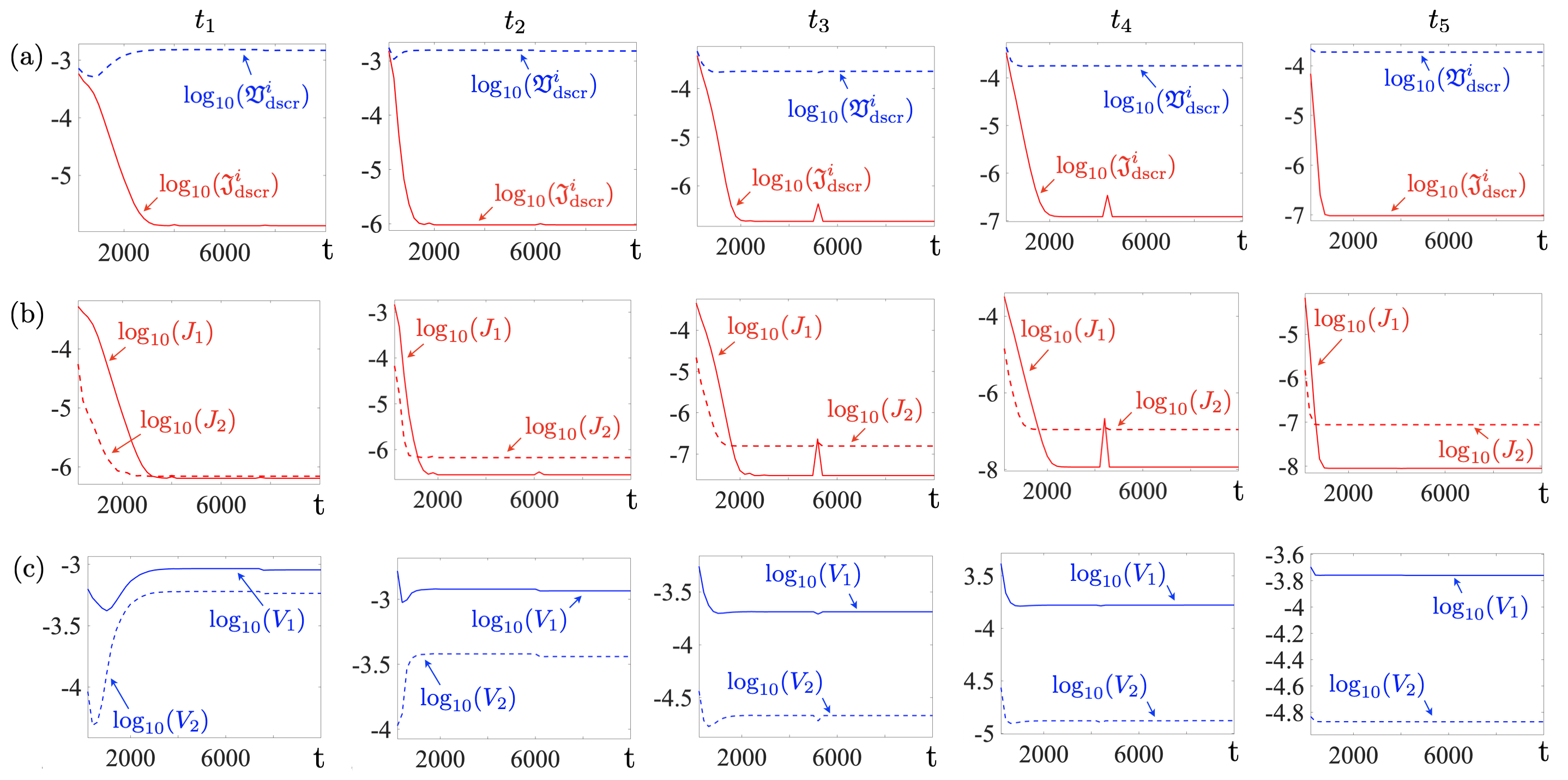} \vspace*{-7mm} 
\caption{Loss trajectories in Step~1 of training the \emph{discrepancy-informed R-Nets} on $10\%$ noisy data:~(a) total training loss $\log_{10}(\mathfrak{J}_{\text{dscr}}^i)$ (solid red line) and total validation loss $\log_{10}(\mathfrak{V}_{\text{dscr}}^i)$ (dashed blue line) against the number of epochs $\text{\emph{t}}$ at every sensing step $t_k$, $k = 1,2,\ldots, 5$, (b) training loss components $\log_{10}({J}_{1})$ (solid red line) and $\log_{10}({J}_{2})$ (dashed blue line) versus $\text{\emph{t}}$, and (c) validation loss components $\log_{10}({V}_{1})$ (solid red line) and $\log_{10}({V}_{2})$ (dashed blue line) versus $\text{\emph{t}}$. } \lb{lvi10pn}
\vspace*{-2mm}
\end{figure} 

 \begin{figure}[!tp]
\center\includegraphics[width=0.999\linewidth]{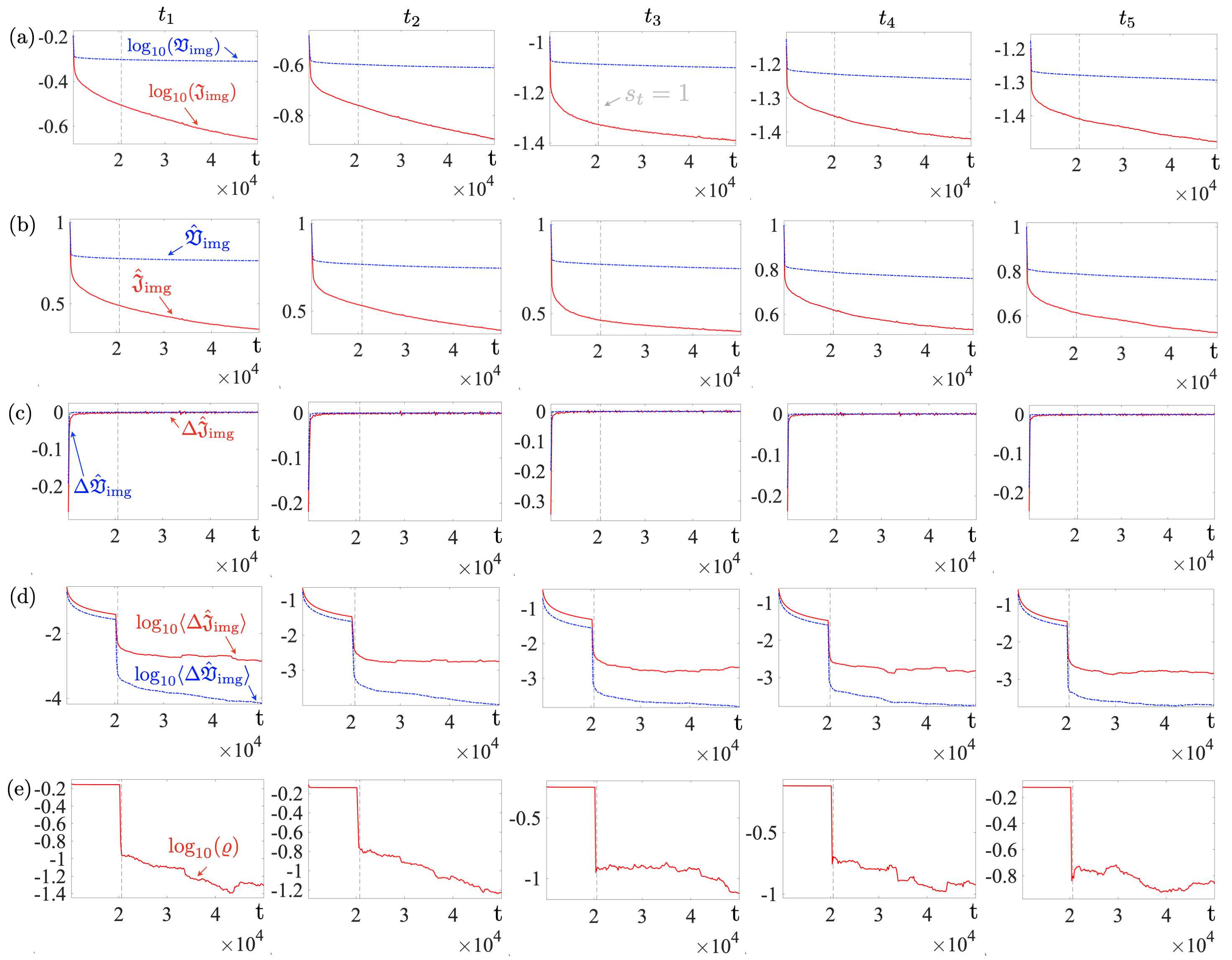} \vspace*{-7mm} 
\caption{Loss trajectories in Step~2 of training the \emph{discrepancy-informed R-Nets} on $10\%$ noisy data:~(a) training loss $\log_{10}(\mathfrak{J}_{\text{img}})$ (solid red line) and validation loss $\log_{10}(\mathfrak{V}_{\text{img}})$ (dash-dotted blue line) against the number of epochs $\text{\emph{t}}$ at every sensing step $t_k$, $k = 1,2,\ldots, 5$, (b) normal training loss $\hat{\mathfrak{J}}_{\text{img}\!}$ (solid red line) and normal validation loss $\hat{\mathfrak{V}}_{\text{img}\!}$ (dash-dotted blue line) versus $\text{\emph{t}}$, (c) variation of normal training loss $\Delta\hat{\mathfrak{J}}_{\text{img}\!}$ (solid red line) and variation of normal validation loss $\Delta\hat{\mathfrak{V}}_{\text{img}\!}$ (dash-dotted blue line) against $\text{\emph{t}}$, (d) rms of variations of normal training loss $\log_{10} \langle \Delta\hat{\mathfrak{J}}_{\text{img}} \rangle$ (solid red line) and rms of variations of normal validation loss  $\log_{10} \langle\Delta\hat{\mathfrak{V}}_{\text{img}}\rangle$ (dash-dotted blue line) versus the number of epochs $\text{\emph{t}}$, and (e) relative loss trajectory $\log_{10} ( \varrho ) (\text{\emph{t}})$. In all panels, the vertical dashed line indicates where the stop training criteria per Algorithm~\ref{AL2} is satisfied i.e., where $s_t = 1$.} \lb{iil10pn}
\vspace*{-2mm}
\end{figure} 

 \begin{figure}[!tp]
\center\includegraphics[width=0.999\linewidth]{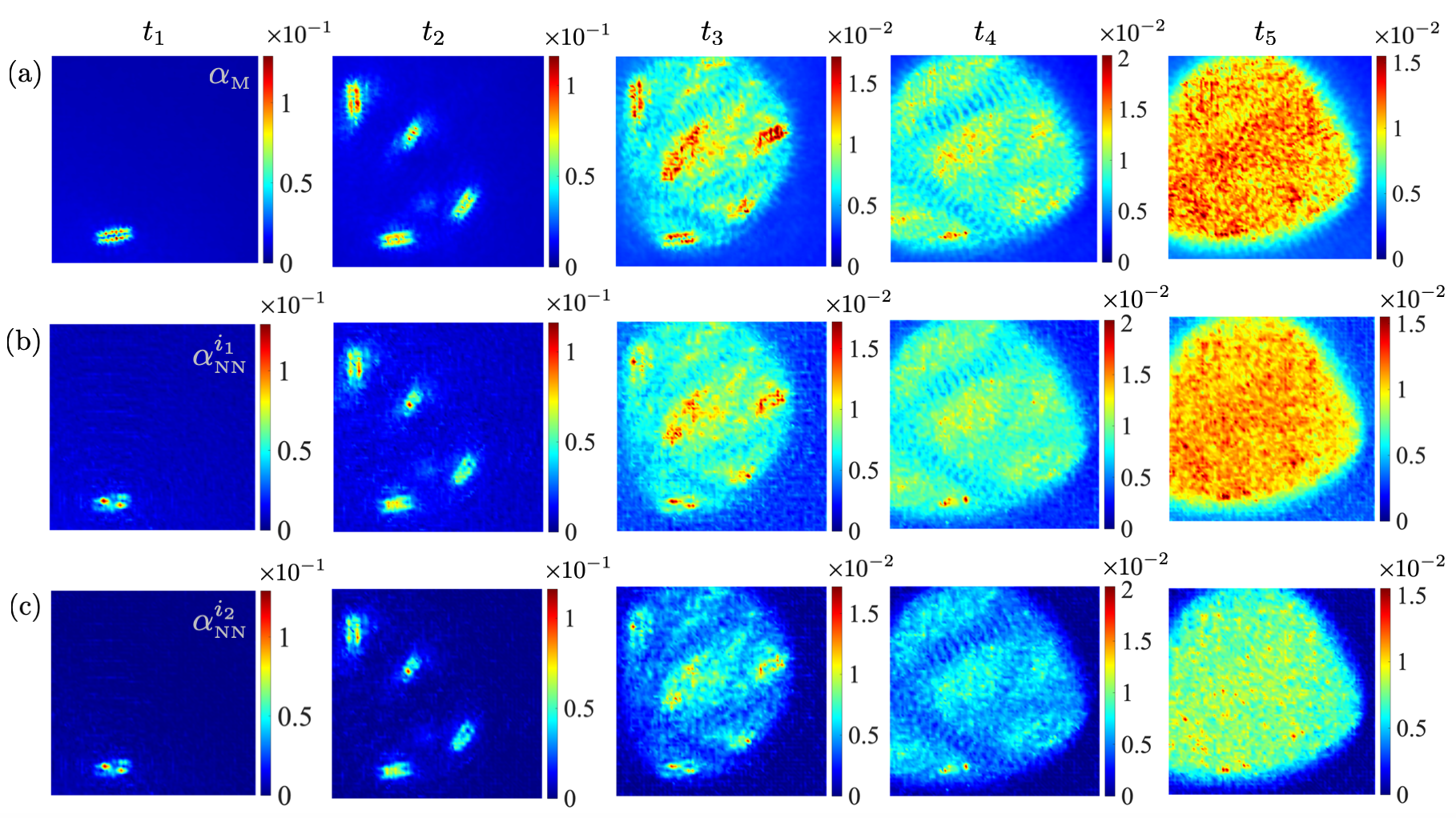} \vspace*{-7mm} 
\caption{Performance of the \emph{discrepancy-informed R-Nets} trained on $10\%$ noisy data:~(a) manually optimized Morozov regularization maps $\alpha_{\text{\tiny \emph{M}}}$ on the dense grid (of $100 \!\times\! 100$ sampling points) at every sensing step $t_k$, $k = 1,2,\ldots, 5$, (b) network-predicted regularization maps $\alpha_{\text{\tiny \emph{NN}}}^{i_1}$ by the end of Step 1 (at $\text{\emph{t}}=10^4$) on the same sampling grid, and (c) R-Net-generated regularization maps $\alpha_{\text{\tiny \emph{NN}}}^{i_2}$ by the end of Step 2 (when $s_t = 1$ as shown in Fig.~\ref{iil10pn}). Maps corresponding to each sensing step are plotted on the same color scale. } \lb{aiMNN10pn}
\vspace*{0mm}
\center\includegraphics[width=0.999\linewidth]{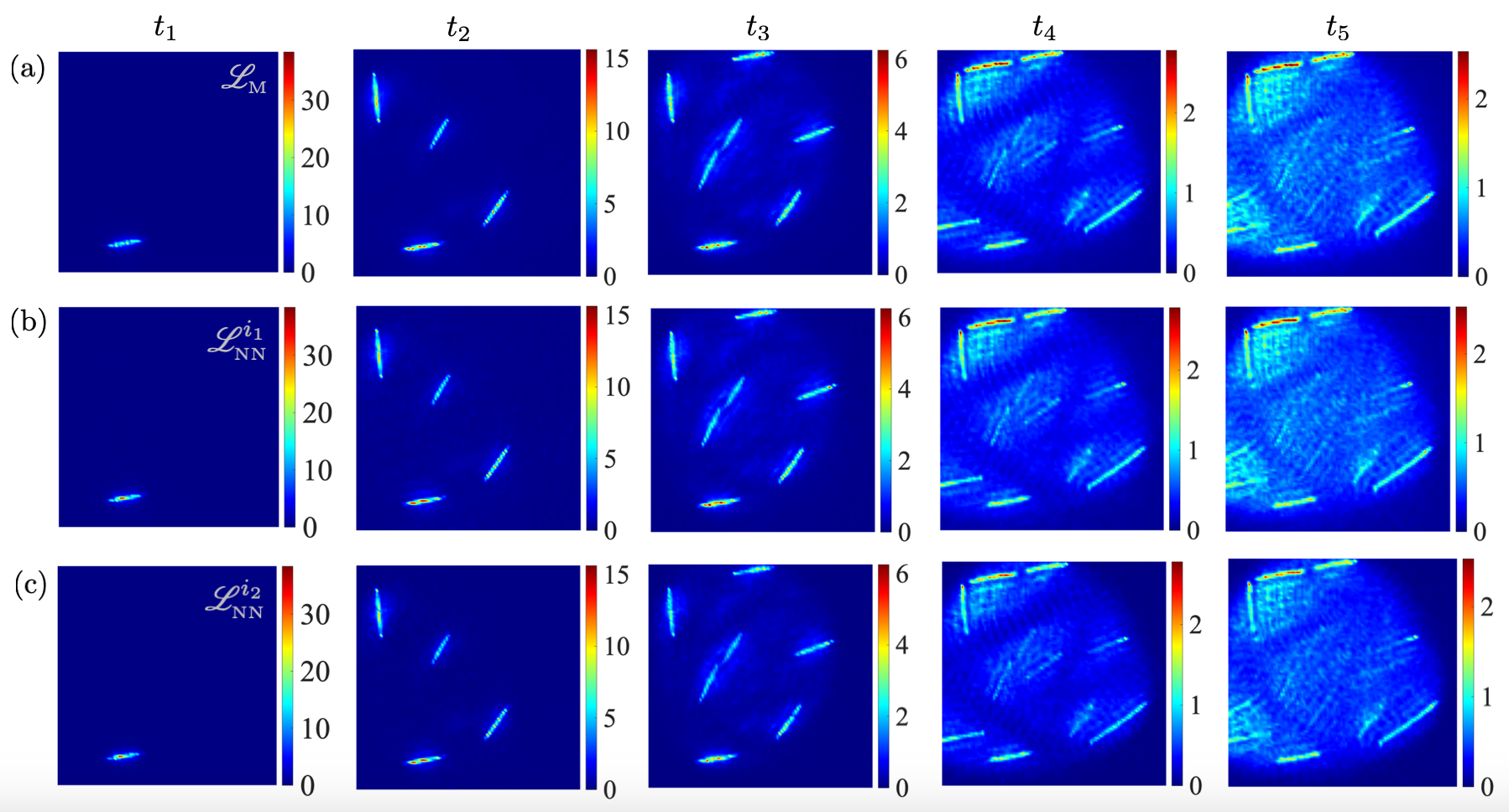} \vspace*{-7mm} 
\caption{LSM reconstructions corresponding to the regularization maps of Fig.~\ref{aiMNN10pn}:~(a) LSM images $\mathcal{L}_{\text{\tiny \emph{M}}}$ based on the manually optimized Morozov maps $\alpha_{\text{\tiny \emph{M}}}$,~(b) LSM reconstructions $\mathcal{L}_{\text{\tiny \emph{NN}}}^{i_1}$ by the network-generated regularization maps $\alpha_{\text{\tiny \emph{NN}}}^{i_1}$ at the end of Step 1 (at $\text{\emph{t}}=10^4$), and~(c) LSM indicator $\mathcal{L}_{\text{\tiny \emph{NN}}}^{i_2}$ via the R-Net regularization maps $\alpha_{\text{\tiny \emph{NN}}}^{i_2}$ at the end of training (when $s_t = 1$ in Step 2). Maps corresponding to each sensing step are plotted on the same color scale.} \lb{liMNN10pn}
\vspace*{-2mm}
\end{figure} 

\renewcommand{\arraystretch}{1.25}
 \begin{table}[!h]
\begin{center}
\caption{\small Contrast metric computed for the LSM reconstructions in Fig.~\ref{liMNN10pn} corresponding to \emph{discrepancy-informed R-Nets} trained on $10\%$ noisy data.} \vspace*{-1.5mm}
\label{T_liMNN10pn}
\begin{tabular}{c||c|c|c|c|c} 
\diagbox{{\small{contrast}}}{{\!\small{time}}} & $t_1$\! & $t_2$\! & $t_3$\! & $t_4$\! & $t_5$\! \\ \hline\hline  
$\mathfrak{C}_{\text{\tiny mn}}(\mathcal{L}_\text{\tiny M})$    & $94.23$  & $28.22$  & $10.61$ & $3.59$ & $2.48$ \\  
 \hline
$\mathfrak{C}_{\text{\tiny mn}}(\mathcal{L}_\text{\tiny NN}^{i_1})$    & \!\!$111.91$\!\!  & \!\!$26.18$\!\!  & \!\!$10.68$\!\! & \!\!$3.43$\!\! & \!\!$2.41$\!\! \\ 
 \hline
$\mathfrak{C}_{\text{\tiny mn}}(\mathcal{L}_\text{\tiny NN}^{i_2})$   & \!\!$137.02$\!\!  & \!\!$36.48$\!\! & \!\!$11.43$\!\! & \!\!$3.67$\!\! & \!\!$2.48$\!\! \\
 \hline
$\mathfrak{C}_{\text{\tiny mx}}(\mathcal{L}_\text{\tiny M})$    & $141.23$  & $48.31$  & $26.52$ & $9.85$ & $7.30$ \\ 
 \hline
$\mathfrak{C}_{\text{\tiny mx}}(\mathcal{L}_\text{\tiny NN}^{i_1})$     & \!\!$215.12$\!\!  & \!\!$50.00$\!\!  & \!\!$26.36$\!\! & \!\!$9.81$\!\! & \!\!$6.92$\!\! \\ 
 \hline
$\mathfrak{C}_{\text{\tiny mx}}(\mathcal{L}_\text{\tiny NN}^{i_2})$   & \!\!$284.63$\!\!  & \!\!$71.76$\!\! & \!\!$27.24$\!\! & \!\!$9.98$\!\! & \!\!$6.80$\!\! \\
\end{tabular}
\end{center}
\vspace*{-1mm}
\end{table}

\begin{figure}[!tp]
\center\includegraphics[width=0.999\linewidth]{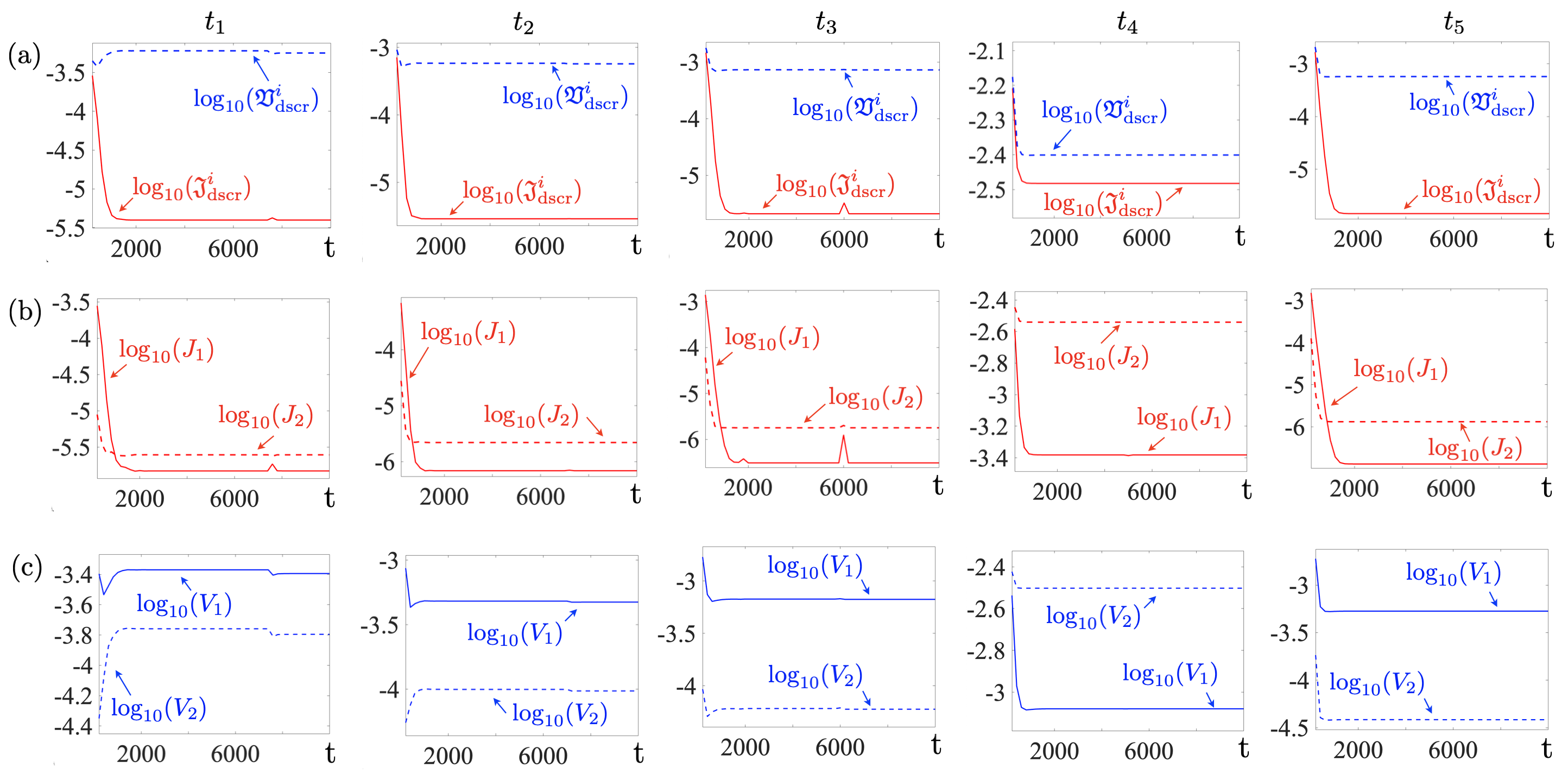} \vspace*{-7mm} 
\caption{Loss trajectories in Step~1 of training the \emph{discrepancy-informed R-Nets} on $25\%$ noisy data:~(a) total training loss $\log_{10}(\mathfrak{J}_{\text{dscr}}^i)$ (solid red line) and total validation loss $\log_{10}(\mathfrak{V}_{\text{dscr}}^i)$ (dashed blue line) against the number of epochs $\text{\emph{t}}$ at every sensing step $t_k$, $k = 1,2,\ldots, 5$, (b) training loss components $\log_{10}({J}_{1})$ (solid red line) and $\log_{10}({J}_{2})$ (dashed blue line) versus $\text{\emph{t}}$, and (c) validation loss components $\log_{10}({V}_{1})$ (solid red line) and $\log_{10}({V}_{2})$ (dashed blue line) versus $\text{\emph{t}}$.} \lb{lvi25pn}
\vspace*{-2mm}
\end{figure} 

 \begin{figure}[!tp]
\center\includegraphics[width=0.999\linewidth]{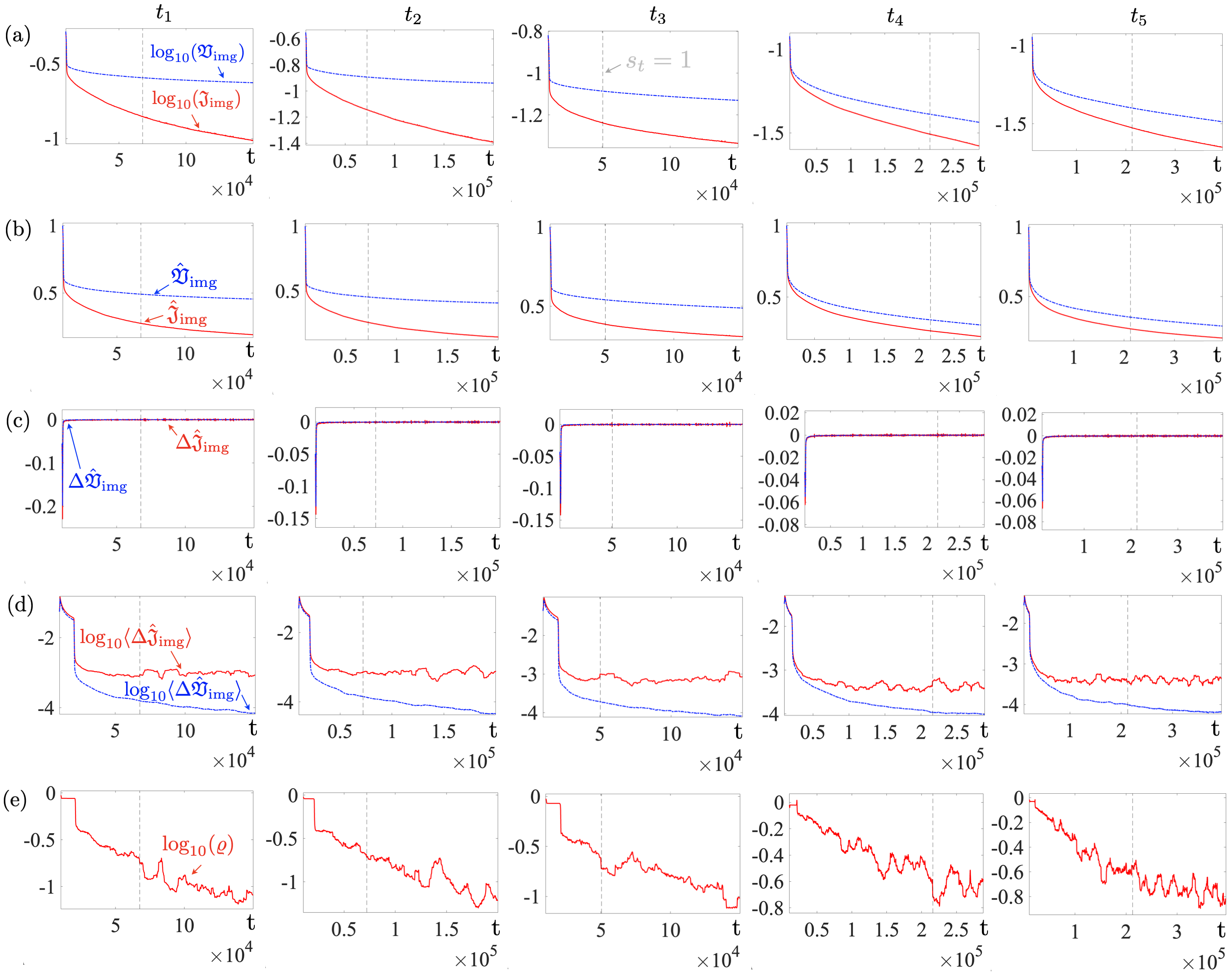} \vspace*{-7mm} 
\caption{Loss trajectories in Step~2 of training the \emph{discrepancy-informed R-Nets} on $25\%$ noisy data:~(a) training loss $\log_{10}(\mathfrak{J}_{\text{img}})$ (solid red line) and validation loss $\log_{10}(\mathfrak{V}_{\text{img}})$ (dash-dotted blue line) against the number of epochs $\text{\emph{t}}$ at every sensing step $t_k$, $k = 1,2,\ldots, 5$, (b) normal training loss $\hat{\mathfrak{J}}_{\text{img}\!}$ (solid red line) and normal validation loss $\hat{\mathfrak{V}}_{\text{img}\!}$ (dash-dotted blue line) versus $\text{\emph{t}}$, (c) variation of normal training loss $\Delta\hat{\mathfrak{J}}_{\text{img}\!}$ (solid red line) and variation of normal validation loss $\Delta\hat{\mathfrak{V}}_{\text{img}\!}$ (dash-dotted blue line) against $\text{\emph{t}}$, (d) rms of variations of normal training loss $\log_{10} \langle \Delta\hat{\mathfrak{J}}_{\text{img}} \rangle$ (solid red line) and rms of variations of normal validation loss  $\log_{10} \langle\Delta\hat{\mathfrak{V}}_{\text{img}}\rangle$ (dash-dotted blue line) versus the number of epochs $\text{\emph{t}}$, and (e) relative loss trajectory $\log_{10} ( \varrho ) (\text{\emph{t}})$. In all panels, the vertical dashed line indicates where the stop training criteria per Algorithm~\ref{AL2} is satisfied i.e., where $s_t = 1$. } \lb{iil25pn}
\vspace*{-2mm}
\end{figure} 

 \begin{figure}[!tp]
\center\includegraphics[width=0.999\linewidth]{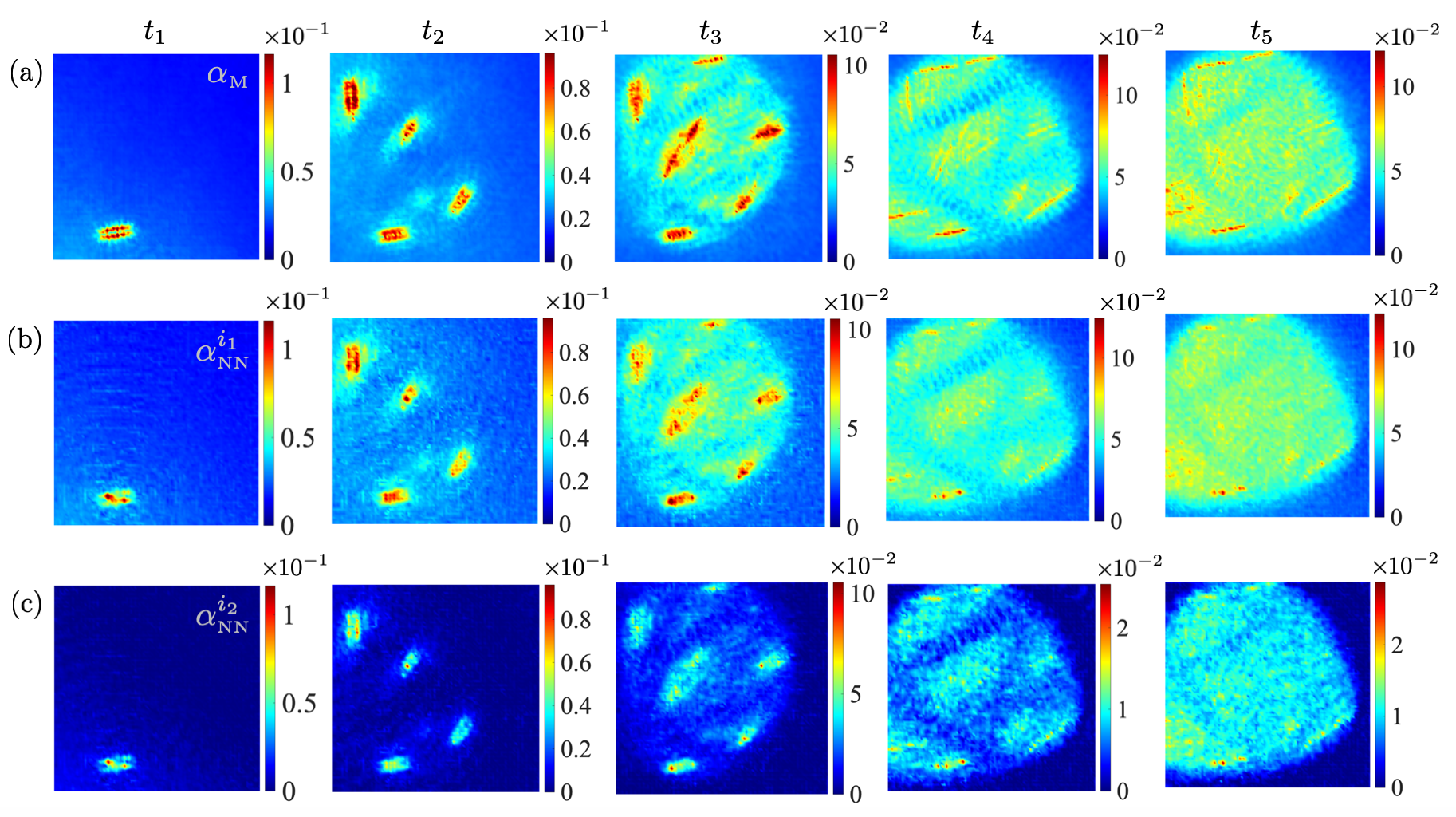} \vspace*{-7mm} 
\caption{Performance of the \emph{discrepancy-informed R-Nets} trained on $25\%$ noisy data:~(a) manually optimized Morozov regularization maps $\alpha_{\text{\tiny \emph{M}}}$ on the dense grid (of $100 \!\times\! 100$ sampling points) at every sensing step $t_k$, $k = 1,2,\ldots, 5$, (b) network-predicted regularization maps $\alpha_{\text{\tiny \emph{NN}}}^{i_1}$ by the end of Step 1 (at $\text{\emph{t}}=10^4$) on the same sampling grid, and (c) R-Net-generated regularization maps $\alpha_{\text{\tiny \emph{NN}}}^{i_2}$ by the end of Step 2 (when $s_t = 1$ as shown in Fig.~\ref{iil25pn}). Maps corresponding to each sensing step are plotted on the same color scale. } \lb{aiMNN25pn}
\vspace*{0mm}
\center\includegraphics[width=0.999\linewidth]{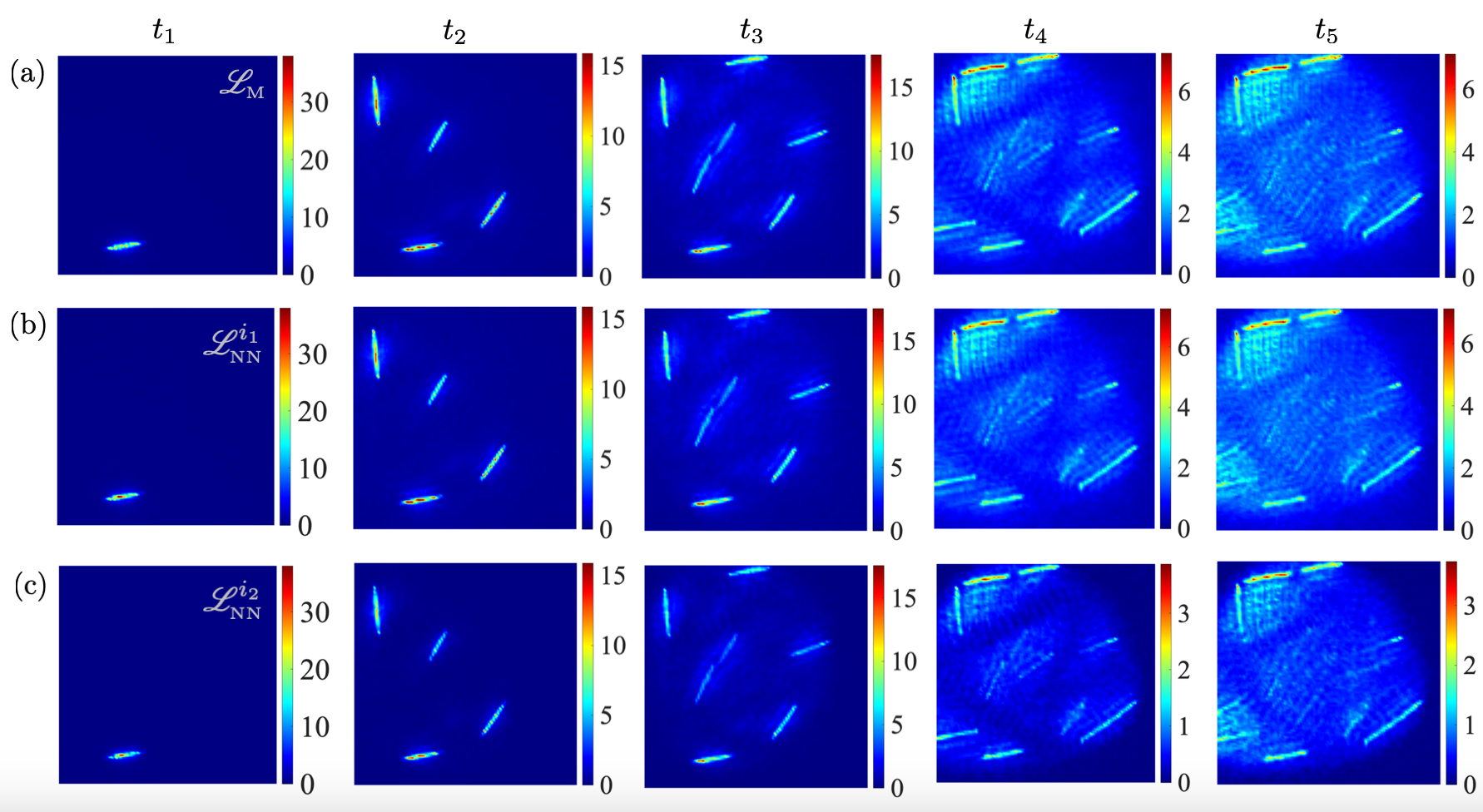} \vspace*{-7mm} 
\caption{LSM reconstructions corresponding to the regularization maps of Fig.~\ref{aiMNN25pn}:~(a) LSM images $\mathcal{L}_{\text{\tiny \emph{M}}}$ based on the manually optimized Morozov maps $\alpha_{\text{\tiny \emph{M}}}$,~(b) LSM reconstructions $\mathcal{L}_{\text{\tiny \emph{NN}}}^{i_1}$ by the network-generated regularization maps $\alpha_{\text{\tiny \emph{NN}}}^{i_1}$ at the end of Step 1 (at $\text{\emph{t}}=10^4$), and~(c) LSM indicator $\mathcal{L}_{\text{\tiny \emph{NN}}}^{i_2}$ via the R-Net regularization maps $\alpha_{\text{\tiny \emph{NN}}}^{i_2}$ at the end of training (when $s_t = 1$ in Step 2). Maps corresponding to each sensing step are plotted on the same color scale.} \lb{liMNN25pn}
\vspace*{-2mm}
\end{figure}

\renewcommand{\arraystretch}{1.25}
 \begin{table}[!h]
\begin{center}
\caption{\small Contrast metric computed for the reconstructions in Fig.~\ref{liMNN25pn} corresponding to discrepancy-informed R-Nets trained on $25\%$ noisy data.} \vspace*{-1.5mm}
\label{T_liMNN25pn}
\begin{tabular}{c||c|c|c|c|c} 
\diagbox{{\small{contrast}}}{{\!\small{time}}} & $t_1$\! & $t_2$\! & $t_3$\! & $t_4$\! & $t_5$\! \\ \hline\hline  
$\mathfrak{C}_{\text{\tiny mn}}(\mathcal{L}_\text{\tiny M})$    & $87.56$  & $24.45$  & $8.48$ & $2.97$ & $2.24$ \\  
 \hline
$\mathfrak{C}_{\text{\tiny mn}}(\mathcal{L}_\text{\tiny NN}^{i_1})$    & \!\!$98.57$\!\!  & \!\!$23.30$\!\!  & \!\!$7.91$\!\! & \!\!$2.74$\!\! & \!\!$2.12$\!\! \\ 
 \hline
$\mathfrak{C}_{\text{\tiny mn}}(\mathcal{L}_\text{\tiny NN}^{i_2})$   & \!\!$132.28$\!\!  & \!\!$38.88$\!\! & \!\!$10.74$\!\! & \!\!$3.81$\!\! & \!\!$2.57$\!\! \\
 \hline
$\mathfrak{C}_{\text{\tiny mx}}(\mathcal{L}_\text{\tiny M})$   & $118.05$  & $43.01$  & $19.26$ & $8.15$ & $6.35$ \\ 
 \hline
$\mathfrak{C}_{\text{\tiny mx}}(\mathcal{L}_\text{\tiny NN}^{i_1})$    & \!\!$173.37$\!\!  & \!\!$42.10$\!\!  & \!\!$19.09$\!\! & \!\!$6.65$\!\! & \!\!$5.53$\!\! \\ 
 \hline
$\mathfrak{C}_{\text{\tiny mx}}(\mathcal{L}_\text{\tiny NN}^{i_2})$   & \!\!$253.38$\!\!  & \!\!$74.64$\!\! & \!\!$29.32$\!\! & \!\!$10.99$\!\! & \!\!$7.38$\!\! \\
\end{tabular}
\end{center}
\vspace*{-2mm}
\end{table}

\begin{table}[!h]
\fontsize{10}{12}\selectfont \caption{Average training time per epoch for every step of training R-Nets.}
\label{computational_cost}
\centering{}%
\vspace{-2.0mm}
\begin{tabular}{cc}
\toprule
method & computational cost(time per epoch)\\*[0.5mm]
\midrule
{\small Step 1 (basic)}\!\!\!\! & 0.25s\\*[0.5mm]
{\small Step 1 (informed)} & 0.48s\\*[0.5mm]
{\small Step 2 (basic/informed)} & 0.72s \\*[0.5mm]
\bottomrule
\end{tabular}
\vspace*{-4.0mm}
\end{table}

\section{Conclusion}\lb{Conc}   


The \emph{discrepancy-informed R-Nets} are introduced to accelerate and enhance the solution of large inverse problems involving noisy operators. In this work, the focus is on \emph{online} supervised learning of Tiknonov regularization which is germane to imaging by the LSM indicator. The idea is to downscale the inverse problem and learn the regularization maps via R-Nets on a low-resolution dataset in a generalizable manner. The trained R-Nets will then solve the original high-dimensional problem and further optimize the inverse solution. As such, training R-Nets entails two steps:~(1) learning the logic of Morozov discrepancy principle to upscale the regularization maps, and (2) optimizing the R-Net predictions through further minimization of the Tikhonov loss within the Bayes risk minimization framework. Step 1 requires optimization of a many-objective loss function. For this purpose, by taking advantage of the GradNorm and DynScl logics, we proposed an adaptive loss balancing technique  that does not require a separate optimization. Step 2 demands careful regulation due to absence of labeled data. This is addressed by introducing a criteria to stop training based on the relative trajectories of training and validation loss functions. The proposed approach is synthetically benchmarked for ultrasonic imaging of progressive damage in an elastic plate using the LSM indicator. The results indicate that discrepancy-informed R-Nets can effectively upscale the regularization maps and remarkably enhance the inverse solution (i.e., image contrast), particularly, in reconstructions from noisy data in complex environments.

\section{Acknowledgements}\lb{Ac}   

The corresponding author kindly acknowledges the support provided by the National Science Foundation (Grant No.~1944812). This work utilized resources from the University of Colorado Boulder Research Computing Group, which is supported by the National Science Foundation (awards ACI-1532235 and ACI-1532236), the University of Colorado Boulder, and Colorado State University. The authors wish to thank Prof.~Luis Tenorio for his insightful comments during the course of this investigation.

\bibliographystyle{elsarticle-num}
\bibliography{inverse_with_doi}

\end{document}